\documentclass{article}
\usepackage[utf8]{inputenc}
\usepackage[T1]{fontenc}
\usepackage{amsmath}
\usepackage{amsfonts}
\usepackage{amsthm}
\usepackage{amssymb}
\usepackage{color}

\usepackage[margin=1in]{geometry}

\providecommand{\bysame}{\leavevmode\hbox to3em{\hrulefill}\thinspace}
\providecommand{\MR}{\relax\ifhmode\unskip\space\fi MR }

\providecommand{\href}[2]{#2}

\newcommand{\red}{\color{red}\tt}
\newcommand{\blue}{\color{blue}}

\newcommand{\bpf}{\begin{proof}}
\newcommand{\epf}{\end{proof} \medskip}

\newcommand{\benum}{\begin{enumerate}}
\newcommand{\eenum}{\end{enumerate}}

\newcommand{\bitem}{\begin{itemize}}
\newcommand{\eitem}{\end{itemize}}

\newcommand{\brmq}{\begin{rmq}}
\newcommand{\ermq}{\end{rmq}}

\newcommand{\brmqs}{\begin{rmqs}}
\newcommand{\ermqs}{\end{rmqs}}

\newcommand{\bapp}{\begin{application}}
\newcommand{\eapp}{\end{application}}

\newcommand{\bapps}{\begin{applications}}
\newcommand{\eapps}{\end{applications}}

\newcommand{\bdefi}{\begin{definition}}
\newcommand{\edefi}{\end{definition}}

\newcommand{\beq}{\begin{equation}}
\newcommand{\eeq}{\end{equation}}

\def\bpm{\begin{pmatrix}}
\def\epm{\end{pmatrix}}

\newcommand{\bcas}{\begin{cases}}
\newcommand{\ecas}{\end{cases}}

\newcommand{\bex}{\begin{exemp}}
\newcommand{\eex}{\end{exemp}}

\newcommand{\bexs}{\begin{exemps}}
\newcommand{\eexs}{\end{exemps}}

\newcommand{\bprop}{\begin{proposition}}
\newcommand{\eprop}{\end{proposition}}

\newcommand{\bthm}{\begin{theoreme}}
\newcommand{\ethm}{\end{theoreme}}

\newcommand{\bcor}{\begin{corollaire}}
\newcommand{\ecor}{\end{corollaire}}

\newcommand{\blem}{\begin{lemme}}
\newcommand{\elem}{\end{lemme}}

\newcommand{\beqna}{\begin{eqnarray}}
\newcommand{\eeqna}{\end{eqnarray}}

\newcommand{\beqnas}{\begin{eqnarray*}}
\newcommand{\eeqnas}{\end{eqnarray*}}

\newcommand{\Var}{{\rm Var}} 

\def\SL{ \rm{SL}}

\def\Id{{\rm{Id}}} 
\def\id{{\rm{id}}} 

\def\vol{{\rm{vol}}}

\def\diag{{\rm diag}}

\def\supp{{\rm supp}}
\def\htt{{\rm ht}}
\def\Cum{{\rm Cum}}

\def\cC{{\mathcal C}}
\def\cD{{\mathcal D}}
\def\cE{{\mathcal E }}
\def\cF{{\mathcal  F}}
\def\cG{{\mathcal  G}}

\def\cJ{{\mathcal  J}}

\def\cN{{\mathcal N }}

\def\cS{{\mathcal S }}
\def\cV{{\mathcal V}}

\def\cQ{{\mathcal Q}}

\def\bbE{{\mathbb{E}}}

\newcommand{\bbN}{{\mathbb {N}}}
\newcommand{\bbP}{{\mathbb P}}
\newcommand{\bbR}{{\mathbb {R}}}
\newcommand{\bbT}{{\mathbb {T}}} 

\newcommand{\bbZ}{\mathbb {Z}}

\def\un{{\mathbf{1}}}
\def\zero{{\mathbf{0}}}

\def\bfb{{\mathbf{b}}}

\def\bfj{{\mathbf{j}}}

\def\bfm{{\mathbf{m}}}

\def\bfp{{\mathbf{p}}}
\def\bfq{{\mathbf{q}}}

\def\bfu{{\mathbf{u}}}

\def\bfv{{\mathbf{v}}}

\def\bfx{{\mathbf{x}}}

\def\bfy{{\mathbf{y}}}

\def\bfz{{\mathbf{z}}}
\def\bf\Sigma{{\mathbf{\Sigma}}}

\def\bdalp{{\boldsymbol{\alpha}}}

\def\bf\lambda{{\mathbf{\lambda}}}

\newtheorem{theoreme}{Theorem}[section]
\newtheorem{lemme}[theoreme]{Lemma}
\newtheorem{definition}[theoreme]{Definition}
\newtheorem{proposition}[theoreme]{Proposition}
\newtheorem{corollaire}[theoreme]{Corollary}

\newenvironment{exemp}{\noindent{\bf Example. --- }}{\par}
\newenvironment{exemps}{\noindent{\bf Examples}\benum}{\eenum\par}
\newtheorem{rmq}[theoreme]{Remark}
\newtheorem{rmqs}[theoreme]{Remarks}

\newenvironment{application}{\noindent{\bf Application. --- }}{\par}
\newenvironment{applications}{\noindent{\bf Applications. ---
}\benum}{\eenum\par}

\theoremstyle{definition}

\title{The invariance principle for inhomogeneous Diophantine approximations}

\author{Songzi Li\thanks{School of Mathematics, Renmin University of China, 59, Zhongguancun Da Jie, Beijing, 100872, China ({\sf sli@ruc.edu.cn}). Supported by NSFC No.~11901569.}}
\date{\today}

\begin{document}
\maketitle

\begin{abstract}
We establish the central limit theorem and the invariance principle for the inhomogeneous Diophantine approximations. The proof employs the cumulant method, which was developed by Björklund and Gorodnik to prove the central limit theorem in the homogeneous setting. Our approach also relies on the effective mixing of expanding translates for high-order correlations on the affine lattice space, extending the previous result by Kim. 
\end{abstract}



\section{Introduction and main results}
In the 1920s, Khintchine\cite{K26} established the celebrated theorem illustrating the Borel-Cantelli property in Diophantine approximations. Since then, various universal limit theorems from probability theory have been demonstrated to hold under various conditions in Diophantine geometry. In this work, we extend the central limit theorem (CLT) and the invariance principle (IP) to the inhomogeneous, non-simultaneous setting. 

Given $u \in M_{m, n}(\bbT)$ and $\bfx \in \bbT^{m}$, define the following inhomogeneous linear forms on $\bbZ^{n}$
\beqna\label{lin.fm}
L^{(i)}_{(u, \bfx)}(q_{1}, \dots, q_{n}) = \sum^{n}_{j=1}u_{ij}q_{j} + x_{i}, 
\eeqna
for $i =1, \dots, m$. For positive constants $\{\theta_{i}, \ i =1, \dots, m\}$ and $\{\omega_{i}, \  i =1, \dots, m\}$ with $\sum^{m}_{i=1}\omega_{i} = n$, we consider the Diophantine inequalities given by
\beqna\label{ineq.di}
|p_{i} + L^{(i)}_{(u, \bfx)}(q_{1}, \dots, q_{n})| < \theta_{i}\|\bfq\|^{-\omega_{i}}, \ i=1, \dots, m, 
\eeqna
for $\bfp =  (p_{1}, \dots, p_{m}) \in \bbZ^{m}$ and $\bfq = (q_{1}, \dots, q_{n}) \in \bbZ^{n} \setminus \{\zero\}$.  Let $\|\cdot\|$ be a norm on $\bbR^{n}$. The special case $\omega_{1} = \dots = \omega_{m}$ corresponds to the simultaneous case. 

Classical results imply that there are infinitely many solutions $(\bfp, \bfq) \in  \bbZ^{m} \times (\bbZ^{n} \setminus  \{\zero\})$ for the Diophantine inequalities. To describe the distribution of the number of solutions, let $\Delta_{T}(u, \bfx)$ denote the counting function, defined as
\beqnas
\nonumber \Delta_{T}(u, \bfx) = \sharp \{(\bfp, \bfq) \in \bbZ^{m} \times (\bbZ^{n} \setminus  \{\zero\}),\  0 < \|\bfq\| < T, \ |p_{i} + L^{(i)}_{(u, \bfx)}(q_{1}, \dots, q_{n})| < \theta_{i}\|\bfq\|^{-\omega_{i}}, \ i=1, \dots, m\}. 
\eeqnas

For the homogeneous case $\bfx = \zero$, Schmidt~\cite{schm1} established an asymptotic formula for $\Delta_{T}(u, \zero)$ as $T \rightarrow \infty$, which can be viewed the law of large numbers. Leveque~\cite{Lev1}, \cite{Lev2}, Philipp~\cite{Phil}, Fuchs~\cite{fuchs} derived the CLT of $\Delta_{T}(u, \zero)$ for the one-dimensional case. An IP was also derived by Fuchs \cite{fuchs2} in this situation.  Dolgopyat-Fayad-Vinogradov~\cite{DFV} studied the CLT for the simultaneous case in high dimensional case, while Björklund and Gorodnik~\cite{BG19} proved the CLT for the non-simultaneous approximations. 

The inhomogeneous analogues of Khintchine's theorem were proven by Cassels~\cite{Cas1} and Sz\"usz~\cite{sz58}. The Schmidt theorem and related definitions have been developed under a much more general circumstances, see for example Beresnevich-Velani \cite{BV06}, Badziahin-Beresnevich-Velani \cite{BBV13}. We mention that Dolgopyat-Fayad-Vinogradov~\cite{DFV} proved the CLT for the inhomogeneous, simultaneous case. 

Inspired by the work of Björklund-Gorodnik~\cite{BG19}, we prove the CLT for the inhomogeneous, non-simultaneous case. Throughout this paper,  we assume that $m \geq 2$, $n \geq 1$.\bthm\label{clt.inDA}
Assume that $(u, \bfx)$ is uniformly distributed on $\bbT^{mn} \times \bbT^{m}$. We have
\beqnas
\frac{\Delta_{T}(u, \bfx) - \sigma^{2}_{m, n}\log T}{\sigma_{m, n}\sqrt{\log T}} \Longrightarrow \cN(0, 1) \ in  \ distribution, 
\eeqnas
as $T \rightarrow \infty$, where 
\beqna\label{var.aff.f}
\sigma_{m, n}^{2} = 2^{m}(\prod^{m}_{i=1}\theta_{i})n\vol_{n}.
\eeqna
and $\vol_{n}$ is the Euclidean volume of the unit ball with respect to the norm $\|\cdot\|$.
\ethm

Via Dani correspondence, the CLT for the counting function in Diophantine approximations reduces to a CLT for flows on the affine lattice space, or more precisely, to the CLT for Siegel transforms of piecewise smooth functions. In \cite{DFV}, the authors employed the martingale approach, which was initiated by Le Borgne~\cite{Borgne}, to prove the CLT for diagonal flows on the lattice space. However, the method fails to extend to the non-simultaneous case. We adopt the cumulant method developed in \cite{BG19}  for the homogeneous, non-simultaneous setting. The cumulant method relies on the estimates for high-order correlations of expanding flows on the homogeneous space, tracing back to the seminal work by Kleinbock-Margulis \cite{KM96}. Based on Kim's recent work \cite{Kim1}, which established an analogue mixing result on the affine lattice space, we derive the exponential mixing for high-order correlations, thus rendering the cumulant method applicable for the inhomogeneous case. 

Furthermore,  for $t \in [0, 1]$, we construct a piecewise linear function $X_{N, t}$ on $\bbT^{mn} \times \bbT^{m}$ that interpolates between $0$ and $\Delta_{N}$ (see Section~\ref{ip} for the precise definition of $X_{N, t}$).  Using the cumulant method, we can more easily verify the moment condition required for tightness of the distributions, which yields the following IP. 
\bthm\label{ip.inDA}
Assume that $v = (u, \bfx)$ is uniformly distributed on $\bbT^{mn} \times \bbT^{m}$. We have
\beqnas
\frac{X_{N, t}(u, \bfx) -  \sigma^{2}_{m, n}t\log N}{\sigma_{m, n}\sqrt{\log N}} \Longrightarrow W_{t} \ in  \ distribution, 
\eeqnas
as $N \rightarrow \infty$ in $\cC([0, 1])$, where $\{W_{t}, \ t \in [0, 1]\}$ is the standard Brownian motion on $\bbR$. 
\ethm
The paper is organized as follows. Section~\ref{sec.2} is devoted to the mixing property of the flows on the affine lattice space. As a consequence, we derive the exponential mixing for high-order correlations, which implies the CLT for the smooth functions. In Section~\ref{sec.3} we extend the CLT to smooth Siegel transforms. We then prove the main result, Theorem~\ref{clt.inDA} in Section~\ref{sec.4}. In the last Section, we present the proof of Theorem~\ref{ip.inDA}.

\section{Mixing property of higher order correlations}\label{sec.2}
Let $G = \SL_{d}(\bbR)$, $\Gamma = \SL_{d}(\bbZ)$,  and let $X =  G / \Gamma $ be the uni-modular lattice space. Let $\hat{G} = \SL_{d}(\bbR) \ltimes \bbR^{d}$, $\hat{\Gamma} = \SL_{d}(\bbZ) \ltimes \bbZ^{d}$,  and $Y =\hat{G}/ \hat{\Gamma}$ be the affine lattice space. The multiplication law on $\hat{G}$ is given by
\beqnas
(g_{1}, \bfv_{1})(g_{2}, \bfv_{2}) = (g_{1}g_{2}, \bfv_{1}+ g_{1}\bfv_{2}). 
\eeqnas
The action of $(g, \bfv) \in \hat{G}$ on $\bbR^{d}$ is defined as
\beqnas
(g, \bfv)\bfz = g\bfz + \bfv.
\eeqnas
Denote by $d_{G}$ and $d_{\hat{G}}$ the right invariant Riemannian metric on $G$ and $\hat{G}$, which induces metrics $d_{X}$ and $d_{Y}$ on $X$ and $Y$, respectively. Define the norm on $G$ as $\|g\| = \max_{1 \leq i, j \leq d}\{|g_{ij}|, |g^{-1}_{ij}|\}$. Denote $m_{X}$, $m_{Y}$ by the normalized Haar measure on $X$, $Y$, respectively.  Let 
\beqna\label{def.h}
H = \{M | M = \begin{pmatrix}
\Id_{m}  & u\\
0 & \Id_{n}
\end{pmatrix}, 
u \in M_{m, n}\}
\eeqna
and $m_{H}$ be the Haar measure on $H$. Notice that H is a unipotent abelian subgroup of $G$ which is expanding horospherical with respect to 
\beqna\label{def.aat}
a_{t} = \diag \{e^{nt}\un_{m}, e^{-mt}\un_{n}\}
\eeqna
for $t > 0$. For a subset $V \subset H$ with compact closure, define $m_{V} = \frac{1}{m_{H}(V)}m_{H}$.

The effective ergodicity theorem proved by Kleinbock-Margulis~\cite{KM96} shows the exponential mixing property of the expanding translates of $H$ on $X$. 

\bthm[Kleinbock-Margulis]\label{KM}
Let $V \subset H$ be a fixed neighborhood of the identity in $H$ with smooth boundary and compact closure. For any compact set $L \subset X$ and $x_{0} \in L$, there exist constants $\lambda= \lambda(m, n) > 0$ and $T(L) \geq 0$ such that for any $f \in \cC^{\infty}_{c}(X)$ and $t \geq T(L)$,
\beqna\label{effe}
\int_{V}f(a_{t}ux_{0})dm_{V}(u)  - \int_{X}fdm_{X} = O(\cS(f)e^{-\lambda t}),
\eeqna
where $\cS$ is a Sobolev norm on $\cC^{\infty}_{c}(X)$ and the constant depends on $m$, $n$, $V$.
\ethm

Kim~\cite{Kim1} recently extended Kleinbock-Margulis' result to the space of affine lattices. He introduced the function $\zeta: \bbR^{d} \times \bbR^{+} \rightarrow \bbN$, given by 
\beqna\label{def.zeta}
\zeta(\bfb, T) = \min \{N \in \bbN, \min_{1 \leq |q| \leq N} |q\bfb|_{\bbZ} \leq \frac{N^{2}}{T}\},
\eeqna
where  $|\cdot|_{\bbZ}$ denotes the supremum distance from $0 \in \bbT^{d}$. 

\bthm[Kim]\label{thm.kim}
Let $V \subset H$ be a fixed neighborhood of the identity in $H$ with smooth boundary and compact closure. There exists a constant $\delta = \delta(d) > 0$ such that
\beqna\label{effe.eqdis.aff.1}
\int_{V}f(a_{t}uy_{0}) dm_{V}(u) - \int_{Y}f dm_{Y} = O(\cS(f)\zeta(\bfb_{0}, e^{\frac{nt}{2}})^{-\delta}),
\eeqna
holds for $f \in \cC^{\infty}_{c}(Y)$,  where $\cS$  is a Sobolev norm on $\cC^{\infty}_{c}(Y)$ , $y_{0} = (g_{0}, \zero)(1, \bfb_{0})\hat{\Gamma}\in Y$ with $g_{0} \in \SL(d, \bbR)$, $\bfb_{0} \in \bbT^{d}$ and $t \geq 0$ such that  $\|g_{0}\| \leq \zeta(\bfb_{0}, e^{\frac{nt}{2}})^{\delta}$. The constant depends on $V$ and $d$. 
\ethm
We consider
\beq\label{def.gt}
g_{t} = \diag \{e^{\alpha_{1}t}, \dots, e^{\alpha_{m}t}, e^{-\alpha_{m+1}t}, \dots e^{-\alpha_{m+n}t}\}
\eeq
for $t > 0$ and $\bar{\alpha} = (\alpha_{1}, \dots, \alpha_{m+n})$ with positive numbers $\alpha_{1}, \dots, \alpha_{m+n}$ satisfying
\beqnas
\sum^{m}_{i=1}\alpha_{i} = \sum^{m+n}_{i=m+1}\alpha_{i}.
\eeqnas
Define
\beqnas
\alpha_{0} = \min_{1 \leq i \leq d}\alpha_{i}.
\eeqnas

The effective equidistribution results with respect to $g_{t}$ also hold on $\SL(d, \bbR) / \SL(d, \bbZ)$ and $\SL_{d}(\bbR) \ltimes \bbR^{d} / \SL_{d}(\bbZ) \ltimes \bbZ^{d}$, proved by \cite{KM12} and \cite{Kim1} respectively. 

The mixing property of high-order correlations on $X = \SL(d, \bbR) / \SL(d, \bbZ)$ was obtained in Theorem 2.2 (Corollary 2.4), \cite{BG19}. Let $k \geq 1$. For a $k-$tuple $(t_{1}, \dots, t_{k})$, which we may assume $t_{1} \leq \dots \leq t_{k}$, define
 $$
D(t_{1}, \dots, t_{k}) = \min\{t_{k}, |t_{i} - t_{j}|, 1 \leq i \neq j \leq k\}.
$$
\bthm[Bj\"orklund-Gorodnik]\label{thm.bg}
Let $V \subset H$ be a fixed neighborhood of the identity in $H$ with smooth boundary and compact closure. There exists $\lambda' >0$ such that for $f \in \cC_{c}^{\infty}(V)$, $f_{1}, \dots, f_{k} \in \cC^{\infty}_{c}(X)$, $x_{0} \in X$, and $t_{1}, \dots, t_{k} > 0$, we have
\beqna\label{exp.m.X}
 \int_{V}f(u)\big(\prod^{k}_{i=1}f_{i}(g_{t_{i}}ux_{0})\big)dm_{V}(u) -  \int_{V}f dm_{V} \prod^{k}_{i=1}(\int_{X}f_{i} dm_{X}) = O(\|f\|_{\cC^{k}}\prod^{k}_{i=1}\cS(f_{i})e^{-\lambda' D(t_{1}, \dots, t_{k})}),
\eeqna
where the constant depends on $x_{0}$, $k$ and $V$. 
\ethm
We prove the following effective equidistribution for high-order correlations on $Y = \SL_{d}(\bbR) \ltimes \bbR^{d} / \SL_{d}(\bbZ) \ltimes \bbZ^{d}$.
\bthm\label{thm.2}
Let $V \subset H$ be a fixed neighborhood of the identity in $H$ with smooth boundary and compact closure. There exists $\delta' > 0$ independent of $k$,  such that for $f_{1}, \dots, f_{k} \in \cC^{\infty}_{c}(Y)$ and $t_{1}, \dots, t_{k} > 0$, 
\beqna\label{effe.eqdis.2}
\int_{V}\prod^{k}_{i =1}f_{i}(g_{t_{i}}uy_{0}) dm_{V}(u)  - \prod^{k}_{i =1}\int_{Y}f_{i} dm_{Y} = O(\prod^{k}_{i =1}\cS(f_{i})\zeta(\bfb_{0}, e^{\frac{\alpha_{0}D(t_{1}, \dots, t_{k})}{2}})^{-\delta'}),
\eeqna
holds for $y_{0} =  (g_{0}, \zero)(1, \bfb_{0})\hat{\Gamma} \in Y$ with $\|g_{0}\| \leq \zeta(\bfb_{0}, e^{\frac{\alpha_{0} D(t_{1}, \dots, t_{k})}{2}})^{\delta'}$. The implied constant depends on $k$, $V$ and $d$.
\ethm
As an application, we deduce the following exponential mixing result from Theorem~\ref{thm.2}. 
\bcor\label{thm.3}
Let $V \subset H$ be the same as in Theorem~\ref{thm.2}. Let $m_{\bbT^{d}}$ be the Lebesgue measure on $\bbT^{d}$. Fix $\kappa > d+1$. Denote $y(\bfb) =  (g_{0}, \zero)(1, \bfb)\hat{\Gamma} \in Y$. There exists $\delta_{\kappa} > 0$ such that for $f_{1}, \dots, f_{k} \in \cC^{\infty}_{c}(Y)$ and $t_{1}, \dots, t_{k} > 0$, we have
\beqna\label{effe.eqdis.3}
\int_{\bbT^{d}}\int_{V}\prod^{k}_{i =1}f_{i}(g_{t_{i}}uy(\bfb)) dm_{V}(u)dm_{\bbT^{d}}(\bfb) - \prod^{k}_{i =1}\int_{Y}f_{i} dm_{Y} = O(\prod^{k}_{i =1}\cS(f_{i})e^{-\alpha_{0} D(t_{1}, \dots, t_{k})\delta_{\kappa}}),
\eeqna
holds for $\|g_{0}\| \leq e^{\alpha_{0} D(t_{1}, \dots, t_{k})\delta_{\kappa}}$. The implied constant depends on $k$, $V$ and $d$.
\ecor

The exponential mixing property \eqref{effe.eqdis.3} allows us to apply the cumulant method of \cite{BG19} to prove the central limit theorem below. As the argument closely parallels that of \cite{BG19}, we omit the proof here.
\bthm\label{clt.1}
Let $V \subset H$ be the same as in Theorem~\ref{thm.2}. Let $y(\bfb) =  (g_{0}, \zero)(1, \bfb)\hat{\Gamma} \in Y$ for some $g_{0} \in \SL(d, \bbR)$. Assume that $(u, \bfb)$ is uniformly distributed on $V \times \bbT^{d}$. 
For $f \in \cC^{\infty}_{c}(Y)$, we have
\beqnas
\frac{1}{\sqrt{N}}\sum^{N-1}_{s=0}\big( f(g_{s}uy(\bfb)) - \int_{\bbT^{d}}\int_{V}f(g_{s}uy(\bfb))dm_{V}(u)dm_{\bbT^{d}}(\bfb)\big) \Longrightarrow \cN(0, \sigma^{2}), 
\eeqnas
as $N \rightarrow \infty$, and 
\beqnas
\sigma^{2} = \sum^{+\infty}_{s = -\infty}\big(\int_{Y}(f \cdot g_{s})f dm_{Y} - (\int_{Y}f dm_{Y})^{2} \big).
\eeqnas
\ethm

\subsection{Preliminary facts}

We present some basic facts needed for the proof of Theorem \ref{thm.2} and Corollary \ref{thm.3}. Since these results extend Theorem \ref{thm.kim}, we mainly follow \cite{Kim1} for the preliminaries. The notation $A \ll B$ indicates that $|A| \leq c |B|$ for some constant $c > 0$, where $c$ only depends on the dimensions $m$, $n$, $d$ and the domain $V$. The relation $A \asymp B$ denotes that both $A \ll B$ and $B \ll A$ hold.

\subsubsection{The Sobolev norms}
We recall the Sobolev norm introduced in \cite{Kim1}. For every $V \in \cG$, define the differential operator $D_{V}$ on $\cC^{\infty}_{c}(X)$ by $D_{V}\phi(x) = \frac{d}{dt}|_{t=0}\phi(e^{tV}x)$. For a basis $\{V_{1}, \dots, V_{d^{2}-1}\}$ of $\cG$, every monomial $Z = V^{l_{1}}_{1}\dots V^{l_{r}}_{r}$ defines a differential operator by $D_{Z} = D^{l_{1}}_{V_{1}}\dots D^{l_{r}}_{V_{r}}$ with degree $\deg(Z) = l_{1} + \dots l_{r}$. For $k \in \bbN$, $f \in \cC^{\infty}_{c}(X)$, define the norm $\cS^{X}_{l}$ by
\beqnas\label{def.norm.X}
\cS^{X}_{l}(f)^{2} &=& \sum_{\deg(Z) \leq l}\int_{X}|\htt(x)^{l}D_{Z}f(x)|^{2}dm_{X}(x).
\eeqnas
Similarly, define the Sobolev norm $\cS^{Y}_{l}$ on $\cC^{\infty}_{c}(Y)$ for $l \in \bbN$, $f \in \cC^{\infty}_{c}(Y)$ by
\beqnas\label{def.norm.Y}
\cS^{Y}_{l}(f)^{2} = \sum_{\deg(\hat{Z}) \leq l}\int_{Y}|\htt(\pi(y))^{l}\cD_{\hat{Z}}f(y)|^{2}dm_{Y}(y),
\eeqnas
where $\pi: Y \rightarrow X$ is the natural projection, and $\hat{Z}$ is the monomial generated by the basis of the Lie algebra $\hat{\cG}$. 

We mention the following properties of $\cS^{Y}$ due to \cite{Kim1}: for $f \in \cC^{\infty}_{c}(Y)$, $l$ large enough and $\deg(\hat{Z}) \leq d+2$
\beqna\label{norm.1}
\|\cD_{\hat{Z}}f\|_{L^{\infty}(Y)} \leq \cS^{Y}_{l}(f).
\eeqna
For $g \in G$ and  $f \in \cC^{\infty}_{c}(Y)$, consider $g.f(y) = f((g, \mathbf0)y)$ as a function on $Y$. Then, we can deduce from the properties of $\cS^{X}$ in \cite{EMV} that
\beqna\label{norm.2}
\cS^{Y}_{l}(g.f) \ll \|g\|^{c l}\cS^{Y}_{l}(f),
\eeqna
for some constant $c$, and 
\beqna\label{norm.3}
 \|f - g.f\|_{L^{\infty}} \ll d_{G}(\id, g)\cS^{Y}_{l}(f)
 \eeqna
for $l$ large enough. In this paper, we choose $l_{0}$ such that Theorem~\ref{thm.kim} and Theorem~\ref{thm.bg} hold with $\cS = \cS^{Y}_{l_{0}}$ and $\cS = \cS^{X}_{l_{0}}$ respectively, and such that  the above properties are satisfied by $\cS^{Y}_{l_{0}}$.

\subsubsection{Diophantine type vectors}
A vector $\bfb \in \bbR^{d}$ is said to be of Diophantine type $\kappa \geq 1$ if there exists $c_{\kappa} > 0$ such that
\beqna\label{da.k}
|\bfb - \frac{\bfp}{q}| > c_{\kappa}q^{-\kappa}
\eeqna
for any $\bfp \in \bbZ^{d}$ and $q \in \bbN$, where $|\cdot|$ is the supremum norm of $\bbR^{d}$. 

For $c>0$, define
\beqnas
D(\kappa, c) = \{\bfb \in \bbT^{d} , |q\bfb|_{\bbZ} >  cq^{-\kappa+1}, \ \forall q \in \bbN\}.
\eeqnas
Then for $c_{1} < c_{2}$, $D(\kappa, c_{2}) \subset D(\kappa, c_{1})$. Moreover, we have the following estimate.
\blem
For $\kappa > d+1$, we have
\beqna\label{est.3.1}
m_{\bbT^{d}}(\bbT^{d} \setminus D(\kappa, c) ) \ll c\sum^{\infty}_{q = 1}q^{-\kappa+d}.
\eeqna
\elem

\bpf
For $q \in \bbN$, define
$$
\Pi_{q} = \{\bfb \in \bbT^{d}, |q \bfb|_{\bbZ} \leq  cq^{-\kappa+1}\} = \{\bfb \in \bbT^{d}, |q \bfb - \bfp| \leq  cq^{-\kappa+1},\ for \ some \ \bfp \in \bbZ^{d}\} .
$$
Then 
$$
m_{\bbT^{d}}(\bbT^{d} \setminus D(\kappa, c) ) \leq \sum_{q \in \bbN \setminus \{\zero\}}m_{\bbT^{d}}(\Pi_{q}).
$$
Notice that
\beqnas
m_{\bbT^{d}}(\Pi_{q}) \ll c\sum_{\bfp = O(|q|)}q^{-\kappa} \ll cq^{-\kappa+d} ,
\eeqnas
which leads to \eqref{est.3.1}.
\epf

For $\epsilon > 0$, the set $\bigcup^{\infty}_{c \geq \epsilon}D(\kappa, c)$ consists of vectors of Diophantine type $\kappa$, and 
\beqna\label{est.m}
m_{\bbT^{d}}\big(\bbT^{d} \setminus \bigcup^{\infty}_{c \geq \epsilon} D(\kappa, c) \big) = m_{\bbT^{d}} \big(\bbT^{d} \setminus D(\kappa, \epsilon) \big) \ll \epsilon.
\eeqna

We also recall some facts on $\zeta(\bfb, T)$ in \cite{Kim1}. For a vector $\bfb$ of Diophantine type $\kappa$,  it holds
 \beqna\label{est.z}
\zeta(\bfb, T) \geq c^{\frac{1}{\kappa+1}}_{\kappa} T^{\frac{1}{\kappa+1}}.
\eeqna
Moreover, $\zeta(\bfb, \cdot)$ is non-decreasing, unbounded and
\beq
\label{z.p.1} \zeta(\bfb, cT) \leq \sqrt{c} \zeta(\bfb, T),
\eeq
\beq
\label{z.p.2} \zeta(\bfb, \|\alpha^{-1}\|^{-1}_{op}T) \leq \zeta(\alpha\bfb, T) \leq \zeta(\bfb, \|\alpha\|_{op}T)
\eeq
for $c>0$, $\alpha \in \Gamma$, 
and 
\beq
\label{z.p.3}  \zeta(\bfb, T) \leq T^{\frac{d}{2d+1}}. 
\eeq

\subsubsection{Some facts on $X = \SL(d, \bbR) / \SL(d, \bbZ)$} 
For $x \in X$, set
\beqnas
\htt(x) = \sup \{| g \bfv|^{-1}, x = g\Gamma , \bfv \in \bbZ^{d} \setminus \{0\}\},
\eeqnas
where $|\cdot|$ is the supremum norm of the vector. Note that there exists some constant $c_{1} > 1$ so that
\beqna\label{est.ht.1}
\htt(gx) \leq  c_{1}\|g\|\htt(x)
\eeqna
for any $x \in X$ and $g \in G$.
Let 
$$
K(R) = \{x \in X, \htt(x) \leq R\},
$$
then for all $R > 0$, $K(R)$ is compact due to Mahler's compact criterion. Moreover, 
\beqna\label{est.ht.2}
m_{X}(X \setminus K(R)) \asymp R^{-d}.
\eeqna
The estimate of injective radius in \cite{KM12} indicates that there exists a constant $c_{2} > 0$ such that for $x \in K(c_{2}r^{-\frac{1}{d}})$, $0 < r < \frac{1}{2}$, the map $g \mapsto gx$ is injective on $B^{G}(\id, r)$. 

In the following we recall the fundamental domain $\cF \subset G$ constructed in \cite{Kim1}.  For any $x \in X$, there exists a unique $g \in \cF$ satisfying $x = g\Gamma$. Let $\pi_{X}: G \rightarrow X$ be the canonical projection. Then one can define $\iota: X \rightarrow \cF$ such that $\pi_{X} \cdot \iota = \Id_{X}$. Note that $\iota$ is continuous on $\phi(\cF^{\circ})$ and measure preserving. Moreover, 
\beqna\label{ineq.iota}
\|\iota(x)\| \ll \htt(x)^{d-1}
\eeqna
 for any $x \in X$. For $x \in K(c_{2}r^{-\frac{1}{d}})$, $0 < r < \frac{1}{2}$, $\iota$ is an isometry on $B^{X}(x, r)$.

For $y \in Y$, define $\pi: Y \rightarrow X$ as the natural projection such that  $\pi(y) = x$. Then there is a unique decomposition: $y = (\iota(\pi(y)), \zero)(1, \bfb)\hat{\Gamma} \in Y$, where $\bfb \in \bbT^{d}$. Define $\sigma: Y \rightarrow \bbT^{d}$,  $\sigma(y) = \bfb$ as the projection to $\bbT^{d}$.

For $r > 0$ and $\epsilon > 0$, the author \cite{Kim1} introduced
\beqna\label{set.f}
\cF(r, \epsilon) = \{g \in \cF, \htt(g \Gamma) \leq \epsilon^{-1}, d_{G}(g, \partial \cF) \geq r, d_{G}(g, \partial \cE^{-1}) \geq r^{\frac{1}{20d}}\},
\eeqna
where $\partial \cF$ denotes the boundary of $\cF$ and $\cE^{-1} = \{g \in G, g^{-1} \in \cE\}$. One has
\beqna\label{est.f}
m_{G}(\cF \setminus \cF(r, \epsilon) ) \ll \max \{r^{c_{3}}, \epsilon^{d}\},
\eeqna
where $c_{3} = \frac{1}{100 d^{3}}$. 

We recall the partition of $X$ introduced in \cite{Kim1}. For $0 < r < \frac{1}{2}$, define $B^{G}_{r} = B^{H}(\id, r)B^{H^{0}}(\id, r)B^{H^{-}}(\id, r)$ and  $B_{r}(x) := B^{G}_{r}x$ for $x \in X$.
\bprop\label{def.ptt}
There exist constants $C_{1}, C_{3}>1$, $0 < C_{2} < 1$ such that the following holds. For $0 < r < \frac{1}{2C^{3}_{1}}$, there exist a set $\{z_{1}, \dots, z_{\cN_{r}}\} \subset K(C_{2}r^{-\frac{1}{d}})$ with $\cN_{r}  \asymp r^{-(d^{2}-1)}$ and a partition $\{\omega_{j}\}_{j \in \cJ}$ with $\cJ = \{1, \dots, \cN_{r}\} \cup \{\infty\}$, satisfying 
\beqnas
\sum_{j \in \cJ} \omega_{j} &=& \un_{X},\\
0 &\leq& \omega_{j} \ \leq \ 1, \ for \ j \in \cJ, \\
\un_{B_{r}(z_{j})} &\leq& \omega_{j} \ \leq \ \un_{B_{C^{3}_{1}r}(z_{j})}, \ for \ j \in \cJ \setminus \{\infty\},\\ 
\supp \  \omega_{\infty} &\subset& X \setminus K(\frac{1}{2}C_{2}r^{-\frac{1}{d}}), \\
\|\nabla \omega_{j}\|_{L^{\infty}(X)} &\leq& C_{3}r^{-1}, for \ j \in \cJ \setminus \{\infty\}.
\eeqnas
\eprop

It is pointed out in \cite{Kim1} that  one can apply the effective equidistribution~\eqref{effe} on $X$ to estimate $\pi_{*}\mu_{y_{0}, t}(B_{r}(x))$ for $y_{0} \in Y$. More precisely, there exists $0 < \kappa_{1} \leq \frac{1}{2}$ such that for $\htt(\pi(y_{0})) \leq e^{\kappa_{1}t}$, $e^{-\kappa_{1}t} < r < \frac{1}{2}$, and $x \in K(c_{2}r^{-\frac{1}{d}})$, 
\beqna\label{est.v.1}
\pi_{*}\mu_{y_{0}, t}(B_{r}(x)) \asymp m_{X}(B_{r}(x)) \asymp r^{d^{2}-1}. 
\eeqna 
Since $\pi_{*}\mu_{y_{0}, t}(B_{r}(z_{j})) \leq \pi_{*}\mu_{y_{0}, t}(\omega_{j}) \leq \pi_{*}\mu_{y_{0}, t}(B_{C^{3}_{1}r}(z_{j}))$ for $j \in \cJ \setminus \{\infty\}$, and $z_{j} \in K(C_{2}r^{-\frac{1}{d}})$ due to Proposition~\ref{def.ptt}, one also obtains
\beqna\label{est.v.o.1}
\pi_{*}\mu_{y_{0}, t}(\omega_{j}) \asymp m_{X}(B_{r}(z_{j})) \asymp r^{d^{2}-1},
\eeqna
under the same conditions on $y_{0}$ and $t$. 

The above estimates rely on the observation, as pointed out in \cite{Kim1}, that the error term in~\eqref{effe} depends on $\htt(x_{0})^{\kappa}$ with some $\kappa > 0$. We may also check that the dependence on $x_{0}$ in the error term in \eqref{exp.m.X} can be expressed explicitly as
\beqnas
O(\htt(x_{0})^{\kappa}e^{-\lambda' D(t_{1}, \dots, t_{k})}\|f\|_{\cC^{k}}\prod^{k}_{i=1}\cS(f_{i})),
\eeqnas 
which allows us to estimate the measure of $\prod^{k}_{i=1}B_{r_{i}}(x_{i})$ with respect to $\pi_{*}\mu_{y_{0}, t_{1}, \dots, t_{k}}$. More precisely, there exists $\kappa_{2} > 0$ small enough such that for $\htt(\pi(y_{0})) \leq e^{\kappa_{2}D(t_{1}, \dots, t_{k})}$, $e^{-\kappa_{2}D(t_{1}, \dots, t_{k})} < r < \frac{1}{2}$, and each $x_{i} \in K(c_{2}r^{-\frac{1}{d}})$,  we have
\beqna\label{est.v.k}
\pi_{*}\mu_{y_{0}, t_{1}, \dots, t_{k}}(\prod^{k}_{i=1}B_{r_{i}}(x_{i})) \asymp \prod^{k}_{i=1}m_{X}(B_{r_{i}}(x_{i})) \asymp \prod^{k}_{i=1}r^{d^{2}-1}_{i}.
\eeqna
Similarly, we obtain the estimate of $\{\omega_{j_{i}}\}_{\bfj_{k}}$, where $\bfj_{k} = (j_{1}, \dots, j_{k})$ and each $j_{i} \in \cJ \setminus \{\infty\}$, under the conditions that $\htt(\pi(y_{0})) \leq e^{\kappa_{2}D(t_{1}, \dots, t_{k})}$, $e^{-\kappa_{2}D(t_{1}, \dots, t_{k})} < r < \frac{1}{2}$,
\beqna
\label{est.v.o.k.1}\pi_{*}\mu_{y_{0}, t_{1}, \dots, t_{k}}(\prod^{k}_{i=1}\omega_{j_{i}}) & \asymp& r^{k(d^{2}-1)}.
\eeqna

We recall the estimates in Lemma 5.6 and Proposition 5.7 in \cite{Kim1}, which will be applied in our proof. Let $\cS^{+} = \{\bar{s} = (s_{1}, \dots, s_{d}), s_{1}, \dots, s_{d} > 0, \sum^{m}_{i=1}s_{i} = \sum^{d}_{i=m+1}s_{i}\}$. For $\bar{s} \in \cS^{+}$, set
\beqnas
\lfloor \bar{s} \rfloor = \min_{1 \leq i \leq d}s_{i}.
\eeqnas
Define
\beqnas
g_{\bar{s}} = \diag \{e^{s_{1}}, \dots, e^{s_{m}}, e^{-s_{m+1}}, \dots e^{-s_{d}}\},
\eeqnas
and $\xi_{\bar{s}}: X \times V \rightarrow \cF$ and  $\gamma_{\bar{s}}: X \times V \rightarrow \Gamma$ such that for any $x \in X$ and $u \in V$, there exists unique $\xi_{\bar{s}}(x, u) \in \cF$ and $\gamma_{\bar{s}}(x, u) \in \Gamma$ such that
\beqnas
g_{\bar{s}}u\iota(x) = \xi_{\bar{s}}(x, u)\gamma_{\bar{s}}(x, u).
\eeqnas
By definition, one has
\beqnas
g_{\bar{s}}ux &=&  \xi_{\bar{s}}(x, u)\Gamma, \\
\sigma(g_{\bar{s}}uy) &=& \gamma_{\bar{s}}(x, u)\sigma(y). 
\eeqnas 
For $\bfm_{0} \in \bbZ^{d} \setminus \{0\}$, let $\bfx: X \times V \rightarrow \bbR^{m}$ and $\bfy: X \times V \rightarrow \bbR^{n}$ such that
\beqnas
(\xi_{\bar{s}}(x, u)^{t})^{-1}\bfm_{0} = \left(\begin{array}{cc}
\bfx_{\bar{s}}(x, u) \\
\bfy_{\bar{s}}(x, u) \\
\end{array}
\right).
\eeqnas
\bprop
For $\bar{s} = (s_{1}, \dots, s_{d}) \in \cS^{+} $, $\bfm_{0} \in \bbZ^{d} \setminus \{0\}$, $0 < \epsilon \leq \frac{1}{2}$ and $x \in X$, define the set
\beqna\label{def.ve}
V_{x, \epsilon} = \{u \in V, \|\xi_{\bar{s}}(x, u)\| \leq \epsilon^{-1}, \ \|\bfx_{1}(x, u)\| \geq \epsilon^{2}\|\bfm_{0}\|\}.
\eeqna
Then if 
\beqna
\label{asp.ve}e^{-\kappa_{1}\frac{\lfloor \bar{s} \rfloor}{d}} <  \epsilon  \leq \frac{1}{2},\ \htt(x) < e^{\kappa_{1}\frac{\lfloor \bar{s} \rfloor}{d}},
\eeqna 
we have
 \beq\label{est.gamma.1.s}
 m_{H}(V \setminus V_{x, \epsilon}) \ll \epsilon^{\frac{c_{3}}{2}}
 \eeq
 Moreover, for any $\bfm \in \bbZ^{d}$, 
 \beq\label{est.gamma.3.s}
m_{H}(u \in V_{x, \epsilon}, \gamma_{\bar{s}}(x, u)^{t}\bfm_{0} = \bfm) \ll \epsilon^{-3n}e^{-ns_{1} - (s_{m+1} + \dots + s_{d})}.
\eeq
\eprop

\subsubsection{The time tuple}
Let
\beqnas
\bar{s}_{t} = \diag\{\alpha_{1}t, \dots, \alpha_{d}t\},
\eeqnas
then $g_{\bar{s}_{t}} = g_{t}$.  We  quote a slight adaptation of Lemma 2.9, \cite{BG19}, based on its proof. 
\blem\label{time}
Given any $\{(t_{1}, \dots, t_{k}), t_{i} > 0, t_{k} > t_{k-1}\}$, there exists $\bar{s} \in \cS^{+}$ satisfying
\beqna
\label{t.1}\lfloor \bar{s} \rfloor &\geq& \alpha_{0}(t_{k} - t_{k-1}), \\
\label{t.2}\lfloor \bar{s} - \bar{s}_{t_{k-1}} \rfloor &\geq& \frac{\beta_{0}}{d}(t_{k} - t_{k-1}), \ \beta_{0} = \alpha_{0}\min\{m, n\},
\eeqna
and
\beq\label{t.3}
\bar{s}_{t_{k}} - \bar{s} = \diag\{\frac{z}{m}, \dots, \frac{z}{m}, \frac{z}{n}, \dots, \frac{z}{n}\}
\eeq
for some $z  \geq \frac{mn}{d}\alpha_{0}D(t_{1}, \dots, t_{k})$. 
\elem

\subsection{Estimates of Fourier coefficients}

In this section, we establish the multidimensional analogue of the Fourier decay estimate from Proposition 4.10 in \cite{Kim1}, which plays a key role in the proof of the effective equidistribution. 

We start with the following lemma, which is in the spirit of Lemma 3.4 in \cite{Kim1}. For $\{(t_{1}, \dots, t_{k}), t_{i} > 0, t_{k} > t_{k-1} > \dots > t_{1}\}$, define $\bar{s} \in \cS^{+}$ as in Lemma~\ref{time}. Let
\beqnas
D:= D(t_{1}, \dots, t_{k}) = \min\{t_{k}, |t_{i} - t_{j}|, 1 \leq i \neq j \leq k\}.
\eeqnas

\blem
For $v \in V$, define the action $A_{\bar{s}, v}$ on $V$ by
\beq
A_{\bar{s}, v}u = g_{-\bar{s}}vg_{\bar{s}}u.
\eeq
Then we have for any $y \in Y$, $1 \leq i \leq k-1$, 
\beqna\label{est.s1}
d_{X}(g_{\bar{s}_{t_{i-1}}}A_{\bar{s}, v}u\pi(y), g_{\bar{s}_{t_{i-1}}}u\pi(y))  \ll e^{-\beta_{0}D}d_{G}(v, \id).
\eeqna
Moreover, if $\xi_{\bar{s}_{t_{i-1}}}(\pi(y), u) \in  \cF(r, r^{\frac{1}{d-1}})$ for some $0 < r < \frac{1}{2}$, then for sufficiently large $D$,
\beqna\label{est.g}
\gamma_{\bar{s}_{t_{i-1}}}(\pi(y), u) = \gamma_{\bar{s}_{t_{i-1}}}(\pi(y), A_{\bar{s}, v}u).
\eeqna
\elem

\bpf
By the non-expanding property of $\Phi_{g_{\bar{s}}}(u) = g_{-\bar{s}}ua_{\bar{s}}$ on $H$, we have for $1 \leq i \leq k-1$, $v \in V$,
$$
d_{G}(g_{-(\bar{s}-\bar{s}_{t_{i-1}})}vg_{\bar{s} - \bar{s}_{t_{i-1}}}, \id) \ll e^{- d\lfloor \bar{s} - \bar{s}_{t_{i}} \rfloor }d_{G}(v, \id) \leq e^{- \beta_{0}D}d_{G}(v, \id), 
$$
which implies that 
\beqna\label{est.s2}
\nonumber & &d_{X}(g_{\bar{s}_{t_{i-1}}}A_{\bar{s}, v}u\pi(y), g_{\bar{s}_{t_{i-1}}}u\pi(y)) = d_{X}(g_{\bar{s}_{t_{i-1}}}g_{-\bar{s}}vg_{\bar{s}}u\pi(y), g_{\bar{s}_{t_{i-1}}}u\pi(y))\\
\nonumber&\leq& d_{G}(g_{-(\bar{s}-\bar{s}_{t_{i-1}})}vg_{\bar{s}-\bar{s}_{t_{i-1}}}g_{\bar{s}_{t_{i-1}}}u\pi(y), g_{\bar{s}_{t_{i-1}}}u\pi(y)) \\
&\ll& e^{- \beta_{0}D}d_{G}(v, \id),
\eeqna
where the third line is due to the right invariance of $d_{G}$.

The assumption $ \xi_{\bar{s}_{t_{i-1}}}(\pi(y), u) \in \cF(r, r^{\frac{1}{d-1}})$ implies that $\xi_{\bar{s}_{t_{i-1}}}(\pi(y), u)\Gamma \in K(r^{-\frac{1}{d-1}})$.
By \eqref{ineq.iota}, we have
\beqnas
\|\xi_{\bar{s}_{t_{i-1}}}(\pi(y), u)\| =  \|\iota(\xi_{\bar{s}_{t_{i-1}}}(\pi(y), u)\Gamma)\| \ll \htt(\xi_{\bar{s}_{t_{i-1}}}(\pi(y), u)\Gamma)^{d-1} < r^{-1}, 
\eeqnas
and $\|\xi_{\bar{s}_{t_{i-1}}}^{-1}(\pi(y), u)\| = \|\xi_{\bar{s}_{t_{i-1}}}(\pi(y), u)\| \ll r^{-1}$.
By definitions, 
\beqnas
g_{\bar{s}_{t_{i-1}}}u\iota(\pi(y)) &=& \xi_{\bar{s}_{t_{i-1}}}(\pi(y), u)\gamma_{\bar{s}_{t_{i-1}}}(\pi(y), u),\\
g_{\bar{s}_{t_{i-1}}}A_{\bar{s}, v}u \iota(\pi(y)) &=& \xi_{\bar{s}_{t_{i-1}}}(\pi(y), A_{\bar{s}, v}u )\gamma_{\bar{s}_{t_{i-1}}}(\pi(y), A_{\bar{s}, v}u),
\eeqnas 
we derive that
\beqna\label{est.3}
\nonumber& &d_{G}(\gamma_{\bar{s}_{t_{i-1}}}(\pi(y), A_{\bar{s}, v}u), \gamma_{\bar{s}_{t_{i-1}}}(\pi(y), u)) \\
\nonumber&=& d_{G}(\xi_{\bar{s}_{t_{i-1}}}^{-1}(\pi(y), A_{\bar{s}, v}u)g_{\bar{s}_{t_{i-1}}}A_{\bar{s}, v}u \iota(\pi(y)) , \xi_{\bar{s}_{t_{i-1}}}^{-1}(\pi(y), u)g_{\bar{s}_{t_{i-1}}}u\iota(\pi(y))) \\
\nonumber &\ll& d_{G}(\xi_{\bar{s}_{t_{i-1}}}^{-1}(\pi(y), A_{\bar{s}, v}u), \xi_{\bar{s}_{t_{i-1}}}^{-1}(\pi(y),u)) \\
& &+ \|\xi_{\bar{s}_{t_{i-1}}}^{-1}(\pi(y), u)\|^{2} d_{G}(g_{\bar{s}_{t_{i-1}}}A_{\bar{s}, v}u \iota(\pi(y)) , g_{\bar{s}_{t_{i-1}}}u\iota(\pi(y))).
\eeqna
The second term on the last line is estimated by 
\beqna\label{est.2}
\|\xi_{\bar{s}_{t_{i-1}}}^{-1}(\pi(y), u)\|^{2} d_{G}(g_{\bar{s}_{t_{i-1}}}A_{\bar{s}, v}u \iota(\pi(y)) , g_{\bar{s}_{t_{i-1}}}u\iota(\pi(y))) 
\ll r^{-2}e^{-\beta_{0}D}d_{G}(v, \id).
\eeqna
Sincet $g_{\bar{s}_{t_{i-1}}}u \pi(y) = \xi_{\bar{s}_{t_{i-1}}}(\pi(y), u)\Gamma \in K((r^{\frac{d}{d-1}})^{-\frac{1}{d}})$, $\iota$ is an isometry on $B^{X}(g_{\bar{s}_{t_{i-1}}}u \pi(y), r^{{\frac{d}{d-1}}})$. Also \eqref{est.s2} implies that $g_{\bar{s}_{t_{i-1}}}A_{\bar{s}, v}u \pi(y) \in B^{X}(g_{\bar{s}_{t_{i-1}}}u \pi(y), r^{{\frac{d}{d-1}}})$ for $D$ large enough. Thus, 
\beqna\label{est.4}
d_{G}(\xi_{\bar{s}_{t_{i-1}}}(\pi(y), A_{\bar{s}, v}u), \xi_{\bar{s}_{t_{i-1}}}(\pi(y),u)) =  d_{X}(\xi_{\bar{s}_{t_{i-1}}}(\pi(y), A_{\bar{s}, v}u)\Gamma, \xi_{\bar{s}_{t_{i-1}}}(\pi(y),u)\Gamma).
\eeqna
By the right-invariance of $d_{G}$ and \eqref{est.4}, we deduce that
\beqna\label{est.1}
\nonumber & &d_{G}(\xi_{\bar{s}_{t_{i-1}}}^{-1}(\pi(y), A_{\bar{s}, v}u), \xi_{\bar{s}_{t_{i-1}}}^{-1}(\pi(y),u)) \\
\nonumber&\leq& \|\xi_{\bar{s}_{t_{i-1}}}^{-1}(\pi(y),u)\|^{2}d_{G}(\xi_{\bar{s}_{t_{i-1}}}(\pi(y),u), \xi_{\bar{s}_{t_{i-1}}}(\pi(y), A_{\bar{s}, v}u)) \\
\nonumber&\leq& r^{-2}d_{G}(\xi_{\bar{s}_{t_{i-1}}}(\pi(y), A_{\bar{s}, v}u), \xi_{\bar{s}_{t_{i-1}}}(\pi(y),u)) \\
&\leq& r^{-2}d_{G}(g_{\bar{s}_{t_{i-1}}}A_{\bar{s}, v}u\pi(y), g_{\bar{s}_{t_{i-1}}}u\pi(y)) \ll r^{-2}e^{-\beta_{0}D}d_{G}(v, \id).
\eeqna
Inserting \eqref{est.1} and \eqref{est.2} into  \eqref{est.3}, we obtain
\beqnas
d_{G}(\gamma_{\bar{s}_{t_{i-1}}}(\pi(y), A_{\bar{s}, v}u), \gamma_{\bar{s}_{t_{i-1}}}(\pi(y), u)) \ll r^{-2}e^{-\beta_{0}D}d_{G}(v, \id).
\eeqnas 
By the discreteness of $\Gamma$, there exists $\gamma_{0} := \inf\{d_{G}(\Id, \gamma), \gamma \in \Gamma \setminus \Id\} > 0$. Then for $D$ large enough, we have 
$$
d_{G}(\gamma_{\bar{s}_{t_{i-1}}}(\pi(y), A_{\bar{s}, v}u), \gamma_{\bar{s}_{t_{i-1}}}(\pi(y), u)) < \gamma_{0},
$$
which implies \eqref{est.g}. 
\epf

For  $\{(t_{1}, \dots, t_{k}), t_{i} > 0, t_{k} > t_{k-1} > \dots > t_{1}\}$, define the measure on orbits of flow $(g_{t_{1}}uy_{0}, \dots, g_{t_{k}}uy_{0})$ as
\beqnas
 \mu_{y_{0}, t_{1}, \dots, t_{k}}(\prod^{k}_{i =1}f_{i}) = \int_{V}\prod^{k}_{i =1}f_{i}(g_{t_{i}}uy_{0})dm_{V}(u).
\eeqnas
Define the measure on the projection of the flow on $(\bbT^{d})^{k}$ as
\beqnas
\nu_{y_{0}, t_{1}, \dots, t_{k}} = \sigma_{\ast}\mu_{y_{0}, t_{1}, \dots, t_{k}}.
\eeqnas
Define the probability measure on $Y^{k}$ with respect to the partitions $\{\omega_{j_{i}}\}_{\bfj_{k}}$ introduced in Proposition~\ref{def.ptt}, 
\beqnas
 \mu_{t_{1}, \dots, t_{k}, \bfj_{k}}(\prod^{k}_{i =1}f_{i}) = \pi_{*}\mu_{y_{0}, t_{1}, \dots, t_{k}}(\prod^{k}_{i=1}\omega_{j_{i}})^{-1}\int_{Y^{k}}\prod^{k}_{i=1}f_{i}(y_{i})\omega_{j_{i}}(\pi(y_{i}))d\mu_{y_{0}, t_{1}, \dots, t_{k}}(y_{1}, \dots, y_{k}).
\eeqnas
and
\beqnas
\nu_{t_{1}, \dots t_{k} , \bfj_{k}} = \sigma_{\ast}\mu_{t_{1}, \dots, t_{k},  \bfj_{k}}.
\eeqnas
The Fourier transform of $\nu_{t_{1}, \dots t_{k} , \bfj_{k}}$ is given by
\beqnas
& &\widehat{\nu}_{t_{1}, \dots t_{k} , \bfj_{k}}(\bfm_{1}, \dots, \bfm_{k})  \\
&=& \int_{(\bbT^{d})^{k}} \prod^{k}_{i=1}e^{-2\pi i \bfm_{i} \cdot b_{i}}d\nu_{t_{1}, \dots t_{k} , \bfj_{k}}(b_{1}, \dots, b_{k})\\
&=& (\pi_{*}\mu_{t_{1}, \dots t_{k}}(\prod^{k}_{i =1}\omega_{j_{i}}))^{-1}\int_{V} \prod^{k}_{i=1}e^{-2\pi i \bfm_{i} \cdot \sigma(g_{t_{i}}uy_{0})}\omega_{j_{i}}(\pi(g_{t_{i}}uy_{0})) dm_{V}(u).
\eeqnas


Now we proceed to the estimate of $\widehat{\nu}_{t_{1}, \dots t_{k} , \bfj_{k}}(\bfm_{1}, \dots, \bfm_{k})$. 
\bprop\label{prop.k}
For $y_{0} = (g_{0}, \zero)(1, \bfb_{0})\hat{\Gamma} \in Y$ with $\|g_{0}\| \leq e^{\frac{\alpha_{0} D(t_{1}, \dots, t_{k})}{4}}$, define  
$$
\rho = \max \big( e^{-c_{4}\frac{\alpha_{0}  D(t_{1}, \dots, t_{k})}{d}}, c_{11}^{-\frac{1}{2d}}\zeta(\bfb_{0}, e^{\frac{\alpha_{0} D(t_{1}, \dots, t_{k})}{2}})^{-\frac{1}{2d}}\big)
$$ and $r = \rho^{c^{2}_{4}}$, where $c_{4} = \frac{\min\{\kappa_{1}, c_{3}\}}{2000d^{3}}$ and $c_{11}$ are the same constants as in Proposition 4.10, \cite{Kim1}.  Assume $\htt(\pi(y_{0})) \leq c^{-1}_{1}\rho^{-\frac{1}{10d^{2}(d-1)}}$. 

Define $\{z_{i}\}$, $\{\omega_{i}\}$ as in Proposition~\ref{def.ptt}. Let $\bfj_{k}= (j_{1}, \dots, j_{k}) \in \cJ^{k}$. Assume $\iota(z_{j_{i}}) \in \cF(C^{4}_{1}r, 2C_{2}r^{\frac{1}{d}})$ for $1 \leq i \leq k$. 

Then for any $\overline{\bfm}_{k} = (\bfm_{1}, \dots, \bfm_{k})$ with $0 < \|\overline{\bfm}_{k}\| < \rho^{-c_{4}}$, we have
\beqna\label{est.key.3}
|\hat{ \nu}_{t_{1}, \dots t_{k} , \bfj_{k}}(\bfm_{1}, \dots, \bfm_{k})| = O(\rho^{c^{2}_{4}}).
\eeqna
 \eprop

We take the supremum norm $\|\overline{\bfm}_{k}\| = \max_{i \leq k}\|\bfm_{i}\|$.

\begin{proof}
We prove the estimate by induction. While the case $k=1$ follows from Proposition 5.9, \cite{Kim1}, we now suppose \eqref{est.key.3} is valid for $k=l$, and prove it for $k=l+1$: namely, for $0 < \|\overline{\bfm}_{l+1}\| < \rho_{l+1}^{-c_{4}}$, 
\beqnas
|\hat{ \nu}_{t_{1}, \dots t_{l+1} , \bfj_{l+1}}(\bfm_{1}, \dots, \bfm_{l+1})| = O(\rho_{l+1}^{c^{2}_{4}}), 
\eeqnas
where
\beqna\label{def.rho.l}
\rho_{l+1} =  \max \big( e^{-c_{4}\frac{\alpha_{0}D(t_{1}, \dots, t_{l+1})}{d}}, c_{11}^{-\frac{1}{2d}}\zeta(\sigma(y_{0}), e^{\frac{\alpha_{0}D(t_{1}, \dots, t_{l+1})}{2}})^{-\frac{1}{2d}}\big).
\eeqna
Notice that  $\rho_{l} \leq \rho_{l+1} \leq \rho_{k} = \rho$ for any $l+1 \leq k$. As seen from the proof of Proposition 4.10 in \cite{Kim1}, the estimates holds for $\rho^{c^{2}_{4}} \leq r < C$ for some constant $C > 0$. In our setting, $\{z_{i}\}$, $\{\omega_{i}\}$ are defined with respect to $r = \rho^{c^{2}_{4}} > \rho^{c^{2}_{4}}_{l} > \rho^{c^{2}_{4}}_{1}$; hence, they satisfy the conditions for the case $k=1$ .

Without loss of generality, we assume that $0 < t_{1} < \dots <  t_{l} < t_{l+1}$. By Lemma \ref{time}, there exists $\bar{s} \in \cS^{+}$ associated with $(t_{1}, \dots, t_{l+1})$ such that \eqref{t.1},  \eqref{t.2} and \eqref{t.3} hold.  Moreover, there exists $z > 0$ satisfying
$$
g_{\bar{s}_{t_{l+1}} - \bar{s}} = a_{\frac{z}{mn}}. 
$$

For $v \in H$, define the action $A_{\bar{s}, v}: H \rightarrow H$ by
\beq
A_{\bar{s}, v}u = g_{-\bar{s}}vg_{\bar{s}}u.
\eeq

Let $\omega =  \frac{1}{m_{H}(V)}\un_{V}$, thus $dm_{V} = \omega dm_{H}$. By the invariance of the Haar measure $dm_{H}$, we have
 \beqna\label{eq.3.5.2}
\nonumber & &(\pi_{*}\mu_{t_{1}, \dots t_{l+1}}(\prod^{l+1}_{i =1}\omega_{j_{i}}))\hat{\nu}_{t_{1}, \dots t_{l+1} , \bfj_{l+1}}(\bfm_{1}, \dots, \bfm_{l+1}) \\
\nonumber &=& \int_{H} \prod^{l+1}_{i=1}e^{-2\pi i \bfm_{i} \cdot \sigma(g_{t_{i}}uy_{0})}\omega_{j_{i}}(\pi(g_{t_{i}}uy_{0}))\omega(u) dm_{H}(u)\int_{H} \omega(v)dm_{H}(v)\\
\nonumber &=&  \int_{V}\int_{H} \big( \prod^{l}_{i=1}e^{-2\pi i \bfm_{i} \cdot \sigma(g_{\bar{s}_{t_{i}}}A_{\bar{s}, v}uy_{0})}\omega_{j_{i}}(\pi(g_{\bar{s}_{t_{i}}}A_{\bar{s}, v}uy_{0}))\big) \\
& & \qquad \cdot e^{-2\pi i \bfm_{l+1} \cdot \sigma(g_{\bar{s}_{t_{l+1}}}A_{\bar{s}, v}uy_{0})}\omega_{j_{l+1}}(\pi(g_{\bar{s}_{t_{l+1}}}A_{\bar{s}, v}uy_{0}))\omega(A_{\bar{s}, v}u)dm_{H}(u) dm_{V}(v).
\eeqna
By the the non-expanding property of $\Phi_{g_{\bar{s}}}(v) = g_{-\bar{s}}vg_{\bar{s}}$ on $H$,  we have
\beqnas
d_{G}(A_{\bar{s}, v}u, u) = d_{G}(g_{-\bar{s}}vg_{\bar{s}}, \id) \ll e^{-d\lfloor \bar{s} \rfloor }d_{G}(v, \id).
\eeqnas
Thus there exists a subset $V' \subset H$ with smooth boundary and compact closure, which also contains the identity in $H$, such that $A_{\bar{s}, v}u \in V$ implies $u \in V'$. Moreover, $m_{H}(V') \asymp m_{H}(V)$.  Then we write
\beqnas\label{eq.3.5.2}
\nonumber & &(\pi_{*}\mu_{t_{1}, \dots t_{l+1}}(\prod^{l+1}_{i =1}\omega_{j_{i}}))\hat{\nu}_{t_{1}, \dots t_{l+1} , \bfj_{l+1}}(\bfm_{1}, \dots, \bfm_{l+1}) \\
\nonumber &=&  \int_{V}\int_{V'} \big( \prod^{l}_{i=1}e^{-2\pi i \bfm_{i} \cdot \sigma(g_{\bar{s}_{t_{i}}}A_{\bar{s}, v}uy_{0})}\omega_{j_{i}}(\pi(g_{\bar{s}_{t_{i}}}A_{\bar{s}, v}uy_{0}))\big) \\
& & \qquad \cdot e^{-2\pi i \bfm_{l+1} \cdot \sigma(g_{\bar{s}_{t_{l+1}}}A_{\bar{s}, v}uy_{0})}\omega_{j_{l+1}}(\pi(g_{\bar{s}_{t_{l+1}}}A_{\bar{s}, v}uy_{0}))\frac{1}{m_{H}(V)}dm_{H}(u) dm_{V}(v).
\eeqnas
%
For $1 \leq i \leq l$, by \eqref{est.s1}
\beqnas
d_{X}(g_{\bar{s}_{t_{i}}}A_{\bar{s}, v}u \pi(y_{0}), g_{\bar{s}_{t_{i}}}u\pi(y_{0})) \ll e^{-\beta_{0}D(t_{1}, \dots, t_{l+1})}d_{G}(v, \id),
\eeqnas
which implies that
\beqnas
|\omega_{j_{i}}(g_{\bar{s}_{t_{i}}}A_{\bar{s}, v}u \pi(y_{0})) - \omega_{j_{i}}(g_{\bar{s}_{t_{i}}}u \pi(y_{0}))| \ll e^{-\beta_{0}D(t_{1}, \dots, t_{l+1})}\|\nabla \omega_{j_{i}}\|, 
\eeqnas
and by the assumption $\|\nabla \omega_{j_{i}}\|_{L^{\infty}} \leq r^{-1}$,
\beqna\label{est.3.5.1}
|\prod^{l}_{i=1}\omega_{j_{i}}(g_{\bar{s}_{t_{i}}}A_{\bar{s}, v}u \pi(y_{0})) - \prod^{l}_{i=1}\omega_{j_{i}}(g_{\bar{s}_{t_{i}}}u \pi(y_{0}))| \ll e^{-\beta_{0}D(t_{1}, \dots, t_{l+1})l}r^{-l}.
\eeqna
Observe that the integral is taken over $\{u \in V, \prod^{l+1}_{i=1}\omega_{j_{i}}(g_{\bar{s}_{t_{i}}}u \pi(y_{0})) \neq 0\} \cup \{u \in V', \prod^{l}_{i=1}\omega_{j_{i}}(g_{\bar{s}_{t_{i}}}u \pi(y_{0})) \neq 0\}$, whose measure is controlled by $\pi_{*}\mu_{y_{0}, t_{1}, \dots t_{l+1}}(\prod^{l+1}_{i =1}\omega_{j_{i}})$, up to a multiplicative constant. We assume that $\pi(g_{\bar{s}_{t_{i}}}uy_{0}) = \xi_{\bar{s}_{t_{i}}}(\pi(y_{0}), u)\Gamma \in \supp \ \omega_{j_{i}}$ for each $1 \leq i \leq l$, otherwise the integral is trivial. Since $\supp \ \omega_{j_{i}} \subset B_{C^{3}_{1}r}(z_{j_{i}})$ and $\iota(z_{j_{i}}) \in \cF(C^{4}_{1}r, 2C_{2}r^{\frac{1}{d}})$, we deduce $\iota(\supp \ \omega_{j_{i}}) \subset \cF(r, r^{\frac{1}{d-1}})$ such that $\xi_{\bar{s}_{t_{i}}}(\pi(y_{0}), u) \in  \cF(r, r^{\frac{1}{d-1}})$. 
Then due to \eqref{est.g}, for $1 \leq i \leq l$ and $D$ large enough,
\beqnas
\gamma_{\bar{s}_{t_{i}}}(\pi(y_{0}), A_{\bar{s}, v}u) = \gamma_{\bar{s}_{t_{i}}}(\pi(y_{0}), u),
\eeqnas
such that
\beqna\label{eq.3.5.3}
 \sigma(g_{\bar{s}_{t_{i}}}A_{\bar{s}, v}uy_{0}) = \gamma_{\bar{s}_{t_{i}}}(\pi(y_{0}), A_{\bar{s}, v}u)\sigma(y_{0}) = \gamma_{\bar{s}_{t_{i}}}(\pi(y_{0}), u)\sigma(y_{0}).
\eeqna
By \eqref{eq.3.5.2}, \eqref{est.3.5.1} and \eqref{eq.3.5.3}, we derive that
\beqna\label{eq.3.5.4}
\nonumber& &(\pi_{*}\mu_{t_{1}, \dots t_{l+1}}(\prod^{l+1}_{i =1}\omega_{j_{i}}))\hat{\nu}_{t_{1}, \dots t_{l+1} , \bfj_{l+1}}(\bfm_{1}, \dots, \bfm_{l+1}) \\
\nonumber &=& \int_{V'} \big( \prod^{l}_{i=1}e^{-2\pi i \bfm_{i} \cdot \gamma_{\bar{s}_{t_{i}}}(\pi(y_{0}), u) \sigma(y_{0})}\omega_{j_{i}}(\pi(g_{\bar{s}_{t_{i}}}uy_{0}))\big)I(y_{\bar{s}, u}, \bfm_{l+1})\frac{1}{m_{H}(V)}dm_{H}(u) \\
& &  \qquad   \qquad  + O_{V}\big(\pi_{*}\mu_{y_{0}, t_{1}, \dots t_{l+1}}(\prod^{l+1}_{i =1}\omega_{j_{i}})e^{-\beta_{0}D(t_{1}, \dots, t_{l+1})l}r^{-l}\big),
\eeqna
where $y_{\bar{s}, u} = g_{\bar{s}}uy_{0}$ and
\beqna\label{def.l1}
\nonumber I(y_{\bar{s}, u}, \bfm_{l+1}) &=& \int_{V}e^{-2\pi i \bfm_{l+1} \cdot \sigma(g_{\bar{s}_{t_{l+1}}-\bar{s}}vy_{\bar{s}, u})}\omega_{j_{l+1}}(\pi(g_{\bar{s}_{t_{l+1}}-\bar{s}}vy_{\bar{s}, u}))dm_{V}(v)\\
&=& \int_{V}e^{-2\pi i \bfm_{l+1} \cdot \sigma(a_{\frac{z}{mn}}vy_{\bar{s}, u})}\omega_{j_{l+1}}(\pi(a_{\frac{z}{mn}}vy_{\bar{s}, u}))dm_{V}(v),
\eeqna 
for $z \geq \frac{mn}{d}\alpha_{0}D(t_{1}, \dots, t_{l+1})$ by \eqref{t.3}.
%
Let $\epsilon_{l+1} = \rho_{l+1}^{\frac{1}{10d^{2}(d-1)}}$, and define the set  $V_{x_{0}, c_{1}\epsilon_{l+1}} \subset V'$ as \eqref{def.ve}.  The assumptions~\eqref{asp.ve} are satisfied: $c_{1}\epsilon_{l+1} \geq \epsilon_{l+1} > e^{-\kappa_{1}\frac{\alpha_{0}D(t_{1}, \dots, t_{l+1})}{d}}$,
and $\htt(\pi(y_{0})) \leq c^{-1}_{1}\rho^{-\frac{1}{10d^{2}(d-1)}} \leq c^{-1}_{1}\rho_{l+1}^{-\frac{1}{10d^{2}(d-1)}} = c^{-1}_{1}\epsilon^{-1}_{l+1}$. Moreover, we have $\htt(\pi(y_{\bar{s}, u})) = \htt(g_{\bar{s}}ux_{0}) \leq c_{1}\|\xi_{\bar{s}}(x_{0}, u)\| \leq \epsilon_{l+1}^{-1}$ for $u \in V_{x_{0}, c_{1}\epsilon_{l+1}}$.

Define
\beqnas
\rho_{\bar{s}, u} &=&  \max \big( e^{-c_{4}\frac{z}{mn}}, c_{11}^{-\frac{1}{2d}}\zeta(\sigma(y_{\bar{s}, u}), e^{\frac{z}{2m}})^{-\frac{1}{2d}}\big), \\
V_{1} &=& \{u \in V_{x_{0}, c_{1}\epsilon_{l+1}}, \rho_{\bar{s}, u} \leq \rho_{l+1}\}.
\eeqnas
Then for $u \in V_{1}$,  the assumptions in Prop 4.10, \cite{Kim1} hold: $\htt(\pi(y_{\bar{s}, u})) \leq \epsilon_{l+1}^{-1} = \rho_{l+1}^{-\frac{1}{10d^{2}(d-1)}} \leq \rho_{\bar{s}, u}^{-\frac{1}{10d^{2}(d-1)}}$, $r > \rho_{l+1}^{c^{2}_{4}} \geq \rho_{\bar{s}, u}^{c^{2}_{4}} $.  Consequently, for $0 \leq \|\bfm_{l+1}\| \leq \rho_{l+1}^{-c_{4}} \leq \rho_{\bar{s}, u}^{-c_{4}}$, there exists some constant $C > 0$ such that
$$
|\hat{\nu}_{y_{\bar{s}, u},\frac{z}{mn} }(\bfm_{l+1})| \leq C\rho^{c^{2}_{4}}_{\bar{s}, u} \leq C\rho_{l+1}^{c^{2}_{4}}. 
$$
Since $\htt(\pi(y_{\bar{s}, u})) \leq  e^{\kappa_{1}\frac{z}{mn}}$, we deduce by \eqref{est.v.o.1} that
$$
\pi_{*}\mu_{y_{\bar{s}, u}, \frac{z}{mn}}(\omega_{j_{l+1}}) \asymp m_{X}(B_{r}(z_{l+1}))  \asymp r^{d^{2}-1}. 
$$
Thus for $u \in V_{1}$, 
\beqna\label{est.l.1}
|I(y_{\bar{s}, u}, \bfm_{l+1})| = \pi_{*}\mu_{y_{\bar{s}, u}, \frac{z}{mn}}(\omega_{j_{l+1}})|\hat{\nu}_{y_{\bar{s}, u}, \frac{z}{mn}}(\bfm_{l+1})| \ll \rho_{l+1}^{c^{2}_{4}}m_{X}(B_{r}(z_{l+1})).
\eeqna
Now we estimate the measure of $V_{x_{0},  c_{1}\epsilon_{l+1}} \setminus V_{1} = \{u \in  V_{x_{0}, c_{1}\epsilon_{l+1}}, \rho_{\bar{s}, u} > \rho_{l+1}\}$. Since $e^{-c_{4}\frac{z}{mn}} \leq e^{-c_{4}\alpha_{0}\frac{D(t_{1}, \dots, t_{l+1})}{d}}$, we have
\beqnas
V_{x_{0},  c_{1}\epsilon_{l+1}} \setminus V_{1}  \subset  \{u \in  V_{x_{0}, c_{1}\epsilon_{l+1}}, c_{11}^{-\frac{1}{2d}}\zeta(\sigma(y_{\bar{s}, u}), e^{\frac{z}{2m}})^{-\frac{1}{2d}} > \rho_{l+1}\}.
\eeqnas
By its definition, $\zeta(\sigma(y_{\bar{s}, u}), e^{\frac{z}{2m}}) < c_{11}^{-1}\rho_{l+1}^{-2d}$ implies that there exists a positive integer $1 \leq q_{0} \leq c_{11}^{-1}\rho_{l+1}^{-2d}$, such that
\beqnas
\|q_{0}  \sigma(y_{\bar{s}, u})\|_{\bbZ} = \|q_{0}  \gamma_{\bar{s}}(x_{0}, u) \cdot \sigma(y_{0})\|_{\bbZ} < \frac{c_{11}^{-2}\rho_{l+1}^{-4d}}{e^{\frac{z}{2m}}} \leq c_{11}^{-2}\rho_{l+1}^{-4d + \frac{n}{2c_{4}}},
\eeqnas
where $\sigma(y_{\bar{s}, u}) = \gamma_{\bar{s}}(x_{0}, u)\sigma(y_{0})$. 

Write $\sigma(y_{0}) = \gamma_{0}b_{0}$ with $\gamma_{0} \in \SL(d, \bbZ)$ and $b_{0} \in \bbT^{d}$. We derive that
\beqnas
& & m_{H}\big(u \in V_{x_{0}, c_{1}\epsilon_{l+1}}, \zeta(\sigma(y_{\bar{s}, u}), e^{\frac{z}{2m}}) < c_{11}^{-1}\rho_{l+1}^{-2d}\big)\\
&\leq& m_{H}\big(u \in V_{x_{0}, c_{1}\epsilon_{l+1}}, \|q_{0}  \gamma_{\bar{s}}(x_{0}, u) \cdot \sigma(y_{0})\|_{\bbZ} < c_{11}^{-2}\rho_{l+1}^{-4d + \frac{n}{2\kappa_{4}}},  1 \leq q_{0} \leq c_{11}^{-1}\rho_{l+1}^{-2d}\big)\\
&\leq& m_{H}\big(u \in V_{x_{0}, c_{1}\epsilon_{l+1}}, \gamma^{tr}_{\bar{s}}(x_{0}, u) \sigma(y_{0})= \bfm, \|q_{0}\bfm\|_{\bbZ} < c_{11}^{-2}\rho_{l+1}^{-4d + \frac{n}{2\kappa_{4}}}, 1 \leq q_{0} \leq c_{11}^{-1}\rho_{l+1}^{-2d}\big)\\
&\leq& \sum_{1 \leq q \leq c_{11}^{-1}\rho_{l+1}^{-2d}}\sum_{\|q\bfm\|_{\bbZ} < c_{11}^{-2}\rho_{l+1}^{-4d + \frac{n}{2c_{4}}}}m_{H}\big(u \in V_{x_{0}, c_{1}\epsilon_{l+1}}, \gamma^{tr}_{\bar{s}}(x_{0}, u)\sigma(y_{0})  = \bfm \big)\\
&\ll& \rho_{l+1}^{ -2d-4d^{2} + \frac{nd}{2c_{4}}-\frac{3n}{10d^{2}(d-1)}}e^{-ns_{1} - (s_{m+1} + \dots + s_{d})},
\eeqnas
where the last line is due to \eqref{est.gamma.3.s}. By \eqref{t.1}, one has 
$$
e^{-ns_{1} - (s_{m+1} + \dots + s_{d})} \leq e^{-d \lfloor \bar{s} \rfloor } \leq e^{-d\alpha_{0}D(t_{1}, \dots, t_{l+1})} \leq \rho_{l+1}^{\frac{d^{2}}{c_{4}}}, 
$$
such that $$
m_{H}\big(u \in V_{x_{0}, c_{1}\epsilon_{l+1}}, \zeta(\sigma(y_{\bar{s}, u}), e^{\frac{z}{2m}}) < c_{11}^{-1}\rho_{l+1}^{-2d}\big) \ll \rho_{l+1}, 
$$
and 
\beqnas
m_{H}\big( V_{x_{0}, c_{1}\epsilon_{l+1}} \setminus V_{1} \big) \ll \rho_{l+1}. 
\eeqnas
Thus with \eqref{est.gamma.1.s}, we have
\beqna\label{est.exp1}
m_{H}\big( V' \setminus V_{1} \big) \ll \rho_{l+1} + \epsilon_{l+1}^{\frac{c_{3}}{2}}.
\eeqna
We write
\beqna\label{ineq.3.5.5}
\nonumber  \int_{V'} \big( \prod^{l}_{i=1}e^{-2\pi i \bfm_{i} \cdot \gamma_{\bar{s}_{t_{i}}}(\pi(y_{0}), u) \sigma(y_{0})}\omega_{j_{i}}(\pi(g_{\bar{s}_{t_{i}}}uy_{0}))\big)I(y_{\bar{s}, u}, \bfm_{l+1})\frac{1}{m_{H}(V)}dm_{H}(u) := I_{1} + I_{2}, 
\eeqna
where 
\beqnas
I_{1} &=&  \int_{V_{1}} \big( \prod^{l}_{i=1}e^{-2\pi i \bfm_{i} \cdot \gamma_{\bar{s}_{t_{i}}}(\pi(y_{0}), u) \sigma(y_{0})}\omega_{j_{i}}(\pi(g_{\bar{s}_{t_{i}}}uy_{0}))\big)I(y_{\bar{s}, u}, \bfm_{l+1})\frac{1}{m_{H}(V)}dm_{H}(u), \\
I_{2} &=&   \int_{V'  \setminus V_{1}}  \big( \prod^{l}_{i=1}e^{-2\pi i \bfm_{i} \cdot \gamma_{\bar{s}_{t_{i}}}(\pi(y_{0}), u) \sigma(y_{0})}\omega_{j_{i}}(\pi(g_{\bar{s}_{t_{i}}}uy_{0}))\big)I(y_{\bar{s}, u}, \bfm_{l+1})\frac{1}{m_{H}(V)}dm_{H}(u).
\eeqnas
By \eqref{est.l.1}, 
\beqnas
I_{1} \ll  \int_{V_{1}} \prod^{l}_{i=1}e^{-2\pi i \bfm_{i} \cdot \gamma_{\bar{s}_{t_{i}}}(\pi(y_{0}), u) \sigma(y_{0})}\omega_{j_{i}}(\pi(g_{\bar{s}_{t_{i}}}uy_{0})) dm_{H}(u) \cdot \rho_{l+1}^{c^{2}_{4}}m_{X}(B_{r}(z_{l+1})).
\eeqnas
The assumption that \eqref{est.key.3} holds for $k = l$ implies that for $\|(\bfm_{1}, \dots, \bfm_{l})\| \leq \rho_{l}^{-c_{4}}$
\beqnas
\int_{V} \big( \prod^{l}_{i=1}e^{-2\pi i \bfm_{i} \cdot \gamma_{\bar{s}_{t_{i}}}(\pi(y_{0}), u) \sigma(y_{0})}\omega_{j_{i}}(\pi(g_{\bar{s}_{t_{i}}}uy_{0}))\big)dm_{V}(u) \ll \pi_{*}\mu_{y_{0}, t_{1}, \dots t_{l}}(\prod^{l}_{i =1}\omega_{j_{i}})\rho^{c^{2}_{4}}_{l}.
\eeqnas
On the other hand, by \eqref{eq.3.5.3} and \eqref{eq.3.5.4}, we have
\beqnas
& &\int_{V} \big( \prod^{l}_{i=1}e^{-2\pi i \bfm_{i} \cdot \gamma_{\bar{s}_{t_{i}}}(\pi(y_{0}), u) \sigma(y_{0})}\omega_{j_{i}}(\pi(g_{\bar{s}_{t_{i}}}uy_{0}))\big)dm_{V}(u)\\
&=& \int_{V} \int_{V'} \big( \prod^{l}_{i=1}e^{-2\pi i \bfm_{i} \cdot \gamma_{\bar{s}_{t_{i}}}(\pi(y_{0}), A_{\bar{s}, v}u) \sigma(y_{0})}\omega_{j_{i}}(\pi(g_{\bar{s}_{t_{i}}}A_{\bar{s}, v}uy_{0}))\big)\frac{1}{m_{H}(V)}dm_{H}(u)dm_{V}(v)\\
&=& \int_{V'} \big( \prod^{l}_{i=1}e^{-2\pi i \bfm_{i} \cdot \gamma_{\bar{s}_{t_{i}}}(\pi(y_{0}), u) \sigma(y_{0})}\omega_{j_{i}}(\pi(g_{\bar{s}_{t_{i}}}uy_{0}))\big)\frac{1}{m_{H}(V)}dm_{H}(u)\\
& & + O_{V}\big(\pi_{*}\mu_{y_{0}, t_{1}, \dots t_{l}}(\prod^{l}_{i =1}\omega_{j_{i}})e^{-\beta_{0}D(t_{1}, \dots, t_{l})l}r^{-l}\big),
\eeqnas
which implies that
\beqnas
\int_{V'} \big( \prod^{l}_{i=1}e^{-2\pi i \bfm_{i} \cdot \gamma_{\bar{s}_{t_{i}}}(\pi(y_{0}), u) \sigma(y_{0})}\omega_{j_{i}}(\pi(g_{\bar{s}_{t_{i}}}uy_{0}))\big)dm_{H}(u) \ll \pi_{*}\mu_{y_{0}, t_{1}, \dots t_{l}}(\prod^{l}_{i =1}\omega_{j_{i}})\big(\rho^{c^{2}_{4}}_{l} + e^{-\beta_{0}D(t_{1}, \dots, t_{l})l}r^{-l}\big).
\eeqnas
Hence, 
\beqna\label{est.e1.0}
 I_{1}  \ll \pi_{*}\mu_{y_{0}, t_{1}, \dots t_{l}}(\prod^{l}_{i =1}\omega_{j_{i}})m_{X}(B_{r}(z_{l+1}))(\rho^{2c^{2}_{4}}_{l+1} + \rho^{c^{2}_{4}}_{l+1} e^{-\beta_{0}D(t_{1}, \dots, t_{l})l}r^{-l}).
\eeqna
For $I_{2}$, observe that
\beqnas
\nonumber I_{2} 
&\leq& \int_{H} \un_{V' \setminus V_{1}}(u) \prod^{l}_{i=1}\omega_{j_{i}}(\pi(g_{\bar{s}_{t_{i}}}uy_{0}))\int_{V}\omega_{j_{l+1}}(\pi(g_{\bar{s}_{t_{l+1}}}A_{\bar{s}, v}uy_{0}))dm_{V}(v)\frac{1}{m_{H}(V)}dm_{H}(u)\\
\nonumber &\leq& \int_{H} \un_{V' \setminus V_{1}}(u) \prod^{l}_{i=1}\omega_{j_{i}}(\pi(g_{\bar{s}_{t_{i}}}uy_{0}))\frac{1}{m_{H}(V)}dm_{H}(u).
 \eeqnas
Applying \eqref{exp.m.X} with an approximation of $\un_{V' \setminus V_{1}}(u)$, we derive that
\beqnas\label{est.e1.2}
\nonumber  \int_{H}\un_{V' \setminus V_{1}}(u)  \prod^{l}_{i=1}\omega_{j_{i}}(\pi(g_{t_{i}}uy_{0}))\frac{1}{m_{H}(V)}dm_{H}(u) &\asymp& \frac{1}{m_{H}(V)}m_{H}(V' \setminus V_{1}) \prod^{l}_{i=1}m_{X}(B_{r}(z_{j_{i}})) \\
 &\asymp&  (\rho_{l+1} + \epsilon_{l+1}^{\frac{c_{3}}{2}}) \prod^{l}_{i=1}m_{X}(B_{r}(z_{j_{i}})),
\eeqnas
where the last line is due to \eqref{est.exp1}.
Thus, 
\beqna\label{est.e2}
I_{2} = O\big( (\rho_{l+1} + \epsilon_{l+1}^{\frac{\kappa_{3}}{2}})\prod^{l}_{i=1}m_{X}(B_{r}(z_{j_{i}}))\big).
\eeqna
Notice that by \eqref{est.v.o.k.1}, 
\beqnas
\pi_{*}\mu_{y_{0}, t_{1}, \dots t_{l+1}}(\prod^{l+1}_{i =1}\omega_{j_{i}}) \asymp r^{(d^{2}-1)(l+1)}, \  \pi_{*}\mu_{y_{0}, t_{1}, \dots t_{l+1}}(\prod^{l}_{i =1}\omega_{j_{i}}) \asymp r^{(d^{2}-1)l}.
\eeqnas
Then combining \eqref{eq.3.5.4}, \eqref{est.e1.0} and \eqref{est.e2} together, we derive that for $\|(\bfm_{1}, \dots, \bfm_{l+1})\| \leq \rho_{l+1}^{-c_{4}}$, 
\beqnas
|\hat{\nu}_{t_{1}, \dots t_{l+1} , \bfj_{l+1}}(\bfm_{1}, \dots, \bfm_{l+1})| 
\ll  O_{V}\big( e^{- \beta_{0}D(t_{1}, \dots, t_{l+1})l}r^{-l} + \rho_{l+1}^{2c^{2}_{4}} + (\rho_{l+1} + \epsilon_{l+1}^{\frac{c_{3}}{2}})r^{-(d^{2}-1)} \big).
\eeqnas
Notice that
\beqnas
(e^{- \beta_{0}D(t_{1}, \dots, t_{l+1})}r^{-1})^{l} &\leq& (\rho_{l+1}^{\frac{d}{2c_{4}} - c^{2}_{4}})^{l} \ll \rho_{l+1}, \\
\epsilon_{l+1}^{\frac{c_{3}}{2}}r^{-(d^{2}-1)} &\leq& \rho_{l+1}^{\frac{c_{3}}{20d^{2}(d-1)}- (d^{2}-1)c^{2}_{4}} \leq  \rho_{l+1}^{c^{2}_{4}}, \\
\rho_{l+1} r^{-(d^{2}-1)} &\leq&  \rho_{l+1}^{1- d^{2}c^{2}_{4}} \rho_{l+1}^{c^{2}_{4}} \leq  \rho_{l+1}^{c^{2}_{4}},
\eeqnas
which leads to \eqref{est.key.3}. Thus we finish the proof.
\end{proof}

\subsection{Proof of Theorem~\ref{thm.2}}
With Proposition~\ref{prop.k}, we prove the mixing property for multiple functions by modifying the approach in \cite{Kim1}. 

\bpf
Let
\beqnas
\rho =  \max \big( e^{-c_{4}\frac{\alpha_{0}D(t_{1}, \dots, t_{k})}{d}}, c_{11}^{-\frac{1}{2d}}\zeta(\sigma(y_{0}), e^{\frac{\alpha_{0}D(t_{1}, \dots, t_{k})}{2}})^{-\frac{1}{2d}}\big).
\eeqnas
By \eqref{z.p.3}, $\zeta(\sigma(y_{0}), e^{\frac{\alpha_{0}D(t_{1}, \dots, t_{k})}{2}}) \leq (e^{\frac{\alpha_{0}D(t_{1}, \dots, t_{k})}{2}})^{\frac{d}{2d+1}}$. Then we can choose $\delta_{1} > 0$ such that  
 \beqnas
 \zeta(\sigma(y_{0}), e^{\frac{\alpha_{0}D(t_{1}, \dots, t_{k})}{2}})^{\delta_{1}} \leq e^{\frac{\alpha_{0}D(t_{1}, \dots, t_{k})}{4}}, \ \zeta(\sigma(y_{0}), e^{\frac{\alpha_{0}D(t_{1}, \dots, t_{k})}{2}})^{\delta_{1}} \leq \rho^{-\frac{1}{10d^{2}(d-1)}},
\eeqnas
then $\|g_{0}\| \leq \zeta(\sigma(y_{0}), e^{\frac{\alpha_{0}D(t_{1}, \dots, t_{k})}{2}})^{\delta_{2}}$ for some $0 < \delta_{2} < \delta_{1}$ implies that 
\beqnas
\|g_{0}\| \leq e^{\frac{\alpha_{0}D(t_{1}, \dots, t_{k})}{4}}, \ \htt(\pi(y_{0})) \leq c^{-1}_{1}\rho^{-\frac{1}{10d^{2}(d-1)}},
\eeqnas
so the assumptions in Proposition~\ref{prop.k} are satisfied. 

Let $\pi: Y \rightarrow X$ be the natural projection, and $m_{\pi^{-1}(x)}$ be the normalized Haar measure of $\pi^{-1}(x)$, which can viewed as the Lebesgue measure on the torus $\bbT^{d}$. For any $y \in Y$, $x = \pi(y)$, define 
\beqnas
\bar{f}_{i}(x) &=& \int_{\pi^{-1}(x)}f_{i}(y)dm_{\pi^{-1}(x)}(y),\\ 
h_{i}(y) &=& f_{i}(y) - \bar{f}_{i}(\pi(y)),
\eeqnas 
then $\int_{\pi^{-1}(x)}h_{i}(y)dm_{\pi^{-1}(x)}(y) = 0$. Moreover, $h_{i} \in \cC^{\infty}_{c}(Y)$, satisfying 
\beqnas
\cS(h_{i}) \ll \cS(f_{i}), \ \ \|h_{i}\|_{L^{\infty}(Y)} \ll \|f_{i}\|_{L^{\infty}(Y)}, 
\eeqnas
and for $g \in B^{G}(\id, r)$
\beqnas
|h_{i}(gy) - h_{i}(y)| \ll r\cS(f_{i}).
\eeqnas
%
%
Let $K = \{1, 2, \dots, k\}$, and denote $J_{k'}$ as a subset of $K$ with $|J_{k'}| = k'$. Then we can write
\beqna\label{eq.1}
\nonumber & &\mu_{t_{1}, \dots, t_{k}}(\prod^{k}_{i =1}f_{i}) = \int_{V}\prod^{k}_{i =1}f_{i}(g_{t_{i}}uy_{0}) dm_{V}(u) = \mu_{t_{1}, \dots, t_{k}}(\prod^{k}_{i =1}(\bar{f}_{i} \circ \pi + h_{i}))\\
 &=&\pi_{*}\mu_{t_{1}, \dots, t_{k}}(\prod^{k}_{i =1}\bar{f}_{i})  +  \sum^{k-1}_{k'= 0}\sum_{J_{k'} \subset K}\int_{V}\prod_{i \notin J_{k'}}h_{i}(g_{t_{i}}uy_{0}) \prod_{i \in J_{k'}}\bar{f}_{i}(\pi(g_{t_{i}}uy_{0})) dm_{V}(u).
\eeqna

For the first term, 
 by \eqref{exp.m.X}, we have for some positive $l$
\beqna\label{eq.2}
\nonumber  \pi_{*}\mu_{t_{1}, \dots, t_{k}}(\prod^{k}_{i =1}\bar{f}_{i})  &=& \prod^{k}_{i =1}\int_{X}\bar{f}_{i}dm_{X} + O(e^{-\lambda' D(t_{1}, \dots, t_{k})}\prod^{k}_{i =1}\cS_{l}(\bar{f}_{i}))\\
 &\leq& \prod^{k}_{i =1}\int_{Y}f_{i}dm_{Y} + O(e^{- \lambda' D(t_{1}, \dots, t_{k})}\prod^{k}_{i =1}\cS_{l}(f_{i})).
\eeqna


Let $r = \rho^{c^{2}_{4}}$, $c_{5} = \frac{c^{2}_{4}}{3dk}$, and $\epsilon = \rho^{\frac{c_{5}}{2l_{0}(d-1)}} = \rho^{\frac{c^{2}_{4}}{6dkl_{0}(d-1)}}$. We adopt the partition in Proposition~\ref{def.ptt}, such that $\un_{X} = \sum_{j \in \cJ}\omega_{j}$ and $|\cJ| = \cN_{r} \asymp r^{-(d^{2}-1)}$. Define $\cJ_{int} \subset \cJ$ to be the subset of $\cJ$ with $\iota(z_{j}) \in \cF(C^{4}_{1}r, \epsilon) \subset \cF(C^{4}_{1}r, C_{2}r^{\frac{1}{d}}) $.     
Then by the same arguments as in \cite{Kim1},  we have $\iota(B_{r}z_{j}) \in \cF \setminus \cF(C^{5}_{1}r, \frac{\epsilon}{2})$ for each $j \in \cJ \setminus \cJ_{int}$, such that
\beqna\label{int.q}
|\cJ \setminus \cJ_{int}| \asymp \frac{m_{G}(\cF \setminus \cF(C^{5}_{1}r, \frac{\epsilon}{2}))}{m_{G}(\iota(B_{r}z_{j}))} \leq \max\{r^{c_{3}}, \epsilon^{d}\}r^{-(d^{2}-1)}. 
\eeqna

Define
\beqnas
\cJ^{k} = \{\bfj_{k}= (j_{1}, j_{2}, \dots, j_{k}), \ j_{i} \in \cJ, \ 1 \leq i \leq k\},
\eeqnas
and for $\bfj_{k} \in \cJ^{k}$, define
\beqnas
\mu_{(t_{1}, \dots, t_{k}), \bfj_{k}}(\prod^{k}_{i =1}f_{i}) = \pi_{*}\mu_{(t_{1}, \dots, t_{k})}(\prod_{\bfj_{k}}\omega_{j_{i}})^{-1}\int_{Y}\prod^{k}_{i = 1}\big(\omega_{j_{i}} \circ \pi(y_{i})f_{i}(y_{i})\big)d \mu_{t_{1}, \dots, t_{k}}(y_{1}, \dots, y_{k}).
\eeqnas
Then we can write 
\beqnas
\int_{V}\prod_{i \notin J_{k'}}h_{i}(g_{t_{i}}uy_{0}) \prod_{i \in J_{k'}}\bar{f}_{i}(\pi(g_{t_{i}}uy_{0}))dm_{V}(u) = \sum_{\bfj_{k} \in \cJ^{k}} \pi_{*}\mu_{t_{1}, \dots, t_{k}}(\prod_{\bfj_{k}}  \omega_{j_{i}}) \mu_{(t_{1}, \dots, t_{k}), \bfj_{k} }(\prod_{i \notin J_{k'}}h_{i}\prod_{i \in J_{k'}}\bar{f}_{i} \circ \pi).
\eeqnas

Now we decompose $\cJ^{k}$ into two parts, $\cJ^{k}_{1}$ and $\cJ^{k}_{2}$, defined by
\beqnas
\cJ^{k}_{1} &=& \{\bfj_{k}=(j_{1}, \dots, j_{k}) \in \cJ^{k}, j_{i} \in \cJ_{int} \}, \\
\cJ^{k}_{2} &=&  \{\bfj_{k}=(j_{1}, \dots, j_{k}) \in \cJ^{k}, j_{i} \in \cJ \setminus \cJ_{int} \},\\
\cJ^{k}_{3} &=& \cJ^{k} \setminus (\cJ^{k}_{1} \cap \cJ^{k}_{2}), 
\eeqnas
such that
\beqnas
& &\int_{V}\prod_{i \notin J_{k'}}h_{i}(g_{t_{i}}uy_{0}) \prod_{i \in J_{k'}}\bar{f}_{i}(\pi(g_{t_{i}}uy_{0})) dm_{V}(u)\\
&=& \sum_{\bfj_{k, 1} \in \cJ^{k}_{1}} \pi_{*}\mu_{t_{1}, \dots, t_{k}}(\prod_{\bfj_{k, 1}}  \omega_{j_{i}}) \mu_{(t_{1}, \dots, t_{k}), \bfj_{k, 1}}(\prod_{i \notin J_{k'}}h_{i}\prod_{i \in J_{k'}}\bar{f}_{i} \circ \pi)\\
& & +  \sum_{\bfj_{k, 2} \in \cJ^{k}_{2}} \pi_{*}\mu_{t_{1}, \dots, t_{k}}(\prod_{\bfj_{k, 2}}  \omega_{j_{i}}) \mu_{(t_{1}, \dots, t_{k}), \bfj_{k, 2}}(\prod_{i \notin J_{k'}}h_{i}\prod_{i \in J_{k'}}\bar{f}_{i}  \circ \pi)\\
& & +  \sum_{\bfj_{k, 3} \in \cJ^{k}_{3}} \pi_{*}\mu_{t_{1}, \dots, t_{k}}(\prod_{\bfj_{k, 3}}  \omega_{j_{i}}) \mu_{(t_{1}, \dots, t_{k}), \bfj_{k, 3}}(\prod_{i \notin J_{k'}}h_{i}\prod_{i \in J_{k'}}\bar{f}_{i}  \circ \pi).
\eeqnas

Now we estimate the terms in $\cJ^{k}_{1}$. As in \cite{Kim1}, for each $j \in \cJ_{int}$ define the bijective Lipschitz function $\theta_{j}: \bbT^{d} \rightarrow \pi^{-1}(z_{j})$, $\theta_{j}(\bfb) = (\iota(z_{j}), \zero)(1, \bfb)\hat{\Gamma}$. Let $h_{i, j} = h_{i} \circ \theta_{j}$ and $(\bar{f}_{i} \circ \pi)_{j} = \bar{f}_{i} \circ  \pi \circ \theta_{j}$, then $h_{i, j} \in \cC_{c}^{\infty}(\bbT^{d})$, $(\bar{f}_{i} \circ \pi)_{j}  \in \cC_{c}^{\infty}(\bbT^{d})$. For any $g \in \supp \ \omega_{j} \circ \pi_{X}$, we have
\beqnas
d_{G}(g,  \iota(z_{j})) < C^{3}_{1}r, 
\eeqnas
and for $y \in \supp \ \omega_{j} \circ \pi$, 
\beqnas
|h_{i}(y) - h_{i, j}(\sigma(y))| &\ll r\cS(f_{i}),\\
|\bar{f}_{i}  \circ \pi (y) - (\bar{f}_{i} \circ \pi)_{j}(\sigma(y))| &\ll r\cS(f_{i}).
\eeqnas
Thus for $\bfj_{k, 1} \in \cJ^{k}_{1}$, we have
\beqna\label{ineq.4.5.1}
|\mu_{t_{1}, \dots t_{k}, \bfj_{k, 1}}(\prod_{i \notin J_{k'}}h_{i}\prod_{i \in J_{k'}}\bar{f}_{i}\circ \pi) - \mu_{t_{1}, \dots t_{k}, \bfj_{k, 1}}(\prod_{i \notin J_{k'}}h_{i, j_{i}} \circ \sigma \prod_{i \in J_{k'}}(\bar{f}_{i} \circ \pi)_{j_{i}} \circ \sigma)| \ll r^{k}\prod^{k}_{i=1}\cS(f_{i}).
\eeqna

Now we turn to estimate 
$$ 
 \nu_{t_{1}, \dots t_{k}, \bfj_{k, 1}}(\prod_{i \notin J_{k'}}h_{i, j_{i}} \prod_{i \in J_{k'}}(\bar{f}_{i} \circ \pi)_{j_{i}}) := \mu_{t_{1}, \dots t_{k},  \bfj_{k, 1}}(\prod_{i \notin J_{k'}}h_{i, j_{i}} \circ \sigma \prod_{i \in J_{k'}}(\bar{f}_{i} \circ \pi)_{j_{i}} \circ \sigma).
$$

Write $h_{i, j_{i}}$, $(\bar{f}_{i} \circ \pi)_{j_{i}}$ in their Fourier expansions
\beqnas
h_{i, j_{i}}(\bfb) &=& \sum_{\bfm_{i} \in \bbZ^{d}}\hat{h}_{i, j_{i}}(\bfm_{i})e^{2\pi i \bfm_{i} \cdot \bfb}, \\
(\bar{f}_{i} \circ \pi)_{j_{i}}(\bfb) &=& \sum_{\bfm_{i} \in \bbZ^{d}}\widehat{(\bar{f}_{i} \circ \pi)}_{j_{i}}(\bfm_{i})e^{2\pi i \bfm_{i} \cdot \bfb},
\eeqnas
then we derive that
\beqnas
& &| \nu_{t_{1}, \dots t_{k}, \bfj_{k,1}}(\prod_{i \notin J_{k'}}h_{i, j_{i}} \prod_{i \in J_{k'}}(\bar{f}_{i} \circ \pi)_{j_{i}})|\\
&=& | \sum_{(\bfm_{1}, \dots, \bfm_{k}) \in (\bbZ^{d} \setminus \{0\})^{k}}\hat{\nu}_{t_{1}, \dots t_{k} , \bfj_{k, 1}}(-\bfm_{1}, \dots, -\bfm_{k})\prod_{i \notin J_{k'}}\hat{h}_{i, j_{i}}(\bfm_{i})\prod_{i \in J_{k'}}\widehat{(\bar{f}_{i} \circ \pi)}_{j_{i}}(\bfm_{i})|\\
&\leq& \sum_{0<\|\overline{\bfm}_{k} \| < \rho^{-2\kappa_{5}}}|\hat{\nu}_{t_{1}, \dots t_{k} , \bfj_{k, 1}}(-\bfm_{1}, \dots, -\bfm_{k})|\prod_{i \notin J_{k'}}|\hat{h}_{i, j_{i}}(\bfm_{i})|\prod_{i \in J_{k'}}|\widehat{(\bar{f}_{i} \circ \pi)}_{j_{i}}(\bfm_{i})| \\
& & \qquad + \sum_{\|\overline{\bfm}_{k}\| \geq \rho^{-2\kappa_{5}}}|\hat{\nu}_{t_{1}, \dots t_{k} , \bfj_{k, 1}}(-\bfm_{1}, \dots, -\bfm_{k})|\prod_{i \notin J_{k'}}|\hat{h}_{i, j_{i}}(\bfm_{i})|\prod_{i \in J_{k'}}|\widehat{(\bar{f}_{i} \circ \pi)}_{j_{i}}(\bfm_{i})|\\
&:=& I_{1} + I_{2}, 
 \eeqnas
 where $\overline{\bfm}_{k} = (\bfm_{1}, \dots, \bfm_{k})$. 

As we mentioned before, the assumptions of Proposition~\ref{prop.k} are satisfied. Notice that  $\rho^{-2c_{5}} = \rho^{- \frac{2c^{2}_{4}}{3dk}} \leq \rho^{-c_{4}}$, then by \eqref{est.key.3} we obtain
(choosing $r = \rho^{c^{2}_{4}} = \rho^{3dkc_{5}}$, $c_{5} = \frac{c^{2}_{4}}{3dk}$) 
\beqna\label{est.int.1}
\nonumber I_{1} &\leq& \sum_{0<\|\overline{\bfm}_{k}\| < \rho^{-2c_{5}}} |\hat{ \nu}_{t_{1}, \dots t_{k} , \bfj_{k}}(\bfm_{1}, \dots, \bfm_{k})|\prod^{k}_{i=1}\cS(f_{i})\\
&\ll& \rho^{-2dkc_{5}}\rho^{3dkc_{5}}\prod^{k}_{i=1}\cS(f_{i}) =  \rho^{dkc_{5}}\prod^{k}_{i=1}\cS(f_{i}).
\eeqna
%
Notice that
\beqnas
 \widehat{D^{\bdalp}h_{i, j_{i}}}(\bfm_{i}) &=& (2\pi i \bfm_{i})^{\bdalp}\hat{h}_{i, j_{i}}(\bfm_{i}), \\
 \widehat{D^{\bdalp}(\bar{f}_{i} \circ \pi)_{j_{i}}}(\bfm_{i}) &=& (2\pi i \bfm_{i})^{\bdalp}\widehat{(\bar{f}_{i} \circ \pi)}_{j_{i}}(\bfm_{i}).
 \eeqnas
By \eqref{norm.1}, for $|\bdalp| \leq d+2$, $f \in \cC_{c}^{\infty}(\bbT^{d})$ we have
 \beqnas
\sup_{\bfm \in \bbZ^{d}} | \widehat{D^{\bdalp}f}(\bfm)| \leq \|D^{\bdalp}f\|_{L^{\infty}(\bbT^{d})} \leq \cS^{\bbT^{d}}(f).
 \eeqnas
 Thus we derive that for $|\bdalp| = d+2$, 
 \beqnas
  |\hat{h}_{i, j_{i}}(\bfm_{i})| &\ll& |\bfm_{i}^{-\bdalp}|\cS^{\bbT^{d}}(h_{i, j_{i}}) \ll |\bfm_{i}^{-\bdalp}|\htt(z_{j_{i}})^{l_{0}(d-1)}\cS(f_{i}), \\
 |\widehat{(\bar{f}_{i} \circ \pi)_{j_{i}}}(\bfm_{i})| &\ll& |\bfm_{i}^{-\bdalp}|\cS^{\bbT^{d}}((\bar{f}_{i} \circ \pi)_{j_{i}}) \ll |\bfm_{i}^{-\bdalp}|\htt(z_{j_{i}})^{l_{0} (d-1)}\cS(f_{i}),
 \eeqnas
 implying that
 \beqnas
\prod_{i \notin J_{k'}}|\hat{h}_{i, j_{i}}(\bfm_{i})|\prod_{i \in J_{k'}}|\widehat{(\bar{f}_{i} \circ \pi)}_{j_{i}}(\bfm_{i})| \leq \prod^{k}_{i=1}(|\bfm_{i}^{-\bdalp}|\htt(z_{j_{i}})^{l_{0}(d-1)}\cS(f_{i})).
 \eeqnas
Thus we deduce that 
\beqna\label{est.int.2}
\nonumber I_{2} &\leq& \sum_{\|\overline{\bfm}_{k}\| \geq \rho^{-2c_{5}}}\prod^{k}_{i=1}\big(|\bfm_{i}^{-\bdalp}|\htt(z_{j_{i}})^{l_{0}(d-1)}\cS(f_{i}) \big)\\
&\leq& \rho^{4kc_{5}}\prod^{k}_{i=1}\big(\htt(z_{j_{i}})^{l_{0}(d-1)}\cS(f_{i}) \big) \ll  \rho^{3kc_{5}}\prod^{k}_{i=1}\cS(f_{i}),
\eeqna
where the last line is due to $\sum_{\|\overline{\bfm}_{k}\| \geq \rho^{-2c_{5}}}\prod^{k}_{i=1}|\bfm_{i}^{-\bdalp}| \ll \rho^{4kc_{5}}$
and the fact that $\htt(z_{j_{i}}) \leq \epsilon^{-1} = \rho^{-\frac{c_{5}}{2l_{0}(d-1)}}$ for $z_{j_{i}} \in \cJ_{int}$.
Combining \eqref{est.int.1} and \eqref{est.int.2} together, we derive
 \beqnas
|\nu_{t_{1}, \dots t_{k}, \bfj_{k,1}}(\prod_{i \notin J_{k'}}h_{i, j_{i}} \prod_{i \in J_{k'}}(\bar{f}_{i} \circ \pi)_{j_{i}})| = O(\rho^{kc_{5}}\prod^{k}_{i=1}\cS(f_{i})), 
 \eeqnas 
such that with \eqref{ineq.4.5.1} and $\sum_{\bfj_{k, 1} \in \cJ^{k}_{1}} \pi_{*}\mu_{t_{1}, \dots, t_{k}}(\prod_{\bfj_{k, 1}}  \omega_{j_{i}}) =   \pi_{*}\mu_{t_{1}, \dots, t_{k}}\big(\prod^{k}_{i=1}(\sum_{j_{i} \in \cJ_{int}}  \omega_{j_{i}})\big) < 1$, 
 \beqna\label{est.j1}
\nonumber & &\sum_{\bfj_{k, 1} \in \cJ^{k}_{1}} \pi_{*}\mu_{t_{1}, \dots, t_{k}}(\prod_{\bfj_{k, 1}}  \omega_{j_{i}})\mu_{t_{1}, \dots t_{k}, \bfj_{k, 1}}(\prod_{i \notin J_{k'}}h_{i, j_{i}} \prod_{i \in J_{k'}}(\bar{f}_{i} \circ \pi)_{j_{i}}) \\
 &=& O(r^{k} + \rho^{kc_{5}})\prod^{k}_{i=1}\cS(f_{i}) = O(\rho^{\frac{c^{2}_{4}}{3d}})\prod^{k}_{i=1}\cS(f_{i}).
 \eeqna
 
Now we turn to the terms in $\cJ^{k}_{2}$.  By \eqref{int.q}, $|\cJ^{k}_{2}| \leq \max\{r^{c_{3}}, \epsilon^{d}\}^{k}r^{-k(d^{2}-1)}$. Then, 
\beqna\label{est.j2}
\nonumber & &\sum_{\bfj_{k, 2} \in \cJ^{k}_{2}} \pi_{*}\mu_{t_{1}, \dots, t_{k}}(\prod_{\bfj_{k, 2}}  \omega_{j_{i}}) \mu_{(t_{1}, \dots, t_{k}), \bfj_{k, 2}}(\prod_{i \notin J_{k'}}h_{i}\prod_{i \in J_{k'}}\bar{f}_{i}  \circ \pi)\\
&\leq& \max\{r^{c_{3}}, \epsilon^{d}\}^{k}\prod^{k}_{i=1}\cS(f_{i}) \leq  \max\{\rho^{c_{3}c^{2}_{4}}, \rho^{\frac{c^{2}_{4}}{6l_{0}(d-1)}}\}\prod^{k}_{i=1}\cS(f_{i}).
\eeqna

Let $\cJ^{k}_{3, l}$ denote the subset in $\cJ^{k}_{3}$ which has $l$ elements in $\cJ_{int}$, then $\cJ^{k}_{3} = \bigcup^{k-1}_{l=1}\cJ^{k}_{3, l}$. Let $\ell = \{\ell_{1}, \dots, \ell_{l}\} \subset \{1, \dots, k\}$ with $|\ell| = l$, and $\cJ^{k, \ell}_{3, l} \subset \cJ^{k}_{3, l}$ be the subset with the elements of $\ell$ belong to $\cJ_{int}$ and other $k-l$ elements in $\cJ \setminus \cJ_{int}$. For $\bfj_{k, 3} \in \cJ^{k, \ell}_{3, l}$, write $\bfj_{k, 3} = (\bfj^{\ell}_{l}, \bfj^{\ell}_{k-l})$, where $\bfj^{\ell}_{l} \in \cJ^{l}_{int}$ and $\bfj^{\ell}_{k-l} \in (\cJ \setminus \cJ_{int})^{k-l}$.
\beqnas
& &\sum_{\cJ^{k, \ell}_{3, l}}\pi_{*}\mu_{t_{1}, \dots, t_{k}}(\prod_{j \in \bfj_{k, 3}}  \omega_{j_{i}})\mu_{(t_{1}, \dots, t_{k}), \bfj_{k, 3}}\big(\prod_{i \notin J_{k'}}h_{i}\prod_{i \in J_{k'}}\bar{f}_{i}  \circ \pi \big) \\
&=&\sum_{\bfj^{\ell}_{l} \in \cJ^{l}_{int}}\sum_{\bfj^{\ell}_{k-l} \in (\cJ \setminus \cJ_{int})^{k-l}}\int_{Y}\big(\prod_{i \notin J_{k'}, j_{i} \in \cJ_{int}}h_{i}(y_{i})\omega_{j_{i}}(\pi(y_{i}))\prod_{i \in J_{k'}, j_{i} \in \cJ_{int}}\bar{f}_{i}  \circ \pi(y_{i})\omega_{j_{i}}(\pi(y_{i}))\big) \\
& & \cdot \big(\prod_{i \notin J_{k'},  j_{i} \in \cJ \setminus \cJ_{int}}h_{i}(y_{i})\omega_{j_{i}}(\pi(y_{i}))\prod_{i \in J_{k'},  j_{i} \in \cJ \setminus \cJ_{int}}\bar{f}_{i}  \circ \pi(y_{i})\omega_{j_{i}}(\pi(y_{i})) \big)\mu_{t_{1}, \dots, t_{k}}(y_{1}, \dots, y_{k}) \\
&\leq&\sum_{\bfj^{\ell}_{l} \in \cJ^{l}_{int}}\int_{Y}\prod_{i \notin J_{k'}, j_{i} \in \cJ_{int}}|h_{i}(y_{i})|\omega_{j_{i}}(\pi(y_{i}))\prod_{i \in J_{k'}, j_{i} \in \cJ_{int}}|\bar{f}_{i}  \circ \pi|(y_{i})\omega_{j_{i}}(\pi(y_{i})) \mu_{t_{\ell_{1}}, \dots, t_{\ell_{l}}}(y_{\ell_{1}}, \dots, y_{\ell_{l}}) \\
& & \cdot \prod_{i \notin J_{k'}, j_{i} \in \cJ \setminus \cJ_{int}}\|h_{i}\|_{L^{\infty}}\prod_{i \in J_{k'}, j_{i} \in \cJ \setminus \cJ_{int}}\|\bar{f}_{i}  \circ \pi\|_{\infty}.
\eeqnas 
By \eqref{est.j1}, we have
\beqnas
& &\sum_{\bfj^{\ell}_{l} \in \cJ^{l}_{int}}\int_{Y} \prod_{i \notin J_{k'}, j_{i} \in \cJ_{int}}|h_{i}|(y_{i})\omega_{j_{i}}(\pi(y_{i}))\prod_{i \in J_{k'}, j_{i} \in \cJ_{int}}|\bar{f}_{i}  \circ \pi|(y_{i})\omega_{j_{i}}(\pi(y_{i})) \mu_{t_{\ell_{1}}, \dots, t_{\ell_{l}}}(y_{\ell_{1}}, \dots, y_{\ell_{l}}) \\
&=& O(\rho^{\frac{c^{2}_{4}}{3d}})\prod_{i \notin J_{k'},  j_{i} \in \cJ_{int}}\cS(h_{i})\prod_{i \in J_{k'},  j_{i} \in \cJ_{int}}\cS(\bar{f}_{i}  \circ \pi),
\eeqnas
which implies that
\beqnas
\sum_{\cJ^{k, \ell}_{3, l}}\pi_{*}\mu_{t_{1}, \dots, t_{k}}(\prod_{j \in \bfj_{k, 3}}  \omega_{j_{i}})\mu_{(t_{1}, \dots, t_{k}), \bfj_{k, 3}}\big(\prod_{i \notin J_{k'}}h_{i}\prod_{i \in J_{k'}}\bar{f}_{i}  \circ \pi \big) = O(\rho^{\frac{c^{2}_{4}}{3d}})\prod^{k}_{i=1}\cS(f_{i}).
\eeqnas
Thus we derive that
\beqna\label{est.j3}
\sum_{\cJ^{k}_{3}}\pi_{*}\mu_{t_{1}, \dots, t_{k}}(\prod_{j \in \bfj_{k, 3}}  \omega_{j_{i}})\mu_{(t_{1}, \dots, t_{k}), \bfj_{k, 3}}\big(\prod_{i \notin J_{k'}}h_{i}\prod_{i \in J_{k'}}\bar{f}_{i}  \circ \pi \big) = O_{k}(\rho^{\frac{c^{2}_{4}}{3d}})\prod^{k}_{i=1}\cS(f_{i}).
\eeqna

 Combining \eqref{est.j1}, \eqref{est.j2} and  \eqref{est.j3}, we obtain
 \beqnas
 \int_{V}\prod_{i \notin J_{k'}}h_{i}(g_{t_{i}}uy_{0}) \prod_{i \in J_{k'}}\bar{f}_{i}(\pi(g_{t_{i}}uy_{0})) dm_{V}(u) = O_{k}(\max\{\rho^{c_{3}c^{2}_{4}}, \rho^{\frac{c^{2}_{4}}{6l_{0}(d-1)}}\} +  \rho^{\frac{c^{2}_{4}}{3d}} \big)\prod^{k}_{i=1}\cS(f_{i}), 
 \eeqnas
which leads to 
 \beqna\label{est.mu}
 \mu_{t_{1}, \dots, t_{k}}(\prod^{k}_{i =1}f_{i}) = \prod^{k}_{i =1}\int_{Y}f_{i}dm_{Y} + O_{k}(e^{-\lambda' D(t_{1}, \dots, t_{k})} + \max\{\rho^{c_{3}c^{2}_{4}}, \rho^{\frac{c^{2}_{4}}{6l_{0}(d-1)}}\} +  \rho^{\frac{c^{2}_{4}}{3d}})\prod^{k}_{i =1}\cS(f_{i})
 \eeqna
 with \eqref{eq.1} and \eqref{eq.2}. Then we can choose $0 < \delta' < \delta_{2}$, which is independent of $k$, such that
 \beqnas
 e^{-\lambda' D(t_{1}, \dots, t_{k})} + \max\{\rho^{c_{3}c^{2}_{4}},  \rho^{\frac{c^{2}_{4}}{6l_{0}(d-1)}}\} +  \rho^{\frac{c^{2}_{4}}{3d}} \ll \zeta(\sigma(y_{0}), e^{\frac{\alpha_{0}D(t_{1}, \dots, t_{k})}{2}})^{-\delta'}.
 \eeqnas
 Thus the proof is completed. 
 \epf

\subsection{Proof of Corollary~\ref{thm.3}}
\bpf
Fix $\kappa > d+1$. Let $T_{\epsilon} = \bigcup^{\infty}_{c \geq \epsilon}D(\kappa, c)$. Then by \eqref{est.m}, we have 
$$
m_{\bbT^{d}} \big(\bbT^{d} \setminus T_{\epsilon} \big) \ll \epsilon. 
$$
Moreover, by \eqref{est.z} we have for $\bfb \in T_{\epsilon}$, 
$$
\zeta(\bfb, e^{\frac{\alpha_{0}D(t_{1}, \dots, t_{k})}{2}})^{-1} \leq \epsilon^{-\frac{1}{\kappa+1}} e^{-\frac{\alpha_{0}D(t_{1}, \dots, t_{k})}{2(\kappa+1)}}.
$$
Thus by \eqref{effe.eqdis.2}, for $\|g_{0}\| \leq \epsilon^{\frac{\delta'}{\kappa+1}} e^{\frac{\delta'\alpha_{0}D(t_{1}, \dots, t_{k})}{2(\kappa+1)}}$, 
\beqnas
& &\int_{\bbT^{d}}\int_{V}\prod^{k}_{i =1}f_{i}(g_{t_{i}}uy(\bfb)) dm_{V}(u)dm_{\bbT^{d}}(\bfb) -  \prod^{k}_{i =1}\int_{Y}f_{i} dm_{Y} \\
&=& \int_{T_{\epsilon}}\int_{V}\prod^{k}_{i =1}f_{i}(g_{t_{i}}uy(\bfb)) dm_{V}(u)dm_{\bbT^{d}}(\bfb) +  \int_{\bbT^{d} \setminus T_{\epsilon}}\int_{V}\prod^{k}_{i =1}f_{i}(g_{t_{i}}uy(\bfb)) dm_{V}(u)dm_{\bbT^{d}}(\bfb)  -  \prod^{k}_{i =1}\int_{Y}f_{i} dm_{Y}\\
&\ll&  \prod^{k}_{i =1}\cS(f_{i}) \big(\epsilon^{-\frac{\delta'}{\kappa+1}} e^{-\frac{\alpha_{0} D(t_{1}, \dots, t_{k})\delta'}{2(\kappa+1)}} + \epsilon \big).
\eeqnas
By taking $\epsilon = e^{-\frac{\alpha_{0}D(t_{1}, \dots, t_{k})\delta'}{2(\delta' + \kappa+1)}}$, we prove \eqref{effe.eqdis.3} with $\delta_{\kappa} = \frac{\delta'}{2(\delta' + \kappa+1)}$.
\epf

\section{The central limit theorem for smooth Siegel transforms}\label{sec.3}
We adopt the framework in \cite{BG19} and formulate the counting function $\Delta_{T}(u, \bfx)$ into the Siegel transform on the affine lattice space $Y$. For $u \in M_{m, n}(\bbT)$ and $\bfx \in \bbT^{m}$, define
\beqnas
\Lambda(u, \bfx) = \{(p_{1} + \sum^{n}_{j=1}u_{1j}q_{j} + x_{1} , \dots, p_{m} + \sum^{n}_{j=1}u_{mj}q_{j} + x_{m}, \bfq), \ (\bfp, \bfq) \in \bbZ^{m} \times \bbZ^{n}\}.
\eeqnas
For $T > 1$, let
\beqnas
\Omega_{T} &=& \{(\bfx, \bfy) \in \bbR^{m} \times \bbR^{n}, \ 1 \leq \|\bfy\| < T, \ |x_{i}| < \theta_{i}\|\bfy\|^{-\omega_{i}}, \ i =1, \dots, m\}.
\eeqnas
Then we see that 
$$
\Delta_{T}(u, \bfx) = \sharp \{ \Lambda(u, \bfx) \cap \Omega_{T} \} + O(1),
$$ 
where the constant $O(1)$ depends on the norm $\|\cdot\|$.
Define
\beqnas
A_{s} = \diag\{e^{s\omega_{1}}, \dots, e^{s\omega_{m}}, e^{-s}, \dots, e^{-s}\},
\eeqnas
where $\omega_{i} > 0$, $\sum^{m}_{i=1}\omega_{i} = n$. Let $\omega_{0} = \min_{1 \leq i \leq n}\{\omega_{i}, 1\}$. Note that
\beqnas
A_{-s}\Omega_{e} = \{(\bfx, \bfy) \in \bbR^{m} \times \bbR^{n}, \ e^{s} \leq \|\bfy\| < e^{s+1}, \ |x_{i}| < \theta_{i}\|\bfy\|^{-\omega_{i}}, \ i =1, \dots, m\},
\eeqnas
and
\beqnas
\sharp \{\Lambda(u, \bfx) \cap \Omega_{e^{N}}\} =  \sum^{N-1}_{s=0} \sharp\{\Lambda(u, \bfx) \cap A_{-s}\Omega_{e}\}.
\eeqnas
Recall that for a function  $f$ on $\bbR^{d}$ with compact support, its Siegel transform $\hat{f}$ on the affine lattice space $Y$ is defined by
$$
\hat{f}(y) = \sum_{\bfz \in g\bbZ^{d}+ g\bfb}f(\bfz),
$$
where $y = (g, \zero)(1, \bfb)\hat{\Gamma} \in Y$, $g \in \SL(d, \bbR)$ and $\bfb \in \bbT^{d}$. Also notice that $\Lambda(u, \bfx)$ can be regarded as an affine lattice
\beqnas
\Lambda(u, \bfx) = \left(\left(\begin{array}{ccc}
\un_{m} & u \\
0 & \un_{n} \\
\end{array}
\right), \left(\begin{array}{ccc}
\bfx \\
\mathbf{0}  \\
\end{array}
\right) \right)
\hat{\Gamma}.
\eeqnas
Let $\chi$ be the characteristic function of $\Omega_{e}$, $\hat{\chi}$ be its Siegel transform, so that 
$$
\sharp\{\Lambda(u, \bfx) \cap A_{-s}\Omega_{e}\} = \hat{\chi}(A_{s}\Lambda(u, \bfx)),
$$ 
which enables us to turn the central limit theorem of $\Delta_{T}(u, \bfx)$ into that of the flow $\sum^{[\log T]-1}_{s=0}\hat{\chi}(A_{s}\Lambda(u, \bfx))$ on the affine lattice space $Y$. 

Let
$$
y(\bfx) = (\Id_{d}, \zero)(1,  \left(\begin{array}{ccc}
\bfx \\
\mathbf{0}  \\
\end{array}
\right))\hat{\Gamma} \in Y,
$$
where we set $g_{0} = \Id_{d}$, $\bfx \in \bbT^{m}$, and let  $V \subset H$ given by
\beqnas\label{def.h}
V = \{M | M = \begin{pmatrix}
\Id_{m}  & u\\
0 & \Id_{n}
\end{pmatrix}, 
u \in M_{m, n}([0, 1])\}.
\eeqnas
Notice that every $u \in V$ can be regarded as a matrix $u \in  M_{m, n}([0, 1])$. 
Define $\cV = V \times \bbT^{m}$ and let $m_{\cV} = m_{V} \times m_{\bbT^{m}}$ be the probability measure on $\cV$. Then for $v = (u, \bfx) \in \cV$, we have
\beqnas
\Lambda(v) := uy(\bfx) = \left(\left(\begin{array}{ccc}
\un_{m} & u \\
0 & \un_{n} \\
\end{array}
\right), \left(\begin{array}{ccc}
\bfx \\
\mathbf{0}  \\
\end{array}
\right) \right)
\hat{\Gamma} = \Lambda(u, \bfx).
\eeqnas

In this section we prove the central limit theorem for smooth Siegel transforms. 
\bthm\label{clt.sg}
Let $m \geq 2$. Assume that $v = (u, \bfx)$ is uniformly distributed on $\cV$. For $f \in \cC^{\infty}_{c}(Y)$, $f \geq 0$ and $\supp f \subset \{(x_{m+1}, \dots, x_{m+n}) \neq 0\}$, we have
\beqnas
\frac{1}{\sqrt{N}}\sum^{N-1}_{s=0}\big( \hat{f}(A_{s}\Lambda(v)) - \int_{\cV}\hat{f}(A_{s}\Lambda(v))dm_{\cV}(v) \big) \Longrightarrow \cN(0, \sigma^{2}), 
\eeqnas
as $N \rightarrow \infty$, and 
\beqna\label{var.aff.f.1}
\sigma^{2} 
&=& \sum^{+\infty}_{s = -\infty}\int_{\bbR^{d}}f(A_{s}x)f(x)dx.
\eeqna
\ethm

\subsection{Preliminaries}

To prove Theorem~\ref{clt.sg}, we need some estimates on the Siegel transform $\hat{f}$. Notice that for $y = (g, \zero)(1, \bfb)\hat{\Gamma}$, we have
\beqna\label{shf.af}
\hat{f}(y) = \sum_{\bfz \in g\bbZ^{d}}f(\bfz +  g\bfb) = \sum_{\bfz \in g\bbZ^{d}}\tau_{g\bfb}f(\bfz) = \widehat{\tau_{g\bfb}f}(\pi(y)), 
\eeqna
where $\tau_{g\bfb}f(\bfz) := f(\bfz + g\bfb)$. 
Recall that the Siegel transform on $X$ is controlled by an explicit function $\alpha$ on $X$, i.e., for $f \in \cC_{c}(\bbR^{d})$, 
\beqnas
|\hat{f}(x)| \ll O_{\supp f}(\|f\|_{L^{\infty}}\alpha(x))
\eeqnas
holds for any $x \in X$. Then for $y \in Y$, 
\beqna\label{est.sg.f}
|\hat{f}(y)| = |\widehat{\tau_{g\bfb}f}(\pi(y))| \ll O_{\supp \tau_{g\bfb}f}(\|\tau_{g\bfb}f\|_{L^{\infty}}\alpha(\pi(y))) = O_{\supp f}(\|f\|_{L^{\infty}}\alpha(\pi(y))).
\eeqna
Note that $\alpha \in L^{p}(X)$ for $1 \leq p < d$,  and
\beqna\label{est.alp}
m_{X}(\alpha \geq L) \ll L^{-p}.
\eeqna 
Recall the Rogers formula for the Siegel transform, see Proposition 5.1, \cite{DFV}. For $f \in \cC^{\infty}_{c}(\bbR^{d})$, one has
\beqna\label{rog.1}
\int_{Y}\hat{f}(y)dm_{Y}(y) = \int_{\bbR^{d}}f(x)dx, 
\eeqna
and 
\beqna\label{rog.2}
\int_{Y}|\hat{f}(y)|^{2}dm_{Y}(y) = \left( \int_{\bbR^{d}}f(x)dx \right)^{2} + \int_{\bbR^{d}}f^{2}(x)dx. 
\eeqna
By Proposition 4.5 in \cite{BG19}, there exists $c > 0$ such that for $L \geq 1$ and $s \geq c \log L$
\beqna\label{est.alp.1}
m_{V}(\alpha(A_{s}u) \geq L) = O(L^{-p}),
\eeqna
for all $p < d$.

We have the following estimates with slight modifications of the proofs of Proposition 4.6, Proposition 4.8 in \cite{BG19}.  
\blem
Let $f$ be a bounded measurable function on $\bbR^{d}$ with compact support in $\{(x_{m+1}, \dots, x_{m+n}) \neq 0\}$. Then for $\bfx \in \bbT^{m}$ fixed, 
\beqna\label{est.sg.1}
\sup_{s \geq 0}\int_{V}|\hat{f}(A_{s}uy(\bfx))|dm_{V}(u) < \infty,
\eeqna
and 
\beqna\label{est.sg.2}
\sup_{s \geq 0}(1+s)^{-v_{m}}\int_{V}|\hat{f}(A_{s}uy(\bfx))|^{2}dm_{V}(u) < \infty,
\eeqna
where $v_{1} = 1$ and $v_{m}=0$ when $m \geq 2$. 
\elem

\bpf
Without loss of generality, we assume $f$ to be the characteristic function of the set
\beqnas
 \{(\bfx, \bfy) \in \bbR^{m} \times \bbR^{n}, \ v_{1} \leq \|\bfy\| \leq v_{2}, \ |x_{i}| < \theta \|\bfy\|^{-\omega_{i}}, \ i =1, \dots, m\}.
\eeqnas
for $0 < v_{1} < v_{2}$ and $\theta > 0$. 
Then, 
\beqnas
\int_{V}|\hat{f}(A_{s}uy(\bfx))|dm_{V}(u) = \sum_{ \|\bfq \| \in [e^{s}v_{1}, e^{s}v_{2}]}\prod^{m}_{i=1}\big(\sum_{p_{i} \in \bbZ}\int_{[0, 1]^{n}}\un_{[-\theta \|\bfq\|^{-\omega_{i}}, \theta \|\bfq\|^{-\omega_{i}}]}(p_{i} + \langle \bfu_{i}, \bfq \rangle + x_{i})d\bfu_{i}\big).
\eeqnas
Notice that
\beqnas
\int_{[0, 1]^{n}}\un_{[-\theta \|\bfq\|^{-\omega_{i}}, \theta \|\bfq\|^{-\omega_{i}}]}(p_{i} + \langle \bfu_{i}, \bfq \rangle + x_{i})d\bfu_{i} \ll \theta \|\bfq\|^{-\omega_{i}-1}, 
\eeqnas
and to have a non-trivial integral we need $|p_{i}| \leq O(\|\bfq\|)$. Thus, we derive that
\beqnas
\int_{V}|\hat{f}(A_{s}uy(\bfx))|dm_{V}(u) \ll \sum_{ \|\bfq \| \in [e^{s}v_{1}, e^{s}v_{2}]}\|\bfq\|^{-n} < \infty.
\eeqnas

The estimate \eqref{est.sg.2} follows from a similar modification of the proof of Proposition 4.8, \cite{BG19}. 
\epf

Recall the cut-off function $\eta_{L} \in \cC^{\infty}_{c}(X)$ on $X$: for any $c > 1$, 
\beqna\label{est.eta}
0 \leq \eta_{L} \leq 1, \ \eta_{L} = 1 \ on \  \{\alpha \leq c^{-1}L\},  \eta_{L} = 0 \ on \ \{\alpha >  cL\}, \ \|\eta_{L}\|_{\cC^{k}} \ll 1. 
\eeqna
To control the norm $\cS(\hat{f})$ for $f \in \cC_{c}(\bbR^{d})$, we construct a cut-off function on $X$ following Proposition 4.5 and Lemma 4.11 in \cite{BG19}. 
\blem
For any $c >1$, there exists a family of functions $\{\beta_{L}\} \in \cC^{\infty}_{c}(X)$ satisfying $0 \leq \beta_{L} \leq 1$, 
\beqna\label{beta.1}
\beta_{L} = 1 \ on \ \{\htt \leq c^{-1}L\}, \  \beta_{L} = 0 \ on \ \{\htt > cL\}, \ \|\beta_{L}\|_{\cC^{k}} \ll 1,
\eeqna
and
\beqna\label{beta}
\int_{X}\beta_{L}(x)dm_{X}(x) = m_{X}(X_{L}), 
\eeqna
where $X_{L} =  \{x \in X, \htt(x) \leq L\}$. Moreover, for $L \geq 1$ and $s \geq c'\log L$ where $c'$ is a constant, we have 
\beqna\label{est.ht}
m_{V}(\htt(A_{s}u) \geq L) \ll L^{-d}. 
\eeqna
\elem

\bpf
Let $\un_{L}$ be the characteristic function of the set $X_{L} =  \{x \in X, \htt(x) \leq L\}$. Let $\phi \in \cC^{\infty}_{c}(G)$ be a non-negative function with $\int_{G}\phi dm_{G} =1$, with compact support small enough in the neighborhood of the identity in $G$, such that for all $g \in \supp \ \phi$ and $x \in X$,  
\beqna\label{est.2.9}
c^{-1}\htt(x) \leq \htt(g^{-1}x) \leq c \ \htt(x).
\eeqna
This is ensured by \eqref{est.ht.1} and the fact that $\|g\| = \|g^{-1}\|$. For $x \in X$, define 
\beqnas
\beta_{L}(x) = (\phi \ast \un_{L})(x) = \int_{G} \phi(g)\un_{L}(g^{-1}x)dm_{G}(g).
\eeqnas
Then it is easy to see that $0 \leq \beta_{L} \leq 1$. For $x \in \{\htt \leq c^{-1}L\}$, we have $g^{- 1}x \in  X_{L}$ for $g \in \supp \ \phi$ by \eqref{est.2.9}, leading to $\beta_{L}(x) = 1$. On the other hand, for  $x \in \{\htt > c L\}$, we have $g^{- 1}x \in  X_{L}^{c}$  for $g \in \supp \ \phi$, such that $\beta_{L}(x) = 0$.
For any differential operator $\cD_{Z}$, we have $\|\beta_{L}\|_{\cC^{k}} \ll \sum_{\deg (Z) \leq k}\| \cD_{Z}\phi \|_{L^{1}} \ll 1$.
The invariance of the Haar measure $m_{X}$ leads to
\beqnas
\int_{X}\beta_{L}(x)dm_{X}(x) = \int_{G} \phi(g)\int_{X} \un_{L}(g^{-1}x)dm_{X}(x)dm_{G}(g) = m_{X}(X_{L}).
\eeqnas
Applying \eqref{exp.m.X} to $\beta_{L}$ and notice that $\cS^{X}(\beta_{L}) = \cS^{X}_{l_{0}}(\beta_{L}) \ll L^{l_{0}}$, we derive that
\beqnas
m_{V}(\htt(A_{s}u) \leq cL) \geq \int_{V}\beta_{L}(A_{s}u)dm_{V}(u) =  m_{X}(X_{L}) + O(L^{l_{0}}e^{-\lambda' s}).
\eeqnas
Moreover, by \eqref{est.ht.2}, 
\beqnas
m_{X}(X_{L}) = 1 - m_{X}(\htt(x) \geq L) = 1 + O(L^{-d}),  
\eeqnas
so that combining the above estimates together yields
\beqnas
m_{V}(\htt(A_{s}u) \geq L) =  O(L^{-d} + L^{l_{0}}e^{-\lambda' s}) = O(L^{-d})
\eeqnas
for $s \geq \frac{d+ l_{0}}{\lambda'}\log L$.

\epf

For $f \in \cC_{c}(\bbR^{d})$, define its truncated Siegel transform as $\hat{f}^{L}(y) =  \hat{f}(y)\eta_{L}(\pi(y))\beta_{L}(\pi(y))$. By \eqref{est.sg.f},  we have the following estimates as analogues of Lemma 4.12, \cite{BG19}.
\blem\label{est.Lsg}
\beqna
\label{est.4.4.1}\|\hat{f}^{L}\|_{L^{p}(Y)} &\leq& \|\hat{f}\|_{L^{p}(Y)} = O_{\supp f}(\|f\|_{L^{\infty}}), \ 1 \leq p < d, \\
\label{est.4.4.2}\|\hat{f}^{L}\|_{L^{\infty}(Y)} &=& O_{\supp f}(\|f\|_{L^{\infty}}L), \\
\label{est.4.4.3}\|\hat{f}^{L}\|_{\cC^{k}} &=& O_{\supp f}(\|f\|_{\cC^{k}}L), \\
\label{est.4.4.5}\|\hat{f} - \hat{f}^{L}\|_{L^{1}(Y)} &=& O_{\supp f, p}\big(\|f\|_{L^{\infty}}\|\alpha\|_{L^{p}}L^{-(p-1)} \big), \ 1 \leq p < d, \\
\label{est.4.4.6}\|\hat{f} - \hat{f}^{L}\|_{L^{2}(Y)} &=& O_{\supp f, p}\big(\|f\|_{L^{\infty}}\|\alpha\|_{L^{p}}L^{-\frac{(p-2)}{2}}  \big), \ 1 \leq p < d.
\eeqna
Moreover, 
\beqna\label{est.4.4.4}
\cS_{k}(\hat{f}^{L}) &=&  O_{\supp f}(\|f\|_{\cC^{k}}L^{k+1}). 
\eeqna
\elem
\bpf
By \eqref{est.sg.f}, we have
\beqnas
\|\hat{f}^{L}\|_{L^{p}(Y)} &=&  \|\hat{f}(y)\eta_{L}(\pi(y))\beta_{L}(\pi(y))\|_{L^{p}(Y)} \leq  \|\hat{f}(y)\|_{L^{p}(Y)}\\
&=& O_{\supp f}(\|f\|_{L^{\infty}}\|\alpha\|_{L^{p}(X)}) \leq  O_{\supp f}(\|f\|_{L^{\infty}}).
\eeqnas
Since $\supp \ \eta_{L} \subset  \{\alpha \leq  cL\}$, we derive
\beqnas
\|\hat{f}^{L}\|_{L^{\infty}(Y)} \leq \|f\|_{L^{\infty}}\|\alpha(\pi(y))\eta_{L}(\pi(y))\|_{L^{\infty}(Y)} = O_{\supp f}(\|f\|_{L^{\infty}}L).
\eeqnas
By \eqref{est.eta}, \eqref{beta.1} and  $\cD_{z}\hat{f} = \widehat{\cD_{z}f}$, 
\beqnas
\|\hat{f}^{L}\|_{\cC^{k}(Y)} &=& \sum^{k}_{l=0}\|\cD_{l}(\hat{f} \cdot \eta_{L}\beta_{L})\|_{L^{\infty}} \ll \sum^{k}_{l=0}\|\widehat{\cD_{l}f}\un_{\{\alpha \leq  cL\}}\|_{L^{\infty}} \\
&\leq&  O_{\supp f}(\sum^{k}_{l=0}\|\cD_{l}f\|_{L^{\infty}}L) =   O_{\supp f}(\|f\|_{\cC^{k}}L).
\eeqnas
For $p < d$, let $q$ satisfy $\frac{1}{p} + \frac{1}{q} = 1$, then
\beqnas
\|\hat{f} - \hat{f}^{L}\|_{L^{1}} &=& \int_{Y}|\hat{f} - \hat{f}^{L}|dm_{Y} =  \int_{Y}|\hat{f}(y)||1 - \eta^{L}(\pi(y))\beta^{L}(\pi(y))|dm_{Y}(y)\\
&\ll& O_{\supp f}\big(\|f\|_{L^{\infty}}\int_{\{\alpha \geq c^{-1}L\} \cup \{\htt \geq c^{-1}L \}}\alpha dm_{X} \big) \\
&\leq& O_{\supp f}\big( \|f\|_{L^{\infty}}\|\alpha\|_{L^{p}}(m_{X}(\{\alpha \geq c^{-1}L\})^{\frac{1}{q}} + m_{X}(\{\htt \geq c^{-1}L\})^{\frac{1}{q}}) \big) \\
&\leq& O_{\supp f}\big( \|f\|_{L^{\infty}}\|\alpha\|_{L^{p}}(L^{-(p-1)} + L^{-\frac{d}{p}(p-1)}) \big), 
\eeqnas
where the last line is due to \eqref{est.alp} and \eqref{est.ht.2}. 
The estimate of $\|\hat{f} - \hat{f}^{L}\|_{L^{2}(Y)}$ follows the same argument. 
By the definition~\eqref{def.norm.Y} and \eqref{est.sg.f}, we have for any $k \in \bbN$, 
\beqnas
\cS_{k}(\hat{f}^{L})^{2} &=&  \sum_{\deg(\hat{Z}) \leq k}\int_{Y}|\htt(\pi(y))^{k}\cD_{\hat{Z}}\hat{f}^{L}(y)|^{2}dm_{Y}(y)\\
&\ll& \sum_{\deg(\hat{Z}) \leq k}\int_{Y}|\htt(\pi(y))^{k}\widehat{\cD_{\hat{Z}}f}(y)\un_{\{\alpha (\pi(y)) \leq  cL\}}\un_{\{\htt(\pi(y)) \leq  cL\}}|^{2}dm_{Y}(y)\\
&=& O_{\supp f}(\sum_{\deg(\hat{Z}) \leq k}\|\cD_{\hat{Z}}f\|^{2}_{L^{\infty}}\int_{Y}|\htt(\pi(y))^{k}\alpha(\pi(y)) \un_{\{\alpha (\pi(y)) \leq  cL\}}\un_{\{\htt(\pi(y)) \leq  cL\}}|^{2}dm_{Y}(y))\\
&=&  O_{\supp f}(\|f\|^{2}_{\cC^{k}}L^{2(k+1)}). 
\eeqnas

\epf

\subsection{Proof of Theorem~\ref{clt.sg}}

\bpf
We approximate $\hat{f}$ by $\hat{f}^{L}$.  Since the estimates (Corollary~\ref{thm.3}, Lemma~\ref{est.Lsg}) on $Y$ are of the same character as those on $X$, the proof of the CLT proceeds as in \cite{BG19}. It therefore suffices to verify the variance formula \eqref{var.aff.f.1}. 

The approximation argument as in \cite{BG19} yields the variance $\sigma^{2}$ given by
\beqnas
\sigma^{2} = \sum^{+\infty}_{s = -\infty}\left(\int_{Y}(\hat{f} \cdot A_{s})\hat{f} dm_{Y} - (\int_{Y}\hat{f} dm_{Y})^{2} \right).
\eeqnas
Notice that by the Rogers formula \eqref{rog.2}, 
\beqnas
\int_{Y}|\hat{f}(y) + \hat{f}(A_{s}y)|^{2}dm_{Y}(y) &=& \big(\int_{\bbR^{d}}f(A_{s}x)dx + \int_{\bbR^{d}}f(x)dx\big)^{2} + \int_{\bbR^{d}}(f(A_{s}x) + f(x))^{2}dx,\\
\int_{Y}|\hat{f}(y) - \hat{f}(A_{s}y)|^{2}dm_{Y}(y) &=& \big(\int_{\bbR^{d}}f(A_{s}x)dx - \int_{\bbR^{d}}f(x)dx\big)^{2} + \int_{\bbR^{d}}(f(A_{s}x) - f(x))^{2}dx,
\eeqnas
which lead to
\beqna\label{rog.3}
\nonumber \int_{Y}\hat{f}(y)\hat{f}(A_{s}y)dm_{Y}(y) &=& \int_{\bbR^{d}}f(A_{s}x)dx\int_{\bbR^{d}}f(x)dx + \int_{\bbR^{d}}f(A_{s}x)f(x)dx\\
&=& (\int_{\bbR^{d}}f(x)dx)^{2} + \int_{\bbR^{d}}f(A_{s}x)f(x)dx.
\eeqna
With \eqref{rog.1}, we obtain \eqref{var.aff.f.1}. 
\epf

\section{Proof of Theorem~\ref{clt.inDA}}\label{sec.4}
\bpf
We first claim the following CLT for $\hat{\chi}(A_{s}\Lambda(v))$. 
\bprop\label{clt.chi}
Assume that $v = (u, \bfx)$ is uniformly distributed on $\cV = V \times \bbT^{m}$. Then, 
\beqnas
\frac{1}{\sqrt{N}}\sum^{N-1}_{s=0}\big( \hat{\chi}(A_{s}\Lambda(v)) -  \int_{\cV}\hat{\chi}(A_{s}\Lambda(v))dm_{\cV}(v)\big) \Longrightarrow \cN(0, \sigma_{m, n}^{2}), 
\eeqnas
as $N \rightarrow \infty$, and 
\beqna\label{var.aff.f}
\sigma_{m, n}^{2} = 2^{m}(\prod^{m}_{i=1}\theta_{i})n\vol_{n}. 
\eeqna
\eprop

Moreover, we have the following estimate for the expectation of $\hat{\chi}(A_{s}\Lambda(v))$. 
\bprop\label{est.epc}
For $s \in \bbN$, we have
\beqna\label{est.exp}
\int_{\cV}\hat{\chi}(A_{s}\Lambda(v))dm_{\cV}(v) = C_{m, n} + O(e^{-s}), 
\eeqna
where $C_{m, n} =\sigma_{m, n}^{2} =2^{m}(\prod^{m}_{i}\theta_{i})n\vol_{n}$. Summing over $s = 0, \dots, N-1$ yields
\beqna\label{asym.e}
\sum^{N-1}_{s=0}\int_{\cV}\hat{\chi}(A_{s}\Lambda(v))dm_{\cV}(v) = C_{m, n}N + O(1).
\eeqna
Define
\beqnas
R_{T}(v) :=  \Delta(u, \bfx) - \sum^{[\log T] -1}_{s=0}\int_{\cV}\hat{\chi}(A_{s}\Lambda(v))dm_{\cV}(v)
\eeqnas
Then
\beqna\label{est.R}
\int_{\cV}R_{T}(v)dm_{\cV}(v) = C_{m, n}\{\log T\} + O(1).
\eeqna
\eprop
\bpf
The proof is an analogue of Lemma 6.3 in \cite{BG19}. First notice that
\beqnas
\int_{\cV}\hat{\chi}(A_{s}\Lambda(v))dm_{\cV}(v) = \sum_{e^{s} \leq \|\bfq\| < e^{s+1}}\prod^{m}_{i=1}\sum_{p_{i} \in \bbZ}\int_{[0, 1]}\int_{[0, 1]^{n}}\un_{[-\frac{\theta_{i}}{\|\bfq\|^{\omega_{i}}}, \frac{\theta_{i}}{\|\bfq\|^{\omega_{i}}}]}(p_{i} + \langle \bfu_{i}, \bfq \rangle + x_{i})d\bfu_{i}dx_{i}. 
\eeqnas
Let $\un_{i}(\bfu) = \un_{[-\frac{\theta_{i}}{\|\bfq\|^{\omega_{i}}}, \frac{\theta_{i}}{\|\bfq\|^{\omega_{i}}}]}\cdot e_{1}(\bfu) =  \un_{[-\frac{\theta_{i}}{\|\bfq\|^{\omega_{i}}}, \frac{\theta_{i}}{\|\bfq\|^{\omega_{i}}}]}(u_{1})$. By defining an affine map $S_{x}: \bbR^{n} \rightarrow \bbR^{n}$, $S_{x}(\bfu) = (\langle \bfu, \bfq \rangle + x, u_{2}, \dots, u_{n})$ with $x \in [0, 1]$,  and noticing that $S_{x}$ preserves $m_{\bbT^{n}}$,  one can prove that 
\beqnas
& &\sum_{p_{i} \in \bbZ}\int_{[0, 1]}\int_{[0, 1]^{n}}\un_{[-\frac{\theta_{i}}{\|\bfq\|^{\omega_{i}}}, \frac{\theta_{i}}{\|\bfq\|^{\omega_{i}}}]}(p_{i} + \langle \bfu_{i}, \bfq \rangle + x_{i})d\bfu_{i}dx_{i} \\
&=& \int_{[0, 1]}\int_{\bbT^{n}}\hat{\un}_{i}(S_{x_{i}}(\bfu_{i}))d\bfu_{i}dx_{i} =  \int_{[0, 1]}\int_{\bbT^{n}}\hat{\un}_{i}(\bfu_{i})d\bfu_{i}dx_{i}\\
&=& \int_{[0, 1]}\int_{\bbR}\un_{[-\frac{\theta_{i}}{\|\bfq\|^{\omega_{i}}}, \frac{\theta_{i}}{\|\bfq\|^{\omega_{i}}}]}(u)dudx_{i} = \frac{2\theta_{i}}{\|\bfq\|^{\omega_{i}}}.
\eeqnas
Thus, 
\beqnas
\int_{\cV}\hat{\chi}(A_{s}\Lambda(v))dm_{\cV}(v) = 2^{m}\prod^{m}_{i}\theta_{i}\sum_{e^{s} \leq \|\bfq\| < e^{s+1}}\|\bfq\|^{-n}. 
\eeqnas
Let $N(r) = \sharp\{\bfz \in \bbZ^{n}, \|\bfz\| \leq r\}$, and it is known that $N(r) = \vol_{n}r^{n} + E(r)$, where $\vol_{n}$ is the Euclidean volume of the unit ball $\{z,  \|z\| \leq 1\}$, and $E(r) = O_{n}(r^{n-1})$. Thus, 
\beqnas
\sum_{e^{s} \leq \|\bfq\| < e^{s+1}}\|\bfq\|^{-n} &=& \int^{e^{s+1}}_{e^{s}}r^{-n}dN(r) = n\vol_{n}\int^{e^{s+1}}_{e^{s}}r^{-1}dr + O(\int^{e^{s+1}}_{e^{s}}r^{-2}dr)\\
&=& n\vol_{n} + O(e^{-s}), 
\eeqnas 
which leads to \eqref{est.exp} and \eqref{asym.e}. Notice that
\beqnas
R_{T}(v) &=& \sharp\{ \Lambda(v) \cap \Omega_{T}\} - \sharp\{ \Lambda(v) \cap \Omega_{e^{\log T}}\} + O(1)\\
&=& \sharp \{(\bfp, \bfq) \in \bbZ^{m} \times (\bbZ^{n} \setminus  \{\zero\}),\  e^{[\log T]} < \|\bfq\| < T, \ \eqref{ineq.di} \ holds, \ i=1, \dots, m\} + O(1), 
\eeqnas
which enables us to estimate $\int_{\cV}R_{T}(v)dm_{\cV}(v)$ in the same approach:
\beqnas
\int_{\cV}R_{T}(v)dm_{\cV}(v) = 2^{m}(\prod^{m}_{i}\theta_{i})\sum_{e^{[\log T]} \leq \|\bfq\| < T}\|\bfq\|^{-n} + O(1)  = C_{m, n}\{\log T\} + O(1). 
\eeqnas
\epf

By Proposition~\ref{clt.chi} and Proposition~\ref{est.epc}, the CLT for $\frac{\Delta_{T}(u, \bfx) - C_{m, n}\log T}{\sqrt{\log T}}$ follows from the CLT for \\$\frac{1}{\sqrt{N}}\sum^{N-1}_{s=0}\big( \hat{\chi}(A_{s}\Lambda(v)) -  \int_{\cV}\hat{\chi}(A_{s}\Lambda(v))dm_{\cV}(v)\big) $. The proof of Theorem~\ref{clt.inDA} is then completed as in \cite{BG19}; we omit the details. 

\epf
\subsection{Proof of Proposition~\ref{clt.chi}}
We introduce the smooth approximations $\{f_{\epsilon}\}$ to $\chi$, satisfying
\beqna\label{est.eps.1}
 \chi \leq f_{\epsilon} \leq 1,\ \|f_{\epsilon}\|_{\cC^{k}} \ll \epsilon^{-k}.
 \eeqna
 Assume that $f_{\epsilon}$ is compactly supported in the $\epsilon-$neighborhood of $\Omega_{e}$, denoted by $\Omega^{\epsilon}_{e}$, 
$$
\Omega^{\epsilon}_{e} = \{(\bfx, \bfy) \in \bbR^{m+n}, 1- \epsilon \leq \|\bfy\| \leq e + \epsilon, \ |x_{i}| < \theta_{i}(\epsilon)\|\bfy\|^{-\omega_{i}}, \ i=1, \dots, m\},
$$
where $\theta_{i}(\epsilon) = \theta_{i} + O(\epsilon)$.

The following estimate  follows the same approach as in Proposition 6.2, \cite{BG19}.
\blem
For any $\bfx \in \bbT^{m}$, we have 
\beqna\label{est.eps}
\int_{V}|\hat{f}_{\epsilon}(A_{s}uy(\bfx)) - \hat{\chi}(A_{s}uy(\bfx))|dm_{V}(u) \ll \epsilon + e^{-s}.
\eeqna
\elem

\noindent
\emph{Proof of Proposition~\ref{clt.chi}.}
The idea is to utilize the cumulant method to the smooth cut-off functions to prove the CLT for $\hat{\chi}$. As the proof resembles that in \cite{BG19}, we only present a sketch.

Define
\beqnas
F_{N}(v) = \frac{1}{\sqrt{N}}\sum^{N-1}_{s=0}\big( \hat{\chi}(A_{s}\Lambda(v)) - \int_{\cV}\hat{\chi}(A_{s}\Lambda(v))dm_{\cV}(v) \big),
\eeqnas
\beqnas
\tilde{F}_{N}(v)= \frac{1}{\sqrt{N}}\sum^{N-1}_{s=M}\big( \hat{\chi}(A_{s}\Lambda(v)) -  \int_{\cV}\hat{\chi}(A_{s}\Lambda(v))dm_{\cV}(v) \big),
\eeqnas
for some $M = M(N)$, and
\beqnas
F^{\epsilon}_{N}(v) =  \frac{1}{\sqrt{N}}\sum^{N-1}_{s=M}\big( \hat{f}_{\epsilon}(A_{s}\Lambda(v)) -  \int_{\cV} \hat{f}_{\epsilon}(A_{s}\Lambda(v))dm_{\cV}(v)\big).
\eeqnas
for some  $\epsilon = \epsilon(N)$. Moreover, define
\beqnas
\hat{f}^{L}_{\epsilon}(y) = \hat{f}_{\epsilon}(y)\eta_{L}(\pi(y))\beta_{L}(\pi(y))
\eeqnas
for $y \in Y$, and 
\beqnas
F^{\epsilon, L}_{N}(v) = \frac{1}{\sqrt{N}}\sum^{N-1}_{s = M}\big( \hat{f}_{\epsilon}(A_{s}\Lambda(v)) -  \int_{\cV} \hat{f}_{\epsilon}(A_{s}\Lambda(v))dm_{\cV}(v) \big),
\eeqnas
for some $L = L(N)$. 
 By \eqref{est.4.4.4} and \eqref{est.eps.1}, we derive that
\beqna\label{est.snorm}
\cS(\hat{f}_{\epsilon}^{L}) = \cS_{l_{0}}(\hat{f}_{\epsilon}^{L}) =  O(\epsilon^{-l_{0}}L^{l_{0}+1}) .
\eeqna

We have
\beqna\label{est.sum.1}
\int_{\cV}|F_{N} - \tilde{F}_{N}|dm_{\cV} = O(\frac{M}{\sqrt{N}})
\eeqna
by \eqref{est.sg.1}, 
\beqna\label{est.sum.2}
\int_{\cV}|\tilde{F}_{N} - F^{\epsilon}_{N}|dm_{\cV} = O(\sqrt{N}(\epsilon + e^{-M}))
\eeqna
by \eqref{est.eps}. Notice that by \eqref{est.alp.1} and \eqref{est.ht}, for $p < d$ and $s$ large enough, 
\beqna\label{est.L.3}
\nonumber & &\int_{\cV}|\hat{f}_{\epsilon}(A_{s}\Lambda(v))- \hat{f}_{\epsilon}^{L}(A_{s}\Lambda(v))|dm_{\cV}(v) =  \int_{\cV}|\hat{f}_{\epsilon}(A_{s}\Lambda(v))||1 - \eta_{L}(\pi(A_{s}\Lambda(v)))\beta_{L}(\pi(A_{s}\Lambda(v)))|dm_{\cV}(v)\\ 
\nonumber  &\leq&  (\int_{\cV}|\hat{f}_{\epsilon}(A_{s}\Lambda(v))|^{2}dm_{\cV})^{\frac{1}{2}}\big(m_{\cV}(\alpha(\pi(A_{s}\Lambda(v))) \geq c^{-1}L) + m_{\cV}(\htt(\pi(A_{s}\Lambda(v))) \geq c^{-1}L)\big)^{\frac{1}{2}}\\
 &\ll&  (\int_{\cV}|\hat{f}_{\epsilon}(A_{s}\Lambda(v))|^{2}dm_{\cV})^{\frac{1}{2}}\big(L^{-p} + L^{-d}\big)^{\frac{1}{2}} \leq  (\int_{\cV}|\hat{f}_{\epsilon}(A_{s}\Lambda(v))|^{2}dm_{\cV})^{\frac{1}{2}}L^{-\frac{p}{2}},
\eeqna
where the boundedness of $(\int_{\cV}|\hat{f}_{\epsilon}(A_{s}v)|^{2}dm_{\cV})^{\frac{1}{2}}$ is ensured by \eqref{est.sg.2}. Thus, 
\beqna\label{est.sum.3}
\int_{\cV}|F^{\epsilon}_{N} - F^{\epsilon, L}_{N}|dm_{\cV} = O(\sqrt{N}L^{-\frac{p}{2}}).
\eeqna
Then by choosing proper parameters $\epsilon$, $M$ and $L$ such that as $N \rightarrow 0$, 
\beqna
\label{p}\frac{M}{\sqrt{N}} \rightarrow 0, \ \sqrt{N}(\epsilon + e^{-M})  \rightarrow 0, \ \sqrt{N}L^{-\frac{p}{2}} \rightarrow 0,
\eeqna
it remains to verify that $F^{\epsilon, L}_{N}$ satisfies the CLT. 

Now we adopt the notations and techniques developed in \cite{BG19} to estimate the cumulants and the variance. As in Section 3.2, \cite{BG19}, for $r \geq 3$ define $\{\alpha_{0}, \beta_{1}, \alpha_{1}, \dots, \beta_{r+1}\}$ as 
\beqna\label{pr.ab}
\alpha_{0} = 0, \beta_{1} =  \gamma, \ \alpha_{k} < \beta_{k+1}, \ \alpha_{k} = (3 +r)\beta_{k}, \ k=1, \dots r.
\eeqna
For $\bar{s} = (s_{0}, s_{1}, \dots, s_{r}) \in \bbR^{r+1}_{+}$ with $s_{0} = 0$, $I, J \subset \{0, \dots, r\}$, let
\beqnas
\rho^{I} = \max\{|s_{i} - s_{j}|, i, j \in I\},  \ \rho_{I, J}(\bar{s}) = \min\{|s_{i} - s_{j}|, i \in I, j \in J\}.
\eeqnas
For a partition $\cQ$ of $\{0, \dots, r\}$, let
\beqnas
\rho^{\cQ}(\bar{s}) = \max\{\rho^{I}, I \in \cQ\}, \ \rho_{\cQ}(\bar{s}) = \min\{ \rho_{I, J}(\bar{s}), I \neq J, I, J \in \cQ\},
\eeqnas
and for $\alpha < \beta$, 
\beqnas
\Delta_{\cQ}(\alpha, \beta) &=& \{\bar{s} \in \bbR^{r+1}_{+}, \rho^{\cQ}(\bar{s}) \leq \alpha, \ \rho_{\cQ}(\bar{s}) > \beta\}, \\
\Delta(\alpha) &=& \{\bar{s} \in \bbR^{r+1}_{+}, |s_{i} - s_{j}| \leq \alpha \ for \ all \ i, j\}.
\eeqnas
Define $\Delta(\beta)$ and $\Delta_{\cQ}(\alpha, \beta)$ for a partition $\cQ$ of $\{0, \dots, r\}$. Then $\{M, \dots, N-1\}^{r}$ can be decomposed into the following
\beqnas
\{M, \dots, N-1\}^{r} = \Omega(\beta_{r+1}; M, N) \cup \big( \bigcup^{r}_{j=0}\bigcup_{|\cQ| \geq 2} \Omega_{\cQ}(\alpha_{j}, \beta_{j+1}; M, N)\big),
\eeqnas 
where
\beqnas
\Omega(\beta_{r+1}; M, N) &=& \{M, \dots, N-1\}^{r} \cap \Delta(\beta_{r+1}),\\
\Omega_{\cQ}(\alpha_{j}, \beta_{j+1}; M, N) &=& \{M, \dots, N-1\}^{r} \cap \Delta_{\cQ}(\alpha_{j}, \beta_{j+1}).
\eeqnas
The estimate follows Section 5.1 in \cite{BG19}. Let $M > \beta_{r+1}$ such that $\Omega(\beta_{r+1}; M, N) = \emptyset$. 
Also notice that for $|s_{i} - s_{1}| \leq \alpha_{j}$, we have
\beqnas
\cS_{l_{0}}(\prod_{i \in I}\hat{f}_{\epsilon}^{L} \circ A_{s_{i} - s_{1}}) &\leq& e^{|I|\bar{\omega}\alpha_{j}}\cS_{l_{0}}(\hat{f}_{\epsilon}^{L})^{|I|},
\eeqnas
where $\bar{\omega} = \bar{\omega}(\omega_{1}, \dots, \omega_{m}, n, l_{0}) > 0$ is some constant. Then as $(5.21)$ in  \cite{BG19}, we derive that for $r \geq 3$, 
\beqna\label{est.cum.3}
\nonumber |\Cum^{(r)}(F^{\epsilon, L}_{N})| &\ll& N^{\frac{r}{2}}e^{-\omega_{0}\delta_{\kappa} \beta_{j+1} + r\bar{\omega} \alpha_{j}}\cS(\hat{f}^{L}_{\epsilon})^{r} + N^{1 - \frac{r}{2}}\alpha_{j}^{r -1}L^{(r - d + 1)^{+}}\|f_{\epsilon}\|^{r}_{L^{\infty}}\\
&\ll& N^{\frac{r}{2}}e^{-\omega_{0}\delta_{\kappa} \gamma}L^{r(l_{0}+1)}\epsilon^{-rl_{0}} + N^{1 - \frac{r}{2}}\gamma^{r -1}L^{(r - d + 1)^{+}},
\eeqna
where the last line is due to \eqref{effe.eqdis.3} and \eqref{est.snorm}, and  $\gamma$ in \eqref{pr.ab} is chosen to make sure that $\omega_{0}\delta_{\kappa} \beta_{j+1} - r\bar{\omega} \alpha_{j} > \omega_{0}\delta_{\kappa} \gamma$. Then we yield that   
\beqna\label{cvg.cum}
\lim_{N \rightarrow \infty}\Cum^{(r)}(F^{\epsilon, L}_{N}) = 0
\eeqna
for $r \geq 3$ if
\beqna
\label{p2} N^{\frac{r}{2}}e^{-\omega_{0}\delta_{\kappa} \gamma}L^{r(l_{0}+1)}\epsilon^{-rl_{0}} \rightarrow 0, \  N^{1 - \frac{r}{2}}\gamma^{r -1}L^{(r - d + 1)^{+}} \rightarrow 0.
\eeqna
Next we proceed to prove the convergence of the variance,
\beqnas
\int_{\cV} |F^{\epsilon, L}_{N}|^{2}dm_{\cV} = \Theta^{\epsilon, L}_{N}(0) + 2\sum^{N-M-1}_{s=1}\Theta^{\epsilon, L}_{N}(s), 
\eeqnas
where 
\beqnas
\psi_{s}(v) &:=& \hat{f}^{L}_{\epsilon}(A_{s}\Lambda(v)) -  \int_{\cV}\hat{f}^{L}_{\epsilon}(A_{s}\Lambda(v))dm_{\cV}(v), \\
\Theta^{\epsilon, L}_{N}(s) &:=& \frac{1}{N}\sum^{N-1-s}_{t = M}\int_{\cV}\psi_{s+t}(v)\psi_{t}(v)dm_{\cV}(v). 
\eeqnas
Then for some $K = K(N)$ to be decided later, we have
\beqna\label{est.v.k1}
\sum^{N-M-1}_{s=K}\Theta^{\epsilon, L}_{N}(s)  \ll e^{-\omega_{0}\delta_{\kappa} K}\cS(\hat{f}^{L}_{\epsilon})^{2} \ll  e^{-\omega_{0}\delta_{\kappa} K}L^{2(l_{0}+1)}\epsilon^{-2l_{0}}, 
\eeqna
and for $s < K$, 
\beqna\label{est.v.k2}
\Theta^{\epsilon, L}_{N}(s) = \frac{N-M-s}{N}\Theta^{\epsilon, L}_{\infty}(s) + O(N^{-1}e^{-\omega_{0}\delta_{\kappa} M + \bar{\omega} s}L^{2(l_{0}+1)}\epsilon^{-2l_{0}}), 
\eeqna
where 
\beqnas
\Theta^{\epsilon, L}_{\infty}(s) := \int_{Y} \hat{f}^{L}_{\epsilon}(A_{s}y)  \hat{f}^{L}_{\epsilon}(y)dm_{Y}(y) -  (\int_{Y}\hat{f}^{L}_{\epsilon}(y)dm_{Y}(y))^{2}.
\eeqnas
Moreover, define
\beqnas
\Theta^{\epsilon}_{\infty}(s) := \int_{Y} \hat{f}_{\epsilon}(A_{s}y)  \hat{f}_{\epsilon}(y)dm_{Y}(y) -  (\int_{Y}\hat{f}_{\epsilon}(y)dm_{Y}(y))^{2}, 
\eeqnas
then by \eqref{est.4.4.5} and \eqref{est.4.4.6}, for $1 < p < d$, 
\beqna\label{est.v.k3}
\Theta^{\epsilon, L}_{\infty}(s) = \Theta^{\epsilon}_{\infty}(s) + O(L^{-\frac{p-2}{2}}).
\eeqna
Combining \eqref{est.v.k1}, \eqref{est.v.k2} and \eqref{est.v.k3} together, we derive that
\beqna\label{est.v.k4}
\nonumber & &\int_{\cV} |F^{\epsilon, L}_{N}|^{2}dm_{\cV} 
 = \Theta^{\epsilon}_{\infty}(0) + 2\sum^{K-1}_{s=1}\Theta^{\epsilon}_{\infty}(s) \\
\nonumber & & \qquad + O(KL^{-\frac{p-2}{2}} + N^{-1}(M+K)K  + N^{-1}e^{-\omega_{0}\delta_{\kappa}  M + \bar{\omega} K}L^{2(l_{0}+1)}\epsilon^{-2l_{0}} + e^{-\omega_{0}\delta_{\kappa} K}L^{2(l_{0}+1)}\epsilon^{-2l_{0}}).\\
\eeqna
We claim the following estimate to complete the convergence of the variance.
\blem\label{l5.4}
\beqna\label{est.v.k5}
\lim_{\epsilon \rightarrow 0}\big( \Theta^{\epsilon}_{\infty}(0) + 2\sum^{K-1}_{s=1}\Theta^{\epsilon}_{\infty}(s) \big) = \sigma_{m, n}^{2}, 
\eeqna
where
\beqnas
\sigma_{m, n}^{2} = \int_{\bbR^{d}}\chi(x)dx = 2^{m}(\prod^{m}_{i=1}\theta_{i})n\vol_{n}.
\eeqnas
\elem
\bpf
Notice that by Rogers formula~\eqref{rog.1} and \eqref{rog.2}, 
\beqnas
\Theta^{\epsilon}_{\infty}(s) = \int_{\bbR^{d}}f_{\epsilon}(A_{s}x)f_{\epsilon}(x)dx, 
\eeqnas
such that
\beqnas
\Theta^{\epsilon}_{\infty}(0) + 2\sum^{K-1}_{s=1}\Theta^{\epsilon}_{\infty}(s) =  \sum^{K-1}_{s = -K+1}\int_{\bbR^{d}}f_{\epsilon}(A_{s}x)f_{\epsilon}(x)dx.
\eeqnas
Recall that
\beqnas
\Omega^{\epsilon}_{e} &=& \{(\bfx_{1}, \bfx_{2}) \in \bbR^{m+n}, 1- \epsilon \leq \|\bfx_{2}\| \leq e + \epsilon, \ |x_{1}^{i}| < \theta_{i}(\epsilon)\|\bfx_{2}\|^{-\omega_{i}}, \ i=1, \dots, m\},\\
A^{-s}\Omega^{\epsilon}_{e} &=& \{(\bfx_{1}, \bfx_{2}) \in \bbR^{m+n}, \ e^{s}(1 - \epsilon) \leq \|\bfx_{2}\| < e^{s}(e + \epsilon), \ |x_{1}^{i}| < \theta_{i}(\epsilon)\|\bfx_{2}\|^{-\omega_{i}}, \ i =1, \dots, m\},
\eeqnas
then we deduce that $\{f_{\epsilon}(A_{s}x)f_{\epsilon}(x) \neq 0\} $ only holds for $s = -1, 0, 1$. We compute that
\beqnas
& &|\Theta^{\epsilon}_{\infty}(0) + 2\sum^{K-1}_{s=1}\Theta^{\epsilon}_{\infty}(s) - \sigma_{m, n}^{2}|  = | \sum_{s = -1, 0, 1}\int_{\bbR^{m+n}}f_{\epsilon}(A_{s}x)f_{\epsilon}(x)dx - \int_{\bbR^{m+n}}\chi(x)dx| \\
&<& \int_{1 - \epsilon \leq \|\bfx_{2}\| \leq 1} \int_{|x_{1}^{i}| < \theta_{i}(\epsilon)\|\bfx_{2}\|^{-\omega_{i}}}d\bfx_{1}d\bfx_{2} +  \int_{e \leq \|\bfx_{2}\| \leq e +  \epsilon}\int_{|x_{1}^{i}| < \theta_{i}(\epsilon)\|\bfx_{2}\|^{-\omega_{i}}}d\bfx_{1}d\bfx_{2}\\
&\ll& \epsilon,
\eeqnas
which implies \eqref{est.v.k5}.
%
%
Moreover, we obtain the explicit formula for the variance,
\beqnas
\sigma^{2}_{m, n} &=& \int_{\bbR^{m+n}}\chi(x)dx = 2^{m}\prod^{m}_{i=1}\theta_{i} \int_{1 \leq \|\bfy\| \leq e}\|\bfy\|^{-n}d\bfy \\
&=& 2^{m}\prod^{m}_{i=1}\theta_{i}  \int_{S^{n-1}}\int^{e}_{1} r^{-1}dr  d\bfz = 2^{m}(\prod^{m}_{i=1}\theta_{i})n\vol_{n},
\eeqnas
where $S^{n-1}$ is the unit sphere with respect to the norm $\|\cdot \|$. 
\epf

Combining \eqref{est.v.k4} and \eqref{est.v.k5} together, we derive that
\beqna\label{cvg.var}
\lim_{N \rightarrow \infty}\int_{\cV}|F^{\epsilon, L}_{N}|^{2}dm_{\cV} = \sigma_{m, n}^{2}, 
\eeqna
provided that
\beqna\label{p6}
KL^{-\frac{p-2}{2}} + KN^{-1}(M+K) + KN^{-1}e^{-\omega_{0}\delta_{\kappa} M + \bar{\omega} K}L^{2(l_{0}+1)}\epsilon^{-2l_{0}} + e^{-\omega_{0}\delta_{\kappa} K}L^{2(l_{0}+1)}\epsilon^{-2l_{0}} \rightarrow 0
\eeqna
as $N \rightarrow \infty$. 
Together with the conditions~\eqref{p}, \eqref{p2},  we may take
\beqna\label{para}
M = (\log N)(\log \log N), \ L = N^{q}, \ K = \kappa_{0}\log N, \ \epsilon = N^{-1} ,\  \gamma = \gamma_{0}\log N
\eeqna
for some proper $\frac{1}{d} < q < \frac{(\frac{r}{2} - 1)}{(r - d + 1)^{+}}$, and $\gamma_{0}$, $\kappa_{0}$ large enough. Then with \eqref{cvg.cum} and \eqref{cvg.var}, we complete the proof.
\hfill \qedsymbol

\section{Proof of Theorem~\ref{ip.inDA}}\label{ip}
For $t \in [0, 1]$, define the piecewise linear function
\beqna
\nonumber X_{N, t}(v) =  \left\{
\nonumber \begin{aligned}
\label{def.sN.1} \sum^{[\log N]-1}_{s=0}\hat{\chi}(A_{s}\Lambda(v)) + \frac{\{t \log N\}}{\{\log N\}}R_{N}(v),  &\ t_{0} \leq t \leq 1, \\
\label{def.sN.2} \sum^{[t \log N]-1}_{s=0}\hat{\chi}(A_{s}\Lambda(v)) + \{t \log N\}\hat{\chi}(A_{[t \log N]}\Lambda(v)),  &\ \frac{1}{\log N} \leq t <  t_{0}, \\
\label{def.sN.3} t \log N\hat{\chi}(A_{0}\Lambda(v)),  &\ 0 \leq t < \frac{1}{\log N},
\end{aligned}
\right.
\eeqna
where $t_{0} = \inf\{t < 1, [t\log N] = [\log N]\}$, such that $X_{N, t}(v)$ interpolates $0$ and $\Delta_{N}$. For $0 \leq t \leq 1$, define
\beqnas
S_{N, t}(v) = X_{N, t}(v) -  C_{m, n}t\log N, \ D_{N, t}(v) = \frac{1}{\sigma_{m, n}\sqrt{\log N}}S_{N, t}(v).
\eeqnas
The proof of the IP proceeds by standard methods and consists of two parts: establishing convergence of the finite-dimensional distributions, and proving tightness.

\subsection{Convergence of finite-dimensional distributions}

\bprop
For $0 \leq t_{1} < t_{2} < \dots < t_{k} \leq 1$, we have
\beqnas
\big(D_{N, t_{1}}, D_{N, t_{2}}, \dots, D_{N, t_{k}} \big)  \Longrightarrow \big( W_{t_{1}}, W_{t_{2}}, \dots, W_{t_{k}} \big),
\eeqnas
as $N \rightarrow \infty$.
\eprop

\bpf
Without loss of generality, we may assume $t_{1} \geq \frac{1}{\log N}$. It is equivalent to prove that 
\beqnas
\big(D_{N, t_{2}} - D_{N, t_{1}},  \dots, D_{N, t_{k}} - D_{N, t_{k-1}} \big)  \Longrightarrow \big( W_{t_{2}} - W_{t_{1}}, \dots, W_{t_{k}} - W_{t_{k-1}} \big).
\eeqnas
as $N \rightarrow \infty$. First notice that
\beqnas
& &D_{N, t_{i+1}}(v)  - D_{N, t_{i}}(v) \\
&=& \frac{1}{\sigma_{m, n}\sqrt{\log N}}\big( \sum^{[t_{i+1}\log N] -1}_{s=[t_{i}\log N]}\hat{\chi}(A_{s}\Lambda(v)) + \cE_{N, t_{i}, t_{i+1}}  - C_{m, n}(t_{i+1} - t_{i})\log N \big),
\eeqnas
where 
\beqnas
\cE_{N, t_{i}, t_{i+1}}(v) &=&  \left\{
\begin{aligned}
 \frac{\{t_{i+1} \log N\}}{\{\log N\}}R_{N}(v) -  \frac{\{t_{i} \log N\}}{\{\log N\}}R_{N}(v),  & \  t_{0} \leq t_{i} < t_{i+1}, \\
  \frac{\{t_{i+1} \log N\}}{\{\log N\}}R_{N}(v) -   \{t_{i}\log N\}\hat{\chi}(A_{[t_{i}\log N]}\Lambda(v)),  & \ t_{i} < t_{0} \leq  t_{i+1}, \\
\{t_{i+1}\log N\}\hat{\chi}(A_{[t_{i+1}\log N]}\Lambda(v))  -  \{t_{i}\log N\}\hat{\chi}(A_{[t_{i}\log N]}\Lambda(v)), & \  t_{i} < t_{i+1} <  t_{0}.
\end{aligned}
\right.
\eeqnas
Notice that in each case, we have
\beqnas
|\cE_{N, t_{i}, t_{i+1}}(v)| \leq  \hat{\chi}(A_{[t_{i}\log N]}\Lambda(v)) + \hat{\chi}(A_{[t_{i+1}\log N]}\Lambda(v)),
\eeqnas
such that by \eqref{est.exp}, we have
\beqna\label{est.err}
 \frac{1}{\sigma_{m, n}\sqrt{\log N}}\int_{\cV} |\cE_{N, t_{i}, t_{i+1}}(v)|dm_{\cV}(v) = O( \frac{1}{\sqrt{\log N}}).
\eeqna
Now define
\beqnas
F^{\epsilon, L}_{N, i}(v) = \frac{1}{\sigma_{m, n}\sqrt{\log N}}\sum^{[t_{i+1}\log N] - 1}_{s=[t_{i}\log N]}\big( \hat{f}^{L}_{\epsilon}(A_{s}\Lambda(v)) -  \int_{\cV}\hat{f}^{L}_{\epsilon}(A_{s}\Lambda(v))dm_{\cV}(v) \big),
\eeqnas
where $\hat{f}^{L}_{\epsilon}(y) = \hat{f}_{\epsilon}(y)\eta_{L}(\pi(y))\beta_{L}(\pi(y))$. Then by \eqref{est.err}, \eqref{est.exp}, \eqref{est.eps} and \eqref{est.L.3} for $p < d$, we derive that
\beqna\label{ip.est.1}
\nonumber & &\int_{\cV} |\big( D_{N, t_{i+1}}(v)  - D_{N, t_{i}}(v) \big) - F^{\epsilon, L}_{N, i}(v)|dm_{\cV}(v)  \\
& =&  O\big( \frac{1}{\sqrt{\log N}} + (t_{i+1} - t_{i})(\epsilon + e^{-[t_{i}\log N]}+L^{-\frac{p}{2}}) \big).
\eeqna
By choosing $\epsilon = \epsilon(N)$ and $L= L(N)$ suitably so that $\epsilon +  e^{-[t_{i} \log N]} + L^{-\frac{p}{2}} \rightarrow 0$ as $N \rightarrow +\infty$, it remains to prove
\beqna\label{ip.est.2}
\big(F^{\epsilon, L}_{N, 1}(v),  \dots,F^{\epsilon, L}_{N, k-1}(v)\big)  \Longrightarrow \big( W_{t_{2}} - W_{t_{1}}, \dots, W_{t_{k}} - W_{t_{k-1}} \big).
\eeqna
To this end, we employ the cumulant method and claim the following: estimates \eqref{est.cum.11} and \eqref{est.cum.22} imply that for each $1 < i < k- 1$, $F^{\epsilon, L}_{N, i}$ converges in distribution to $\cN(0, t_{i+1} - t_{i})$, while \eqref{est.cum.33} shows that the components of the limiting distribution are mutually independent.
\blem
For any $1 \leq i \leq k- 1$, we have
\beqna\label{est.cum.11}
\lim_{N \rightarrow \infty}\int_{\cV} |F^{\epsilon, L}_{N, i}(v)|^{2}dm_{\cV}(v) = t_{i+1} - t_{i},
\eeqna
and
\beqna\label{est.cum.22}
\lim_{N \rightarrow \infty}\Cum^{(r)}(F^{\epsilon, L}_{N, i}) = 0
\eeqna
for $r \geq 3$. Moreover, for $2 \leq r \leq k-1$, $(i_{1}, \dots, i_{r}) \subset \{1, \dots, k-1\}$, we have
\beqna\label{est.cum.33}
\lim_{N \rightarrow \infty}\Cum^{(r)}(F^{\epsilon, L}_{N, i_{1}}, F^{\epsilon, L}_{N, i_{2}}, \dots, F^{\epsilon, L}_{N, i_{r}}) = 0.
\eeqna
\elem

\bpf

The proof follows the same approach as Proposition~\ref{clt.chi}. We start with the convergence of the variance. Let $s_{1} = [t_{i}\log N]$, $s_{2} = [t_{i+1}\log N]$. 
\beqnas
\int_{\cV} |F^{\epsilon, L}_{N, i}|^{2}dm_{\cV} = \Theta^{\epsilon, L}_{N}(0) + 2\sum^{s_{2}-s_{1} - 1}_{s=1}\Theta^{\epsilon, L}_{N}(s), 
\eeqnas
where 
\beqnas
\psi_{s}(v) &:=& \hat{f}^{L}_{\epsilon}(A_{s}\Lambda(v)) -  \int_{\cV}\hat{f}^{L}_{\epsilon}(A_{s}\Lambda(v))dm_{\cV}(v), \\
\Theta^{\epsilon, L}_{N}(s) &:=& \frac{1}{\sigma^{2}_{m, n}\log N}\sum^{s_{2}-s-1}_{t = s_{1}}\int_{\cV}\psi_{s+t}(v)\psi_{t}(v)dm_{\cV}(v). 
\eeqnas
Then for some $\tilde{K}= (t_{i+1} - t_{i})K(N) <  s_{2} - s_{1} -1$ to be decided later, we have
\beqna\label{est.var.k1}
\sum^{s_{2}-s_{1}-1}_{s=\tilde{K}}\Theta^{\epsilon, L}_{N}(s)  \ll \frac{s_{2} - s_{1}}{\log N} e^{-\omega_{0}\delta_{\kappa} \tilde{K}}\cS(\hat{f}^{L}_{\epsilon})^{2} \ll  \frac{s_{2} - s_{1}}{\log N}e^{-\omega_{0}\delta_{\kappa} \tilde{K}}L^{2(l_{0}+1)}\epsilon^{-2l_{0}}, 
\eeqna
and for $s < \tilde{K}$, 
\beqna\label{est.var.k2}
\Theta^{\epsilon, L}_{N}(s) = \frac{s_{2}-s_{1}-s}{(t_{i+1} - t_{i})\log N}\Theta^{\epsilon, L}_{\infty}(s) + O(N^{-1}e^{-\omega_{0}\delta_{\kappa}s_{1} + \bar{\omega}s}L^{2(l_{0}+1)}\epsilon^{-2l_{0}}), 
\eeqna
where 
\beqnas
\Theta^{\epsilon, L}_{\infty}(s) := \frac{t_{i+1} - t_{i}}{\sigma^{2}_{m, n}}\big(\int_{Y} \hat{f}^{L}_{\epsilon}(A_{s}y)  \hat{f}^{L}_{\epsilon}(y)dm_{Y}(y) -  (\int_{Y}\hat{f}^{L}_{\epsilon}(y)dm_{Y}(y))^{2}\big),
\eeqnas
\beqnas
\Theta^{\epsilon}_{\infty}(s) := \frac{t_{i+1} -t_{i}}{\sigma^{2}_{m, n}}\big(\int_{Y} \hat{f}_{\epsilon}(A_{s}y)  \hat{f}_{\epsilon}(y)dm_{Y}(y) -  (\int_{Y}\hat{f}_{\epsilon}(y)dm_{Y}(y))^{2}\big).
\eeqnas
Moreover, by \eqref{est.4.4.5} and \eqref{est.4.4.6},
\beqna\label{est.var.k3}
\Theta^{\epsilon, L}_{\infty}(s) = \Theta^{\epsilon}_{\infty}(s) + O(L^{-\frac{p-2}{2}})
\eeqna
holds for $1 < p < d$. Combining \eqref{est.var.k1}, \eqref{est.var.k2} and \eqref{est.var.k3} together, we derive that
\beqna\label{est.var.k4}
\nonumber & &\int_{\cV} |F^{\epsilon, L}_{N}|^{2}dm_{\cV} 
 \nonumber = \Theta^{\epsilon}_{\infty}(0) + \frac{2(s_{2}-s_{1})}{(t_{i+1} - t_{i})\log N}\sum^{\tilde{K}-1}_{s=1}\Theta^{\epsilon}_{\infty}(s) + O(\tilde{K}L^{-\frac{p-2}{2}}) \\
\nonumber& & \qquad + O(\tilde{K}^{2}(\log N)^{-1} + (\log N)^{-1}e^{-\omega_{0}\delta_{\kappa}s_{1} +  \bar{\omega} \tilde{K}}L^{2(l_{0}+1)}\epsilon^{-2l_{0}} + (s_{2} - s_{1})(\log N)^{-1}e^{-\omega_{0}\delta_{\kappa} \tilde{K}}L^{2(l_{0}+1)}\epsilon^{-2l_{0}}).\\
\eeqna
Notice that 
$$
 s_{2}- s_{1}-1 \leq (t_{i+1} - t_{i})\log N \leq s_{2} -s_{1}+1,  
$$
implying
\beqna\label{limit.t}
\lim_{N \rightarrow \infty}\frac{s_{2}-s_{1}}{(t_{i+1} - t_{i})\log N} = 1.
\eeqna
Simliar arguments as in the proof of \eqref{est.v.k5} lead to 
\beqnas
\Theta^{\epsilon}_{\infty}(0) + \frac{2(s_{2}-s_{1})}{(t_{i+1} - t_{i})\log N}\sum^{\tilde{K}-1}_{s=1}\Theta^{\epsilon}_{\infty}(s) = (t_{i+1} - t_{i})(1 + O(\epsilon)). 
\eeqnas
Then we prove \eqref{est.cum.11} once
\beqna\label{cd.1}
\epsilon + \tilde{K}L^{-\frac{p-2}{2}} + \tilde{K}^{2}(\log N)^{-1} + ((\log N)^{-1}e^{-\omega_{0}\delta_{\kappa}s_{1} + \bar{\omega} \tilde{K}} + e^{-\omega_{0}\delta_{\kappa} \tilde{K}})L^{2(l_{0}+1)}\epsilon^{-2l_{0}} \rightarrow 0
\eeqna
holds as $N \rightarrow \infty$.

To estimate the cumulants, we follow the procedure in the proof of Proposition~\ref{clt.chi}. We have the following decomposition
\beqna\label{decomp.2}
\{s_{1}, \dots, s_{2}-1\}^{r} = \Omega(\beta_{r+1}; s_{1}, s_{2}-1) \cup \big( \bigcup^{r}_{j=0} \bigcup_{|\cQ| \geq 2} \Omega_{\cQ}(\alpha_{j}, \beta_{j+1}; s_{1}, s_{2}-1)\big),
\eeqna
where we adopt the same notations $\{\alpha_{0}, \beta_{1}, \dots, \beta_{r+1}\}$ as in \eqref{pr.ab} and take $s_{1} > \beta_{r+1}$ so that $\Omega(\beta_{r+1}; s_{1}, s_{2}-1) = \emptyset$. We set $\gamma = (t_{i+1} -t_{i})\gamma_{0}\log N$ for some $\gamma_{0}$ to be decided later. Notice that for $\cQ = \{\{0\}, \{1, \dots, r\}\}$, 
\beqnas
\sharp\Omega_{\cQ}(\alpha_{j}, \beta_{j+1}; s_{1}, s_{2}-1) \leq (s_{2} -s_{1})\alpha^{r-1}_{j},
\eeqnas
and for $|\cQ| \geq 2$ and $\cQ \neq  \{\{0\}, \{1, \dots, r\}\}$, 
\beqnas
\sharp\Omega_{\cQ}(\alpha_{j}, \beta_{j+1}; s_{1}, s_{2}-1) \leq (s_{2} - s_{1})^{r}. 
\eeqnas
Then by the same estimate as \eqref{est.cum.3}, we obtain
\beqna\label{cum.r}
\nonumber |\Cum^{(r)}(F^{\epsilon, L}_{N, i}) |  &\ll& (s_{2} - s_{1})^{r}(\log N)^{-\frac{r}{2}}e^{-\omega_{0}\delta_{\kappa} \gamma}L^{r(l_{0}+1)}\epsilon^{-rl_{0}} + (s_{2} - s_{1})(\log N)^{ - \frac{r}{2}}\gamma^{r -1}L^{(r - d + 1)^{+}}.\\
\eeqna
Hence we prove \eqref{est.cum.22} if 
\beqna\label{cd.2}
(\log N)^{\frac{r}{2}}e^{-\omega_{0}\delta_{\kappa} \gamma}L^{r(l_{0}+1)}\epsilon^{-rl_{0}} + (\log N)^{1 - \frac{r}{2}}\gamma^{r -1}L^{(r - d + 1)^{+}} \rightarrow 0
\eeqna
as $N \rightarrow \infty$.

Now we proceed to prove \eqref{est.cum.33}. Notice that for any $\{i_{1}, \dots, i_{r}\} \subset \{1, \dots, k-1\}$, $2 \leq r \leq k-1$, we have 
\beqnas
\Cum^{(r)}(F^{\epsilon, L}_{N, i_{1}}, F^{\epsilon, L}_{N, i_{2}}, \dots, F^{\epsilon, L}_{N, i_{r}}) = \frac{1}{\sigma_{m, n}^{r}(\log N)^{r \over 2}}\sum^{[t_{i_{1} +1}\log N]-1}_{z_{1} = [t_{i_{1}}\log N]} \dots \sum^{[t_{i_{r} +1}\log N] -1}_{z_{r} = [t_{i_{r}}\log N]}\Cum^{(r)}(\psi_{z_{1}}, \psi_{z_{2}}, \dots,  \psi_{z_{r}}),
\eeqnas
where $\psi_{z}(v) = \hat{f}^{L}_{\epsilon}(A_{z}\Lambda(v)) -  \int_{\cV}\hat{f}^{L}_{\epsilon}(A_{z}\Lambda(v))dm_{\cV}(v) $. When $r \geq 3$, the estimate follows the same route as above. Let $s_{i_{1}} = [t_{i_{1}}\log N]$ and $s_{i_{r}+1} =  [t_{i_{r}+1}\log N]$, and we decompose $\{s_{i_{1}}, \dots, s_{i_{r}+1}-1\}^{r}$ as in \eqref{decomp.2}. Then once \eqref{cd.2} is satisfied, \eqref{est.cum.33} follows for $r \geq 3$.
When $r = 2$, the cumulant is just the variance. For $(i_{1}, i_{2}) \subset \{1, \dots, k-1\}$ fixed, we have
\beqnas
\Cum^{(2)}(F^{\epsilon, L}_{N, i_{1}}, F^{\epsilon, L}_{N, i_{2}}) 
&=&  \frac{1}{\sigma_{m, n}^{2}\log N}\sum^{[t_{i_{1} +1}\log N]-1}_{z_{1} = [t_{i_{1}}\log N]}\sum^{[t_{i_{2} +1}\log N]-1}_{z_{r} = [t_{i_{2}}\log N]} \int_{\cV}\psi_{z_{1}}(v)\psi_{z_{2}}(v)dm_{\cV}(v),
\eeqnas 
 where
\beqnas
& &\int_{\cV}\psi_{z_{1}}(v)\psi_{z_{2}}(v)dm_{\cV}(v) \\
&=& \int_{\cV} \hat{f}^{L}_{\epsilon}(A_{z_{1}}\Lambda(v))  \hat{f}^{L}_{\epsilon}(A_{z_{2}}\Lambda(v))dm_{\cV}(v) - \int_{\cV}\hat{f}^{L}_{\epsilon}(A_{z_{1}}\Lambda(v))dm_{\cV}(v)\int_{\cV}\hat{f}^{L}_{\epsilon}(A_{z_{2}}\Lambda(v))dm_{\cV}(v).
\eeqnas

We may assume that $i_{1} < i_{2}$. Notice that there exists $c_{k} > 0$, which depends on $(t_{1}, t_{2}, \dots, t_{k})$, such that $D([t_{1}\log N], [t_{2}\log N], \dots, [t_{k}\log N]) > c_{k}\log N$ for $N$ large enough. If $i_{1}+1 < i_{2}$, $z_{2} - z_{1} > c_{k}N$, and
\beqnas
\int_{\cV} \hat{f}^{L}_{\epsilon}(A_{z_{1}}\Lambda(v)) \hat{f}^{L}_{\epsilon}(A_{z_{2}}\Lambda(v)) dm_{\cV} = \big(\int_{Y}\hat{f}^{L}_{\epsilon}dm_{Y} \big)^{2} + O(\cS(\hat{f}^{L}_{\epsilon})^{2}e^{-\omega_{0}\delta_{\kappa}D(z_{1}, z_{2})}).
\eeqnas
Meanwhile, we have for $i =1, 2$
\beqnas
 \int_{\cV}\hat{f}^{L}_{\epsilon}(A_{z_{i}}\Lambda(v))dm_{\cV}(v) = \int_{Y}\hat{f}^{L}_{\epsilon}dm_{Y} + O(\cS(\hat{f}^{L}_{\epsilon})e^{-\omega_{0}\delta_{\kappa}z_{i}}).
\eeqnas
Thus, for $i_{2} - i_{1} > 1$, 
\beqna\label{est.cum.2.1}
\int_{\cV}\psi_{z_{1}}(v)\psi_{z_{2}}(v)dm_{\cV}(v) = O(e^{-\omega_{0}\delta_{\kappa}D(z_{1}, z_{2})}\cS(\hat{f}^{L}_{\epsilon})^{2}) = O(e^{-\omega_{0}\delta_{\kappa}D(z_{1}, z_{2})}L^{2(l_{0}+1)}\epsilon^{-2l_{0}}),
\eeqna
and 
\beqna\label{est.cum.2.1.1}
\Cum^{(2)}(F^{\epsilon, L}_{N, i_{1}}, F^{\epsilon, L}_{N, i_{2}}) = O(e^{-\omega_{0}\delta_{\kappa}c_{k}N}L^{2(l_{0}+1)}\epsilon^{-2l_{0}}).
\eeqna
If $i_{1} +1 = i_{2}$, we have $1 \leq z_{2} - z_{1} \leq (t_{i_{2} + 1} - t_{i_{1}})N$. If $z_{2} - z_{1} > \tilde{K}_{1}$ for some $\tilde{K}_{1} = \min\{t_{i_{1}}, (t_{i_{2} + 1} - t_{i_{1}})\}K_{1}(N)$, $K_{1}(N) < N$ to be decided later, we yield
\beqna\label{est.cum.2.2}
\sum^{[t_{i_{2} +1}N]-1}_{z_{2} > z_{1} + \tilde{K}_{1}}\int_{\cV}\psi_{z_{1}}(v)\psi_{z_{2}}(v)dm_{\cV}(v) = O(\cS(\hat{f}^{L}_{\epsilon})^{2}e^{-\omega_{0}\delta_{\kappa}\tilde{K}_{1}}).
\eeqna
If $z_{2} - z_{1} \leq \tilde{K}_{1}$, 
\beqnas
\int_{\cV} \hat{f}^{L}_{\epsilon}(A_{z_{1}}\Lambda(v)) \hat{f}^{L}_{\epsilon}(A_{z_{2}}\Lambda(v)) dm_{\cV} =  \int_{Y}\hat{f}^{L}_{\epsilon}(\hat{f}^{L}_{\epsilon} \cdot A_{z_{2} - z_{1}})dm_{Y}  + O(\cS(\hat{f}^{L}_{\epsilon})^{2}e^{-\omega_{0}\delta_{\kappa}z_{1} + \bar{\omega} \tilde{K}_{1}}),
\eeqnas
which leads to
\beqna\label{est.cum.2.3}
\int_{\cV}\psi_{z_{1}}(v)\psi_{z_{2}}(v)dm_{\cV}(v) = \int_{Y}\hat{f}^{L}_{\epsilon}(\hat{f}^{L}_{\epsilon} \cdot A_{z_{2} - z_{1}})dm_{Y} -  \big(\int_{Y}\hat{f}^{L}_{\epsilon}dm_{Y} \big)^{2} + O(\cS(\hat{f}^{L}_{\epsilon})^{2}e^{-\omega_{0}\delta_{\kappa}z_{1} +  \bar{\omega} \tilde{K}_{1}}).
\eeqna
Combining \eqref{est.cum.2.2} and \eqref{est.cum.2.3} together, we derive that when $i_{1} +1 = i_{2}$,
\beqna\label{est.cum.2.4}
\nonumber & &\Cum^{(2)}(F^{\epsilon, L}_{N, i_{1}}, F^{\epsilon, L}_{N, i_{2}}) \\
\nonumber &=& \frac{1}{\sigma_{m, n}^{2}\log N}\sum^{[t_{i_{1} +1}\log N]-1}_{z_{1} = [t_{i_{1}}\log N]} \big( \sum^{\tilde{K}_{1}}_{s=1}\big(\int_{Y}\hat{f}^{L}_{\epsilon}(\hat{f}^{L}_{\epsilon} \cdot A_{s})dm_{Y} -  \big(\int_{Y}\hat{f}^{L}_{\epsilon}dm_{Y} \big)^{2}\big) + O(\cS(\hat{f}^{L}_{\epsilon})^{2}( e^{-\omega_{0}\delta_{\kappa}\tilde{K}_{1}} + e^{-\omega_{0}\delta_{\kappa}z_{1} + \bar{\omega} \tilde{K}_{1}})\big)\\
&\ll& \sum^{\tilde{K}_{1}}_{s=1}\big(\int_{Y}\hat{f}^{L}_{\epsilon}(\hat{f}^{L}_{\epsilon} \cdot A_{s})dm_{Y} -  \big(\int_{Y}\hat{f}^{L}_{\epsilon}dm_{Y} \big)^{2}\big) + O(L^{2(l_{0}+1)}\epsilon^{-2l_{0}}( e^{-\omega_{0}\delta_{\kappa}\nonumber \tilde{K}_{1}} + (\log N)^{-1}e^{-\omega_{0}\delta_{\kappa}[t_{i_{1}}\log N] + \bar{\omega} \tilde{K}_{1}})\big).\\
\eeqna
Note that
\beqna\label{est.cum.2.5}
\nonumber \sum^{\tilde{K}_{1}}_{s=1}\big(\int_{Y}\hat{f}^{L}_{\epsilon}(\hat{f}^{L}_{\epsilon} \cdot A_{s})dm_{Y} -  \big(\int_{Y}\hat{f}^{L}_{\epsilon}dm_{Y} \big)^{2}\big) = \sum^{\tilde{K}_{1}}_{s=1}\big(\int_{Y}\hat{f}_{\epsilon}(\hat{f}_{\epsilon} \cdot A_{s})dm_{Y} -  \big(\int_{Y}\hat{f}_{\epsilon}dm_{Y} \big)^{2}\big) + O(\tilde{K}_{1}L^{-\frac{p-2}{2}})\\
\eeqna
 for some $1 < p < d$. From the proof of Lemma~\ref{l5.4}, we deduce that
 \beqna\label{est.cum.2.6}
 \sum^{\tilde{K}_{1}}_{s=1}\big(\int_{Y}\hat{f}_{\epsilon}(\hat{f}_{\epsilon} \cdot A_{s})dm_{Y} -  \big(\int_{Y}\hat{f}_{\epsilon}dm_{Y} \big)^{2}\big) =  \sum^{\tilde{K}_{1}}_{s=1}\int_{\bbR^{d}}f_{\epsilon}(x)f_{\epsilon}(A_{s}x)dx =  \int_{\bbR^{d}}f_{\epsilon}(x)f_{\epsilon}(A_{1}x)dx, 
 \eeqna
and 
\beqna\label{est.cum.2.7}
|\int_{\bbR^{d}}f_{\epsilon}(x)f_{\epsilon}(A_{1}x)dx| \ll \epsilon. 
\eeqna
Therefore, by \eqref{est.cum.2.1.1}  and  \eqref{est.cum.2.4} with \eqref{est.cum.2.5}, \eqref{est.cum.2.6}, \eqref{est.cum.2.7}, we obtain \eqref{est.cum.33} holds for $r =2$ if
\beqna\label{cd.3}
\epsilon + \tilde{K}_{1}L^{-\frac{p-2}{2}} + L^{2(l_{0}+1)}\epsilon^{-2l_{0}}(e^{-\omega_{0}\delta_{\kappa}c_{k}\log N} + e^{-\omega_{0}\delta_{\kappa}\tilde{K}_{1}} + (\log N)^{-1}e^{-\omega_{0}\delta_{\kappa}[t_{i_{1}}\log N] + \bar{\omega} \tilde{K}_{1}}) \rightarrow 0
\eeqna
as $N \rightarrow \infty$. 
Now we choose
\beqna\label{cd.4}
\epsilon = (\log N)^{-1}, \ L = (\log N)^{q}, \ K = \kappa_{0}\log \log N, \ K_{1} = \kappa_{1}\log \log N, \ \gamma = (t_{i+1} - t_{i})\gamma_{0}\log \log N,
\eeqna
for some proper constant $\frac{1}{d} < q < \frac{(\frac{r}{2} - 1)}{(r - d + 1)^{+}}$, $r \geq 3$ and $\kappa_{0}, \kappa_{1}, \gamma_{0}$ large enough to ensure that \eqref{cd.1}, \eqref{cd.2} and \eqref{cd.3} hold. This finishes the proof of the lemma.
\epf
Notice that \eqref{cd.4} also guarantees \eqref{ip.est.1} vanishes as $N \rightarrow \infty$.   Thus we complete the proof of the convergence of finite dimensional distribution. 
\epf

\subsection{Tightness}
%
%
%
%


We quote the criterion (Theorem 8.3 in \cite{billing}) for the tightness of a family of probability measures. 
\bthm\label{c.tight}
A family of probability measures $\{P_{n}\}$ on $\cC([0, 1])$ is tight if the two conditions are satisfied:
\bitem
{\item For each positive $\eta$, there exists an $a$ such that for each $\bbP_{n}$
\beqnas
\bbP_{n}(|x(0)| > a) \leq \eta.
\eeqnas} 
{\item For each positive $\xi$ and $\eta$, there exists a $\delta$, $0 < \delta < 1$, and an integer $n_{0}$ such that
\beqna\label{cond.ip.1}   
\bbP_{n}(\sup_{t \leq s \leq t+\delta}|x(s) - x(t)| \geq \xi) \leq \eta\delta,
\eeqna
for $n \geq n_{0}$ and all $t$. 
}
\eitem
\ethm

%

Define
\beqnas
S_{n}(v) = \sum^{n-1}_{s=0}\hat{\chi}(A_{s}\Lambda(v))  - C_{m, n}n.
\eeqnas
where $n = [\log N]$. Since $D_{N, 0} = 0$, it suffices to verify \eqref{cond.ip.1} for $\{D_{N, t}\}_{N \in \bbN}$, which reduces to the following. 
\bprop\label{tight.s}
Given $\xi> 0$ and $\eta > 0$, there exists $0 < \delta < 1$, 
\beqna\label{cond.ip.2}  
\bbP\big(\max_{i \leq \delta \log N}\frac{1}{\sigma_{m, n}\sqrt{\log N}}|S_{k+i} - S_{k}| \geq \xi \big) \leq \eta \delta
\eeqna
for all $k$ and $N > N_{0}$.
\eprop

 To this end, assume \eqref{cond.ip.2} holds for given $\xi, \eta > 0$. If $s \geq t \geq t_{0}$, then $(s - t)\log N = \{s\log N\} - \{t\log N\} \leq \{\log N\} < 1$, which leads to
 \beqnas
 D_{N, s} - D_{N, t} = \frac{1}{\sigma_{m, n}\sqrt{\log N}}\big(\frac{\{s\log N\} - \{t\log N\}}{\{\log N\}}R_{N} - C_{m, n}(s - t)\log N \big), 
 \eeqnas
 and
 \beqnas
 \bbE |D_{N, s} - D_{N, t}| \leq \frac{1}{\sigma_{m, n}\sqrt{\log N}}(\bbE R_{N} + C_{m, n}) = O(\frac{1}{\sqrt{\log N}}), 
 \eeqnas 
such that there exists $0 < \delta <1$ and $N_{1}$, satisfying
 \beqnas
 \bbP(\sup_{t \leq s \leq t+\delta} |D_{N, s} - D_{N, t}|  \geq \xi) \leq \frac{1}{\xi}\bbE \sup_{t \leq s \leq t+\delta} |D_{N, s} - D_{N, t}| \leq \delta \eta, 
 \eeqnas
for $N  > N_{1}$.  Thus by Theorem~\ref{c.tight}, $\{D_{N, t}\}_{N \in \bbN}$ is tight. 
 
 For $t < t_{0}$, let $k = [t \log N]$. For $\log N > \frac{4}{\delta}$, one has $\frac{\delta}{2}\log N + 2  < \delta \log N$. If $s  \geq t_{0}$, we have
  \beqnas
 |D_{N, s} - D_{N, t}| &\leq&  |\frac{1}{\sigma_{m, n}\sqrt{\log N}}S_{k} - D_{N, s}| + |\frac{1}{\sigma_{m, n}\sqrt{\log N}}S_{k} - D_{N, t}|\\
 &\leq& \frac{2}{\sigma_{m, n}\sqrt{\log N}}\max_{i \leq \delta \log N}|S_{k+i} - S_{k}| + \frac{1}{\sigma_{m, n}\sqrt{\log N}}|\frac{\{s\log N\}}{\{\log N\}}R_{N} - C_{m, n}\{s\log N\}|.
 \eeqnas 
 If $t < s < t_{0}$, we have
  \beqnas
 |D_{N, s} - D_{N, t}| &\leq& |\frac{1}{\sigma_{m, n}\sqrt{\log N}}S_{k} - D_{N, s}| + |\frac{1}{\sigma_{m, n}\sqrt{\log N}}S_{k} - D_{N, t}|\\
 &\leq& \frac{2}{\sigma_{m, n}\sqrt{\log N}}\max_{i \leq \delta \log N}|S_{k+i} - S_{k}|.
 \eeqnas
The above estimates lead to
\beqnas
& &\bbP(\sup_{t \leq s \leq t+ \frac{\delta}{2}}|D_{N, s} - D_{N, t}| \geq 2\xi) \\
&\leq& \bbP(\frac{2}{\sigma_{m, n}\sqrt{N}}\max_{i \leq \delta N}|S_{k+i} - S_{k}| \geq \xi) + \bbP(\frac{1}{\sigma_{m, n}\sqrt{\log N}}|\frac{\{s\log N\}}{\{\log N\}}R_{N} - C_{m, n}\{s\log N\}| \geq \xi)\\
&\leq& 2\eta \delta,  
\eeqnas
for $N > \max\{N_{0}, N_{1}, \frac{4}{\delta}\}$. This implies the tightness of $\{D_{N, t}\}_{N \in \bbN}$ by Theorem~\ref{c.tight}.  

\medskip

\noindent
\emph{Proof of Theorem~\ref{tight.s}.}
For $0 \leq t \leq 1$, we write
\beqnas
S_{[t\log N]} &:=& S^{\epsilon, L}_{[t\log N]} + E^{1}_{[t\log N]} + E^{2}_{[t\log N]}, 
\eeqnas
where
\beqnas
S^{\epsilon, L}_{[t\log N]} &=&  \sum^{[t\log N]-1}_{s=0}\big( \hat{f}^{L}_{\epsilon}(A_{s}\Lambda(v)) - \bbE \hat{f}^{L}_{\epsilon}(A_{s}\Lambda(v))  \big),\\
E^{1}_{[t\log N]} &=& \sum^{[t\log N]-1}_{s=0}\big(\hat{\chi}(A_{s}\Lambda(v)) - \hat{f}^{L}_{\epsilon}(A_{s}\Lambda(v)) \big),\\
E^{2}_{[t\log N]} &=& \sum^{[t\log N]-1}_{s=0} \bbE \hat{f}^{L}_{\epsilon}(A_{s}\Lambda(v)) -  C_{m, n}[t\log N],
\eeqnas
and the parameters $\epsilon$, $L$ are chosen according to \eqref{cd.4}. By \eqref{est.eps} and \eqref{est.L.3}, we derive that 
\beqnas
\bbE \big(\max_{i \leq \delta \log N}|E^{1}_{k+i} - E^{1}_{k}| \big) &=& \bbE \big(\max_{i \leq \delta \log N}|\sum^{i}_{j=1}\hat{\chi}(A_{k+j}\Lambda(v)) - \hat{f}^{L}_{\epsilon}(A_{k+j}\Lambda(v))| \big) \\
&\leq&  \sum^{[\delta \log N]}_{j = 1}\bbE |\hat{\chi}(A_{k+j}\Lambda(v)) - \hat{f}^{L}_{\epsilon}(A_{k+j}\Lambda(v))|\\
&\leq& \delta \log N(\epsilon + L^{-\frac{p}{2}}) + \sum^{[\delta \log N]}_{j = 1}e^{-k-j} \leq \delta(1 + (\log N)^{1 - \frac{pq}{2}}) + e^{-k},
\eeqnas
where we take $\epsilon = (\log N)^{-1}$ and $L = (\log N)^{q}$ as in \eqref{cd.4},  such that $pq > 1$ for some $1 < p < d$. Then, there exists $0 < \delta < 1$ such that
\beqna\label{tight.1}
& &\bbP\big(\frac{1}{\sigma_{m, n}\sqrt{\log N}}\max_{i \leq \delta \log N}|E^{1}_{k+i} - E^{1}_{k}| \geq \xi \big) \leq \frac{2\delta}{\sigma_{m, n}\xi}(\log N)^{-\min\{\frac{pq-1}{2}, \frac{1}{2}\}} + \frac{1}{\sigma_{m, n}\xi}(\log N)^{-\frac{1}{2}}\leq \delta \eta,
\eeqna
where we take $(\log N)^{\min\{\frac{pq-1}{2}, \frac{1}{2}\}} > 4(\sigma_{m, n} \xi \eta)^{-1}$.
By \eqref{est.eps} and \eqref{asym.e}, we estimate that
\beqna\label{tight.2}
\nonumber & &\sup_{t \in [0, 1]}|E^{2}_{[t\log N]}| = \max_{i \leq \log N} \big(\sum^{i}_{j=0}\bbE \hat{f}^{L}_{\epsilon}(A_{j}\Lambda(v)) -C_{m, n}i \big) \\
\nonumber &\leq& \sum^{[\log N]-1}_{j=0}\big(\bbE | \hat{f}^{L}_{\epsilon}(A_{j}\Lambda(v)) -  \hat{\chi}(A_{j}\Lambda(v)) |+  |\bbE \hat{\chi}(A_{j}\Lambda(v)) - C_{m, n}| \big)\\
 &\leq& \epsilon \log N  +\sum^{[\log N]-1}_{j=0}e^{-j}+ O(1) = O(1),
\eeqna
i.e. $sup_{t \in [0, 1]}|E^{2}_{[t\log N]}|$ is bounded. 

Then it suffices to prove that \eqref{cond.ip.2} holds for $\{\frac{1}{\sigma_{m, n}\sqrt{\log N}}S^{\epsilon, L}_{[t\log N]}\}$. To do so, we make use of Theorem 12.2, \cite{billing} as follows.
\bthm\label{thm.mon}
Let $\xi_{1}, \dots, \xi_{n}, \dots$ be random variables and $T_{n} = \sum^{n}_{i=1}\xi_{i}$. If for any $0 \leq i \leq j \leq m$, there exist nonnegative numbers $u_{1}, \dots, u_{m}$ such that
\beqna\label{cond.mon}
\bbE(|T_{j} -T_{i}|^{\gamma}) \leq (\sum^{j}_{l = i+1}u_{l})^{\alpha}
\eeqna
for $\gamma > 0$ and $\alpha > 1$, then for all positive $\lambda$, 
\beqna
\bbP(\max_{0 < k \leq m}|T_{k}| \geq \lambda) \leq \frac{C_{\gamma, \alpha}}{\lambda^{\gamma}}(\sum^{m}_{l = 1}u_{l})^{\alpha}, 
\eeqna
where $C_{\gamma, \alpha}$ is a constant depending only on $\gamma$, $\alpha$.
\ethm

The following moment estimate shows that \eqref{cond.mon} holds for $\{\frac{1}{\sigma_{m, n}\sqrt{\log N}}S^{\epsilon, L}_{[\log N]}\}$ with $\gamma = 4$, $\alpha = 2$, $u_{l} = c_{0}$ for all $l$, where $c_{0}$ is a constant.
\blem
For $0 \leq  s_{1} \leq  s_{2} \leq [\log N]$, we have
\beqna\label{est.m4}
\bbE|S^{\epsilon, L}_{s_{2}} - S^{\epsilon, L}_{s_{1}}|^{4} \leq c^{2}_{0}|s_{2} - s_{1}|^{2}.
\eeqna
\elem

\bpf
Define
\beqnas
F^{\epsilon, L}_{N}(v) &=&  \frac{1}{\sigma_{m, n}\sqrt{\log N}}(S^{\epsilon, L}_{s_{2}} - S^{\epsilon, L}_{s_{1}}) \\
&=& \frac{1}{\sigma_{m, n}\sqrt{\log N}}\sum^{s_{2}-1}_{s=s_{1}}\big( \hat{f}^{L}_{\epsilon}(A_{s}\Lambda(v)) -  \int_{\cV}\hat{f}^{L}_{\epsilon}(A_{s}\Lambda(v))dm_{\cV}(v) \big).
\eeqnas
Recall the relationship between the 4th moment and the cumulants
\beqnas
\bbE |F^{\epsilon, L}_{N}|^{4} &=& \Cum^{(4)}(F^{\epsilon, L}_{N}) + 4\Cum^{(3)}(F^{\epsilon, L}_{N})\bbE(F^{\epsilon, L}_{N}) + 3\Var(F^{\epsilon, L}_{N})^{2} + 6\Var(F^{\epsilon, L}_{N})\bbE(F^{\epsilon, L}_{N}) + (\bbE F^{\epsilon, L}_{N})^{4}\\
&=& \Cum^{(4)}(F^{\epsilon, L}_{N}) + 3\Var(F^{\epsilon, L}_{N})^{2},
\eeqnas
such that
\beqnas
\bbE|S^{\epsilon, L}_{s_{2}} - S^{\epsilon, L}_{s_{1}}|^{4} =  \sigma^{4}_{m, n}(\log N)^{2}\bbE |F^{\epsilon, L}_{N}|^{4} = \sigma^{4}_{m, n}(\log N)^{2}(\Cum^{(4)}(F^{\epsilon, L}_{N}) + 3\Var(F^{\epsilon, L}_{N})^{2}).
\eeqnas
By the estimates~\eqref{est.var.k4} and \eqref{cum.r} and the choice of the parameters~\eqref{cd.4}, we deduce that
\beqnas
\log N\Var(F^{\epsilon, L}_{N}) &=&  (s_{2} - s_{1})(1 + O((\log \log N)(\log N)^{-a_{1}}),\\
(\log N)^{2}|\Cum^{(4)}(F^{\epsilon, L}_{N})| &=& (s_{2} - s_{1})^{2}O((\log \log N)^{3}(\log N)^{-a_{2}}),
\eeqnas
for some positive constants $a_{1}, a_{2}$. Note that the implicit constant in~\eqref{cum.r} only depends on $\supp(f_{\epsilon})$. Thus there exists a constant $c^{2}_{0}$, which is independent of $s_{1}$, $s_{2}$ and $N$, such that \eqref{est.m4} holds.  
\epf

According to Theorem~\ref{thm.mon},  the estimate  \eqref{est.m4} leads to 
\beqnas
\bbP(\max_{0 < i \leq \delta \log N}\frac{1}{\sigma_{m, n}\sqrt{\log N}}|S^{\epsilon, L}_{k+i} - S^{\epsilon, L}_{k}| \geq \xi) \leq \frac{C}{\xi^{4}}\delta^{2}
\eeqnas
for any $k > 0$ and $N$ large enough. By choosing $\delta < \frac{\eta \xi^{4}}{C}$, we prove that \eqref{cond.ip.2} holds for $\{\frac{1}{\sigma_{m, n}\sqrt{\log N}}S^{\epsilon, L}_{[t\log N]}\}$. Thus we complete the proof. 
\hfill \qedsymbol

\bibliographystyle{amsplain}

\begin{thebibliography}{10}

\bibitem{BBV13}
Dzmitry Badziahin, Victor Beresnevich, and Sanju Velani, \emph{Inhomogeneous
  theory of dual Diophantine approximation on manifolds}, Advances in
  Mathematics \textbf{232} (2013), no.~1, 1--35.

\bibitem{BV06}
Victor Beresnevich and Sanju Velani, \emph{Schmidt's theorem, Hausdorff
  measures, and slicing}, International Mathematics Research Notices (2006),
  48794.

\bibitem{billing}
Patrick Billingsley, \emph{Convergence of probability measures}, John Wiley and
  Sons, Inc., 1968.

\bibitem{BG19}
Michael Bj{\"o}rklund and Alexander Gorodnik, \emph{Central limit theorems for
  Diophantine approximants}, Mathematische Annalen \textbf{374} (2019), no.~3,
  1371--1437.

\bibitem{Borgne}
S.~Le Borgne, \emph{Principes d'invariance pour les flots diagonaux sur sl(d,
  r)/sl(d, z)}, Ann. Inst. H. Poincar{\'e} Probab. Stat. \textbf{38} (2002),
  no.~4, 581--612.

\bibitem{Cas1}
J.W.S. Cassels, \emph{Some metrical theorems in Diophantine approximation I},
  Proc. Camb. Phil. Soc \textbf{46} (1950), 209--218.

\bibitem{DFV}
D.I. Dolgopyat, B.~Fayad, and I.~Vinogradov, \emph{Central limit theorem for
  simultaneous diophantine approximations}, J. {\'E}c. Polytech. Math.
  \textbf{4} (2017), 1--36.

\bibitem{EMV}
M.~Einsiedler, G.~Margulis, and A.~Venkatesh, \emph{Effective equidistribution
  for closed orbits of semisimple groups on homogeneous spaces}, Inventiones
  mathematicae \textbf{177} (2009), no.~1, 137--212.

\bibitem{fuchs2}
M.~Fuchs, \emph{Invariance principles in metric Diophantine approximation},
  Monatshefte fur Mathematik \textbf{139} (2003), 177--203.

\bibitem{fuchs}
\bysame, \emph{On a problem of W. J. Leveque concerning metric Diophantine
  approximation}, Trans. Am. Math. Soc. \textbf{355} (2003), no.~5, 1787--1801.

\bibitem{K26}
A.~Khintchine, \emph{Zur metrischen th{\'e}orie der Diophantischen
  approximationen}, Math. Z. \textbf{24} (1926), 706--714.

\bibitem{Kim1}
W.~Kim, \emph{Effective equidistribution of expanding translates in the space
  of affine lattices}, Duke Math. J \textbf{173} (2024), no.~17, 3317--3375.

\bibitem{KM12}
D.~Y. Kleinbock and G.~A. Margulis, \emph{On effective equidistribution of
  expanding translates of certain orbits in the space of lattices}, Number
  Theory, Analysis and Geometry: In Memory of Serge Lang (2012), 385--396.

\bibitem{KM96}
Dmitry Kleinbock and G~Margulis, \emph{Bounded orbits of nonquasiunipotent
  flows on homogeneous spaces}, Amer. Math. Soc. Transl. \textbf{171} (1996).

\bibitem{Lev1}
W.J. Leveque, \emph{On the frequency of small fractional parts in certain real
  sequences I}, Trans. Am. Math. Soc. \textbf{87} (1958), 237--260.

\bibitem{Lev2}
\bysame, \emph{On the frequency of small fractional parts in certain real
  sequences II}, Trans. Am. Math. Soc. \textbf{94} (1959), 130--149.

\bibitem{Phil}
W.~Philipp, \emph{Mixing sequences of random variables and probabilistic number
  theory.}, In: Memoirs of the American Mathematical Society, vol. 114,
  American Mathematical Society, Providence, 1971.

\bibitem{schm1}
W.~Schmidt, \emph{A metrical theorem in Diophantine approximation}, Canad. J.
  Math. \textbf{12} (1960), 619--631.

\bibitem{sz58}
P.~Sz\"usz, \emph{\"Uber die metrische theorie der Diophantischen
  approximation}, Acta Math. Acad. Sci. Hungar. \textbf{9} (1958), no.~177-193.

\end{thebibliography}

\end{document}